\documentclass[10pt]{article}  
\usepackage{latexsym}  
\usepackage[all]{xy}  
\usepackage{amsmath}  
\usepackage{amssymb}  
\usepackage{amsfonts}  
  
\newtheorem{thm}{Theorem}[subsection]  
\newtheorem{prop}[thm]{Proposition}  
\newtheorem{lem}[thm]{Lemma}  
\newtheorem{df}[thm]{Definition}  
\newtheorem{cor}[thm]{Corollary}

\newtheorem{rmk}[thm]{Remark}  
\newtheorem{ex}[thm]{Example}

\begin{document}

\title{\textbf{Homotopical Algebraic Geometry I:\\Topos theory}}   
\bigskip  
\bigskip  
  
\author{\bigskip\\  
Bertrand To\"en \\
\small{Laboratoire Emile Picard}\\
\small{UMR CNRS 5580} \\
\small{Universit\'{e} Paul Sabatier, Toulouse}\\
\small{France}\\
\bigskip \and  
\bigskip \\  
Gabriele Vezzosi \\  
\small{Dipartimento di Matematica Applicata}\\
\small{``G. Sansone''}\\
\small{Universit\`a di Firenze}\\  
\small{Italy}\\  
\bigskip}  
  
\date{\small{November 2002.\\ Final version, to appear in Adv. Math., June 2004.}}


\maketitle

\begin{abstract}  
This is the first of a series of papers devoted to lay the foundations of Algebraic Geometry in homotopical   
and higher categorical contexts. In this first part we investigate a notion of \textit{higher topos}.  
  
For this, we use $S$-categories (i.e. simplicially enriched categories) as models for certain   
kind of $\infty$-categories, and we develop the notions of \textit{$S$-topologies, $S$-sites} and   
\textit{stacks} over them. We prove in particular, that for an $S$-category $T$ endowed with an   
$S$-topology, there exists a model category of stacks over $T$,   
generalizing the model category structure on simplicial presheaves   
over a Grothendieck site of A. Joyal and R. Jardine.  
We also prove some analogs of the relations between topologies and localizing subcategories of the  
categories of presheaves, by proving that there exists a one-to-one correspodence between   
$S$-topologies on an $S$-category $T$, and certain \textit{left exact Bousfield localizations}  
of the model category of pre-stacks on $T$.   
Based on the above results, we study the notion of \textit{model topos} introduced by C.  Rezk, and   
we relate it to our model categories of stacks over $S$-sites.  
  
In the second part of the paper, we present a parallel theory where $S$-categories, $S$-topologies and  
$S$-sites are replaced by \textit{model categories}, \textit{model topologies} and \textit{model sites}. We prove that  
a canonical way to pass from the theory of stacks over model sites to the theory   
of stacks over $S$-sites is provided by the  
simplicial localization construction of Dwyer and Kan. As an example of application, we propose a   
definition of \textit{\'etale $K$-theory of ring spectra}, extending the \'etale $K$-theory of  
commutative rings.  
  
\end{abstract}  
  
  
\medskip  
  
\textsf{MSC-class:} $14$A$20$; $18$G$55$; $55$P$43$; $55$U$40$; $18$F$10$.  
  
  
\tableofcontents  

\bigskip
\bigskip  
  
\setcounter{section}{0}  
  
\begin{section}{Introduction}  
  
This is the first part of a series of papers devoted to the foundations of Algebraic Geometry in homotopical   
and higher categorical contexts, the ultimate goal being a theory of \textit{algebraic  
geometry over monoidal $\infty$-categories}, a higher categorical generalization of     
\textit{algebraic geometry over monoidal categories} (as developed, for example, in \cite{del1,del2,hak}). We refer the  
reader to the Introduction of the research announcement \cite{web} and to \cite{camb}, where motivations and prospective   
applications (mainly to the so-called \textit{derived moduli spaces} of \cite{ko1,ck1,ck2}) are provided.   
These applications, together with the remaining required \textit{monoidal} part of the theory, will be given in \cite{partII}.  
  
In the present work we investigate the required theory of \textit{higher sheaves}, or equivalently \textit{stacks}, 
as well as its associated   
notion of \textit{higher topoi}. \\  
  
\begin{center} \textbf{Topologies, sheaves and topoi} \end{center}   
  
As we will proceed by analogy,   
we will start by recalling some basic   
constructions and results from topos theory, in a way that is 
suited for our generalization. Our references for this overview are \cite{sga4,schu,mm}.  
Throughout this introduction we will neglect any kind of set   
theoretical issues, always assuming that   
categories are \textit{small} when required.  
  
Let us start with a category $C$ and let us denote by $Pr(C)$ the category of  
presheaves of sets on $C$ (i.e. $Pr(C):=Set^{C^{op}}$).   
If $C$ is endowed with a Grothendieck topology $\tau$, one can define the notion of   
$\tau$-\textit{local isomorphisms in $Pr(C)$} by requiring   
injectivity and surjectivity only up to a $\tau$-covering.
We denote by $\Sigma_{\tau}$ the subcategory of $Pr(C)$ consisting of local isomorphisms. One possible way to define  
the category $Sh_{\tau}(C)$, of sheaves (of sets) on the Grothendieck site $(C,\tau)$, is by setting  
$$Sh_{\tau}(C):=\Sigma_{\tau}^{-1}Pr(C),$$  
where $\Sigma_{\tau}^{-1}Pr(C)$ denotes the \textit{localization} of $Pr(C)$ along $\Sigma_{\tau}$ i.e. the category obtained from $Pr(C)$ by formally inverting the   
morphisms in $\Sigma_{\tau}$ (see \cite[19.1, 20.3.6 (a)]{schu}). The main basic properties of the category   
$Sh_{\tau}(C)$ are collected in the following well known theorem.  
  
\begin{thm}\label{thmA}  
Let $(C,\tau)$ be a Grothendieck site and $Sh_{\tau}(C)$ its category of sheaves as defined above.  
\begin{enumerate}  
\item The category $Sh_{\tau}(C)$ has all limits and colimits.  
\item The natural localization morphism $a : Pr(C) \longrightarrow Sh_{\tau}(C)$ is left   
exact (i.e. commutes with finite limits) and  
has a fully faithful right adjoint $j : Sh_{\tau}(C) \longrightarrow Pr(C)$.  
\item The category $Sh_{\tau}(C)$ is cartesian closed (i.e. has internal $Hom$-objects).  
\end{enumerate}  
\end{thm}  
  
Of course, the essential image of the functor $j : Sh_{\tau}(C) \longrightarrow Pr(C)$ is the   
usual subcategory of sheaves, i.e. of presheaves having descent with respect   
to $\tau$-coverings, and the localization functor $a$ becomes   
equivalent to the associated sheaf functor. The definition   
of $Sh_{\tau}(C)$ as $\Sigma_{\tau}^{-1}Pr(C)$ is therefore   
a way to define the category of sheaves without even mentioning what   
a sheaf is precisely. 
  
In particular, Theorem \ref{thmA} shows that the datum of a topology $\tau$ on $C$ gives rise to an   
adjunction  
$$a : Pr(C) \longrightarrow Sh_{\tau}(C) \qquad Pr(C) \longleftarrow Sh_{\tau}(C) : j,$$  
whith $j$ fully faithful and $a$ left exact. Such an adjoint pair will be called an \textit{exact localization} of  
the category $Pr(C)$. Another fundamental result in sheaf theory is the following  
  
\begin{thm}\label{thmB}  
The rule sending a Grothendieck topology $\tau$ on $C$ to the exact localization   
$$a : Pr(C) \longrightarrow Sh_{\tau}(C) \qquad Pr(C) \longleftarrow Sh_{\tau}(C) : j,$$  
defines a bijective correspondence between the set of topologies on $C$ and the set of (equivalences classes) of  
exact localizations of the category $Pr(C)$.  
In particular, for a category $T$ the following two conditions are equivalent  
\begin{itemize}  
\item There exists a category $C$ and a Grothendieck topology $\tau$ on $C$ such that $T$ is equivalent to   
$Sh_{\tau}(C)$.  
\item There exists a category $C$ and a left exact localization  
$$a : Pr(C) \longrightarrow T \qquad Pr(C) \longleftarrow T : j.$$  
\end{itemize}  
A category satisfying one the previous conditions is called a \emph{Grothendieck topos}.  
\end{thm}

Finally, a famous theorem by Giraud (\cite{sga4} Exp. IV, Th\'eor\`eme 1.2) provides an internal 
characterization of Grothendieck topoi.

\begin{thm}\label{thmCgiraud}
\emph{(Giraud's Theorem)} A category $T$ is a Grothendieck topos if and only if
it satisfies the following conditions.
\begin{enumerate}
\item The category $T$ is has a small set 
of strong generators.
\item The category $T$ has small colimits.
\item Sums are disjoint in $T$
(i.e. $x_{j}\times_{\coprod_{i} x_{i}}x_{k}\simeq \emptyset$
for all $j\neq k$).
\item Colimits commute with pull backs.
\item Any equivalence relation
is effective. 
\end{enumerate}
\end{thm}

The main results of this work are generalizations to a homotopical setting  
of the notions of topologies, sites and sheaves satisfying analogs of theorems \ref{thmA}, \ref{thmB}, and \ref{thmCgiraud}. We have chosen to use both the concept of $S$\textit{-categories} (i.e. simplicially enriched categories)
and of \textit{model categories} as our versions of base categories carrying homotopical data. For both we have
developed homotopical notions of \textit{topologies, sites} and \textit{sheaves}, and proved analogs of theorems
\ref{thmA}, \ref{thmB} and \ref{thmCgiraud} which we will now describe in  
more details. \\  
  
\medskip  
  
\begin{center} \textbf{$S$-topologies, $S$-sites and stacks} \end{center}  

Let $T$ be a base $S$-category. We consider the category $SPr(T)$, of $T^{op}$-diagrams in the category $SSet$ of  
simplicial sets. This category can be endowed with an objectwise model structure for which the equivalences are  
defined objectwise on $T$. This model category $SPr(T)$ will be called the \textit{model   
category of pre-stacks on $T$}, and  
will be our higher analog of the category of presheaves of sets. The category $SPr(T)$ comes with a   
natural \textit{Yoneda  
embedding} $L\underline{h} : T \longrightarrow SPr(T)$, a \textit{up to homotopy analog} of the   
usual embedding of a category into the   
category of presheaves on it (see Corollary \ref{c1}).  
  
We now consider $\mathrm{Ho}(T)$, the category having the same objects as $T$ but for which  
the sets of morphisms are the connected components of the simplicial sets of morphisms in $T$.   
Though it might be surprising at first sight, we define an   
\textit{$S$-topology on the $S$-category $T$} to be simply a   
Grothendieck topology on the category $\mathrm{Ho}(T)$   
(see Defintion \ref{d9}). A pair $(T,\tau)$, where $T$ is an   
$S$-category and $\tau$ is an $S$-topology on $T$, will be  
called an \textit{$S$-site}.   
Of course, when $T$ is a usual category (i.e. all its simplicial sets of morphisms are discrete), an 
$S$-topology on $T$ is nothing else than a Grothendieck topology on $T$. Therefore,   
a site is in particular an $S$-site, and   
our definitions are actual generalizations of the usual definitions of topologies and  
sites.  
  
For the category of presheaves of sets on a Grothendieck site, we have already   
mentioned that the topology induces a notion of local isomorphisms. In the case where $(T,\tau)$ is an 
$S$-site we define a notion of   
\textit{local equivalences} in $SPr(T)$ (see Definition \ref{d14}).  
When $T$ is a category, and therefore $(T,\tau)$ is a site in the usual sense,   
our notion of local equivalences specializes to the  
notion introduced by L. Illusie and later by R. Jardine (see \cite{ja}).   
Our first main theorem is a generalization of the existence   
of the local model category structure on the category of simplicial presheaves on a site  (see \cite{ja,bl}).  
  
\begin{thm}\emph{(Thm. \ref{t3}, Prop. \ref{p3} and Cor. \ref{cc2})}\label{thmA'}  
Let $(T,\tau)$ be an $S$-site.  
\begin{enumerate}  
\item  There exists a model structure on the category $SPr(T)$, called the \emph{local model structure},    
for which the equivalences are the local equivalences. This new model category, denoted by $SPr_{\tau}(T)$,   
is furthermore the left Bousfield localization of the model category $SPr(T)$ of pre-stacks along the local equivalences.   
\item The identity functor  
$$\mathrm{Id} : SPr(T) \longrightarrow SPr_{\tau}(T)$$  
commutes with homotopy fibered products.   
\item The homotopy category $\mathrm{Ho}(SPr_{\tau}(T))$ is cartesian closed, or equivalently, it has internal $Hom$-objects.  
\end{enumerate}  
The model category $SPr_{\tau}(T)$ is called the \emph{model category of stacks on the   
$S$-site $(T,\tau)$}.  
\end{thm}  
  
This theorem is our \textit{higher analog} of Theorem \ref{thmA}. Indeed, the existence of the local model structure  
formally implies the existence of \textit{homotopy limits} and \textit{homotopy   
colimits} in $SPr_{\tau}(T)$, which are homotopical generalizations  
of the notion of limits and colimits (see \cite[\S $19$]{hi}).   
Moreover, $SPr_{\tau}(T)$ being a left Bousfield localization of  
$SPr(T)$, the identity functor $\mathrm{Id} : SPr_{\tau}(T) \longrightarrow SPr(T)$ is a right Quillen functor and therefore  
induces an adjunction on the level of homotopy categories  
$$a:=\mathbb{L}\mathrm{Id} : \mathrm{Ho}(SPr(T)) \longrightarrow \mathrm{Ho}(SPr_{\tau}(T)) \qquad  
\mathrm{Ho}(SPr(T)) \longleftarrow \mathrm{Ho}(SPr_{\tau}(T)) : j:=\mathbb{R}\mathrm{Id}.$$  
It is a general property of Bousfield localizations that the functor $j$ is \textit{fully faithful}, and  
Theorem \ref{thmA'} $(2)$ implies that the functor $a$ is  
\textit{homotopically left exact}, i.e. commutes   
with homotopy fibered products. Finally,     
part $(3)$ of Theorem \ref{thmA'} is a homotopical analog of Theorem \ref{thmA} $(3)$.   
  
As in the case of sheaves on a site, it remains to characterize the essential image of the inclusion functor  
$j :  \mathrm{Ho}(SPr_{\tau}(T)) \longrightarrow \mathrm{Ho}(SPr(T))$.    
One possible homotopy analog of the sheaf condition is the \textit{hyperdescent property} for objects in $SPr(T)$ (see  
Definition \ref{d15}). It is a corollary of our   
proof of the existence of the local model structure $SPr_{\tau}(T)$ that   
the essential image of the inclusion functor $j : \mathrm{Ho}(SPr_{\tau}(T)) \longrightarrow \mathrm{Ho}(SPr(T))$ is exactly the   
full subcategory of objects satisfying the hyperdescent condition (see Corollary \ref{c6}).   
We call these objects  
\textit{stacks} over the $S$-site $(T,\tau)$ (Definition \ref{d16}).   
The functor $a : \mathrm{Ho}(SPr(T)) \longrightarrow \mathrm{Ho}(SPr_{\tau}(T))$  
can then be identified with the \textit{associated stack functor} (Definition \ref{d16}). 

Finally, we would like to mention that the model categories $SPr_{\tau}(T)$ are not 
in general Quillen equivalent to model categories of simplicial presheaves on some site. 
Therefore, Theorem \ref{thmA'} is a new result in the sense that neither its statement nor its proof can be reduced
to previously known notions and results in the theory of simplicial presheaves. \\
  
\medskip  
  
\begin{center} \textbf{Model topoi and $S$-topoi} \end{center}  
  
Based on the previously described notions of $S$-sites and stacks, we develop a 
related theory of \textit{topoi}. For this,   
note that Theorem \ref{thmA'} implies that   
an $S$-topology $\tau$ on an $S$-category $T$ gives rise to the model category  
$SPr_{\tau}(T)$, which is a left Bousfield localization of the model category $SPr(T)$. This Bousfield  
localization has moreover the property that the identity functor $\mathrm{Id} : SPr(T) \longrightarrow SPr_{\tau}(T)$  
preserves homotopy fibered products. We call such a localization a \textit{left exact Bousfield localization}  
of $SPr(T)$ (see Definition \ref{d18}). This notion is a homotopical analog of the notion of exact localization  
appearing in topos theory as reviewed before Theorem \ref{thmB}.   
The rule $\tau \mapsto SPr_{\tau}(T)$, defines a map    
from the set of $S$-topologies on a given $S$-category $T$ to the set of left exact Bousfield localizations of the model   
category $SPr(T)$.  
The model category $SPr_{\tau}(T)$ also possesses a natural \textit{additional} property, called \textit{$t$-completeness} which is a new feature of the homotopical context which does not have any counterpart in classical sheaf theory (see Definition \ref{d19}).   
An object $x$ in some model category $M$ is called $n$-\textit{truncated} if for  
any $y\in M$, the mapping space $Map_{M}(y,x)$ is an $n$-truncated simplicial set; an object in $M$ is \textit{truncated} if it is $n$-truncated for some $n\geq 0$. A model category $M$ will then be called $t$\emph{-complete} if truncated objects detect   
isomorphisms in $\mathrm{Ho}(M)$: a morphism $u : a \rightarrow b$  
in $\mathrm{Ho}(M)$ is an isomorphism if and only if, for any truncated object $x$ in $\mathrm{Ho}(M)$, the map   
$u^{*} : [b,x] \longrightarrow [a,x]$ is bijective.

The notion of $t$-completeness is very natural and very often satisfied
as most of the equivalences in model categories are defined
using isomorphisms on certain homotopy groups. The $t$-completeness
assumption simplyy states that an object with 
trivial homotopy groups is homotopically trivial, which is
a very natural and intuitive condition. The usefulness of this notion of $t$-completeness is explained by the following theorem, which is our analog of Theorem \ref{thmB}  
  
\begin{thm}\emph{(Thm. \ref{t4} and Cor. \ref{c7})}\label{thmB'}  
Let $T$ be an $S$-category. The correspondence $\tau \mapsto SPr_{\tau}(T)$ induces a bijection between $S$-topologies  
on $T$ and $t$-complete left exact Bousfield localizations of $SPr(T)$.  
In particular, for a model category $M$ the following two conditions are equivalent  
\begin{itemize}  
\item There exists an $S$-category $T$ and an $S$-topology on $T$ such that $M$ is Quillen equivalent to   
$SPr_{\tau}(T)$.  
\item The model category $M$ is $t$-complete and there exists an $S$-category $T$ such that $M$ is Quillen equivalent  
to a left exact Bousfield localization of $SPr(T)$.   
\end{itemize}  
A model category satisfying one the previous conditions is called a   
\emph{$t$-complete model topos}.  
\end{thm}  

It is important to stress that there are $t$-complete model topoi which are \textit{not} Quillen   
equivalent to any $SPr_{\tau}(C)$, for $C$ a usual category (see Remark \ref{lurie} (1)). Therefore, Theorem \ref{thmB'} also shows the unavoidable relevance of considering topologies on general $S$-categories rather than only on usual categories. In other words there is no way to reduce the theory developed in this paper to the theory of simplicial presheaves over Grothendieck sites
as done in \cite{ja,jo}. 

The above notion of \textit{model topos} was suggested to us by C. Rezk, who   
defined a more general notion of \textit{homotopy topos} (a model topos without the $t$-completness assumption), which is  
a model category Quillen equivalent to an arbitrary left exact Bousfield localization of some $SPr(T)$ 
(see Definition \ref{d18}). 
The relevance of Theorem \ref{thmB'} is that, on one hand it 
shows that the notion of $S$-topology we used is 
correct exactly because it classifies all 
($t$-complete) left exact Bousfield localizations, and, 
on the other hand it provides an answer to 
a question raised by Rezk on which notion of 
topology could be the source of his homotopy topoi. 

It is known that there exist model topoi which are not $t$-complete (see Remark \ref{lurie}), and therefore our notion of  
stacks over $S$-categories does not model \textit{all} of Rezk's homotopy topoi.   
However, we are strongly convinced that   
Theorem \ref{thmB'} has a more general version, in which the $t$-completeness assumption is dropped, involving a corresponding notion  
of \textit{hyper-topology} on $S$-categories as well as the associated notion of \textit{hyper-stack}   
(see Remark \ref{lurie}).   
  
Using the above notion of model topos, we also define the notion of $S$\textit{-topos}. 
An $S$-topos is by definition an $S$-category which is   
equivalent, as an $S$-category, to some $LM$, for $M$ a model topos (see Definition \ref{dstop}).   
Here we have denoted by $LM$ the Dwyer-Kan simplicial localization of $M$ with respect to   
the set of its weak equivalences (see the next paragraph for further explanations on the   
Dwyer-Kan localization). \\
  
\medskip  
  
\begin{center} \textbf{$S$-Categories and model categories} \end{center}  
  
Most of the $S$-categories one encounters in practice come from model categories via   
the Dwyer-Kan \textit{simplicial localization}. The simplicial localization is a refined version of the  
Gabriel-Zisman localization of categories. It associates an $S$-category $L(C,S)$ to any   
category $C$ equipped with a subcategory $S\subset C$ (see $\S 2.2$), such that the homotopy category  
$\mathrm{Ho}(L(C,S))$ is naturally equivalent to the Gabriel-Zisman localization $S^{-1}C$, but in general  
$L(C,S)$ contains non-trivial higher homotopical informations.   
The simplicial localization   
construction is particularly well behaved when applied to a model category $M$ equipped with its  
subcategory of weak equivalences $W \subset M$: in fact in this case the $S$-category $LM:=L(M,W)$ encodes the   
so-called \textit{homotopy mapping spaces} of the model category $M$ (see $\S  2.2$). We will show furthemore that  
the notions of $S$-topologies, $S$-sites and stacks previously described in this introduction, also   
have their analogs in the model category context, and that the simplicial localization construction  
allows one to pass from the theory over model categories to the theory over $S$-categories.  
  
For a model\footnote{Actually, in Section $4$, all the constructions are given for the weaker notion of   
\textit{pseudo-model categories} because we will need this increased   
flexibility in some present and future applications. However,   
the case of model categories will be enough for this introduction.} category $M$,   
we consider the category $SPr(M)$ of simplicial presheaves on $M$, together  
with its objectwise model structure. We define the model category $M^{\wedge}$ to be the   
left Bousfield localization of $SPr(M)$ along the set of equivalences in $M$   
(see Definition \ref{d20}). In particular, unlike that of $SPr(M)$, the   
model structure of $M^{\wedge}$ takes into account the fact that $M$ is   
not just a bare category but has an additional (model) structure.   
The model category $M^{\wedge}$ is called the \textit{model category of pre-stacks on $M$}, and it is  
important to remark that its homotopy category can be identified with the full subcategory  
of $\mathrm{Ho}(SPr(M))$ consisting of functors $F : M^{op} \longrightarrow SSet$ sending   
equivalences in $M$ to equivalences of simplicial sets. We construct a homotopical Yoneda-like functor   
$$\underline{h} : M \longrightarrow M^{\wedge},$$  
roughly speaking by sending an object $x$ to the simplicial presheaf $y \mapsto Map_{M}(y,x)$, where  
$Map_{M}(-,-)$ denotes the homotopy mapping space in the model   
category $M$ (see Definition \ref{d21'}). An easy but fundamental  
result states that the functor $\underline{h}$ possesses a right derived functor  
$$\mathbb{R}\underline{h} : \mathrm{Ho}(M) \longrightarrow \mathrm{Ho}(M^{\wedge})$$  
which is fully faithful (Theorem \ref{t10}). This is a model category version of the Yoneda lemma.   
  
We also define the notion of a \textit{model pre-topology} on the model category $M$ and show that this induces in a natural way a Grothendieck topology on the homotopy category $\mathrm{Ho}(M)$.   
A model category endowed with a model pre-topology will be called a \textit{model site} (see Definition  
\ref{d21}). For a model site $(M,\tau)$, we define a notion of   
\textit{local equivalences} in the category of pre-stacks $M^{\wedge}$.  
The analog of Theorem \ref{thmA} for model categories is then the following   
  
\begin{thm}\emph{(Thm. \ref{t5})}\label{thmA''}  
Let $(M,\tau)$ be a model site.  
\begin{enumerate}  
\item 
There exists a model structure on the category $M^{\wedge}$,   
called the \emph{local model structure},    
for which the equivalences are the local equivalences.   
This new model category, denoted by $M^{\sim,\tau}$,   
is furthermore the left Bousfield localization of the model category of pre-stacks   
$M^{\wedge}$ along the local equivalences.   
\item The identity functor  
$$\mathrm{Id} : M^{\wedge} \longrightarrow M^{\sim,\tau}$$  
commutes with homotopy fibred products.   
\item The homotopy category $\mathrm{Ho}(M^{\sim,\tau})$ is cartesian closed.  
\end{enumerate}  
The model category $M^{\sim,\tau}$ is called the \emph{model category of stacks} on the   
model site $(M,\tau)$.  
\end{thm}  
  
As for stacks over $S$-sites, there exists a notion of object satisfying    
a \textit{hyperdescent condition} with respect to the topology   
$\tau$, and we prove that $\mathrm{Ho}(M^{\sim,\tau})$ can be identified with the full   
subcategory of $\mathrm{Ho}(M^{\wedge})$   
consisting of objects satisfying hyperdescent (see Definition \ref{d152}).  
  
Finally, we compare the two parallel constructions of stacks over $S$-sites and over model sites.  
  
\begin{thm}\emph{(Thm. \ref{t6})}\label{thmC}  
Let $(M,\tau)$ be a model site.   
\begin{enumerate}  
\item The simplicial localization $LM$ possesses an induced $S$-topology $\tau$, and   
is naturally an $S$-site.  
\item The two corresponding model categories of stacks $M^{\sim,\tau}$ and $SPr_{\tau}(LM)$   
are naturally Quillen equivalent. In particular $M^{\sim,\tau}$ is a $t$-complete  
model topos.  
\end{enumerate}  
\end{thm}    
  
The previous comparison theorem finds his pertinence in the fact that   
the two approaches, stacks over model sites and stacks over $S$-sites, seem to possess  
their own advantages and disadvantages, depending of the situation and the goal that one wants to reach.  
On a computational level  
the theory of stacks over model sites seems to be better suited than that of stacks over $S$-sites. On the other hand,  
$S$-categories and $S$-sites are much more intrinsic than model categories and model sites, and   
this has already some consequences, e.g. at the level of functoriality properties of the categories of stacks. We are  
convinced that having the full picture, including the two approaches and the comparison theorem \ref{thmC},   
will be a very friendly setting for the purpose of several future applications. \\  
  
\medskip  
  
\begin{center} \textbf{A Giraud theorem for model topoi} \end{center} 

Our version of Theorem \ref{thmCgiraud} is on the model categories' side of the
theory. The corresponding statement for $S$-categories
would drive us too far away from
the techniques used in this work, and will not be investigated here.

\begin{thm}\emph{(Thm. \ref{tgiraud})}\label{tgiraud'}
A combinatorial model category $M$ is a model topos if and only if it 
satisfies the following conditions.
\begin{enumerate}
\item Homotopy coproducts are disjoints in $M$.
\item Homotopy colimits are stable under homotopy pullbacks.
\item All Segal equivalences relations are homotopy effective.
\end{enumerate}
\end{thm}

The condition of being a combinatorial model
category is a set theoretic condition on $M$
(very often satisfied in practice), 
very similar to the condition of having
a small set of generators (see appendix $A.2$). Conditions $(1)$ and 
$(2)$ are straightforward homotopy theoretic analogs of conditions
$(3)$ and $(4)$ of Theorem \ref{thmCgiraud}: we essentially replace pushouths, pullbacks and colimits by 
homotopy pushouts, homotopy pullbacks and homotopy colimits (see Definition \ref{dgiraud}). Finally, 
condition $(3)$ of Theorem \ref{tgiraud'}, spelled out in Definition \ref{dgiraud}(3) and \ref{dgiraud}(4), is  a homotopical version of 
condition $(5)$ of Giraud's theorem \ref{thmCgiraud}, where groupoids of equivalence relations are replaced
by Segal groupoids and effectivity has to be 
 understood homotopically. 

The most important consequence of Theorem \ref{tgiraud'}
is the following complete characterization 
of $t$-complete model topoi.

\begin{cor}\emph{(Cor. \ref{cgiraud})}\label{cC'}
For a combinatorial model category $M$, the following two conditions
are equivalent.
\begin{enumerate}
\item There exists a small $S$-site $(T,\tau)$, such that
$M$ is Quillen equivalent to $SPr_{\tau}(T)$.
\item $M$ is $t$-complete and satisfies the conditions of 
Theorem \ref{tgiraud'}.
\end{enumerate}
\end{cor}

\medskip

\begin{center} \textbf{A topological application:   
\'etale $K$-theory of commutative $\mathbb{S}$-algebras} \end{center}  
  
As an example of application of our constructions, we give a definition of the  
\textit{\'etale $K$-theory of (commutative) $\mathbb{S}$-algebras}, which is to  
algebraic $K$-theory of $\mathbb{S}$-algebras (as defined for example in \cite[\S VI]{ekmm})   
what \'etale $K$-theory of rings is to algebraic $K$-theory of rings.
For this, we use the notion of   
etale morphisms of $\mathbb{S}$-algebras introduced in \cite{min} (and in \cite{web})   
in order to define an \textit{\'etale} pre-topology  
on the model category of commutative $\mathbb{S}$-algebras (see Definition \ref{etalecoverings}). Associated to this   
model pre-topology, we have the model category of \'etale stacks $(Aff_{\mathbb{S}})^{\sim,\text{\'et}}$;  
the functor $K$ that maps an $\mathbb{S}$-algebra $A$ to its algebraic  
$K$-theory space $K(A)$, defines an object $K \in (Aff_{\mathbb{S}})^{\sim,\text{\'et}}$.  
If $K_{\text{\'et}} \in (Aff_{\mathbb{S}})^{\sim,\text{\'et}}$ is an \'etale fibrant model for $K$, we define  
the space of \'etale $K$-theory of an $\mathbb{S}$-algebra $A$ to be   
the simplicial set $K_{\text{\'et}}(A)$ (see Definition \ref{dket}).   
Our general formalism also allows us to compare $K_{\text{\'et}}(Hk)$ with  
the usual definition of etale $K$-theory of a field $k$ (see Corollary \ref{restrictedcomparison}). 

This definition of \'etale $K$-theory of $\mathbb{S}$-algebras gives a possible answer to a question 
raised by J. Rognes  in \cite{ro}. In the future it might be used as a starting point to develop
\textit{\'etale localization techniques} in $K$-theory of $\mathbb{S}$-algebras, as Thomason's style \'etale descent 
theorem, analog of the Quillen-Lichtenbaum's conjecture, etc. For further applications of the general theory developed in this paper to algebraic geometry over commutative ring spectra, we refer the reader to \cite{partII} and \cite{newton}.\\
 
\medskip    
  
\begin{center} \textbf{Organization of the paper} \end{center}  
  
The paper is organized in five sections and one appendix. In \textsf{Section $2$} we review the   
main definitions and results concerning $S$-categories. Most of the materials can be found in  
the original papers \cite{dk1,dk2,dkh}, with the possible exception of the last two subsections.   
In \textsf{Section $3$} we define the notion of $S$-topologies, $S$-sites, local equivalences and   
stacks over $S$-sites. This section contains the proofs of theorems \ref{thmA'} and \ref{thmB'}.   
We prove in particular the existence of the local model structure  
as well as internal $Hom$'s (or equivalently, stacks of morphisms). We also investigate here the  
relations between Rezk's model topoi and $S$-topologies. \textsf{Section $4$}  
is devoted to the theory of model topologies, model sites and stacks over them. As it follows a pattern very  
similar to the one followed in Section 3 (for $S$-categories), some details have been  omitted. 
It also contains comparison results between the theory of stacks over $S$-sites and the theory of stacks over model sites, as well as the Giraud's style
theorem for model topoi.  
In \textsf{Section $5$} we present one application of the theory to the notion of   
\textit{\'etale $K$-theory of $\mathbb{S}$-algebras}. For this we review briefly the homotopy theory of   
$\mathbb{S}$-modules and $\mathbb{S}$-algebras, and we define an \'etale topology on the model category  
of commutative $\mathbb{S}$-algebra, which is an extension of the \'etale topology on affine schemes.  
Finally we use our general formalism to define the \'etale $K$-theory space of   
a commutative $\mathbb{S}$-algebra.  
  
Finally, in \textsf{Appendix A} we collected some definitions  
and conventions concerning model categories and the use of universes  
in this context. \\  
  
\medskip  
  
\begin{center} \textbf{Related works} \end{center}  
  
There has been several recent works on (higher) stacks theory which use a simplicial and/or   
a model categorical approach (see \cite{dhi,hol,ja2,s2,sh,to2,to3}).   
The present work is strongly  
based on the same idea that simplicial presheaves are after all very good models  
for \textit{stacks in $\infty$-groupoids}, and provide a powerful and rich theory.   
It may also be considered as a natural continuation of the foundational papers \cite{jo,ja}.  
  
A notion of a topology on a $2$-category, as well as a notion of \textit{stack over a  
$2$-site} has already been considered by R. Street in \cite{str}, D. Bourn in \cite{bou} 
and, more recently, by K. Behrend in his work on  
DG-schemes \cite{be}. Using truncation functors (Section \ref{truncations}), a precise 
comparison with these approaches will appear in   
the second part of this work \cite{partII} (the reader is also referred   
to Remark \ref{nella}).  
  
We have already mentioned that the notion of model topos used in Section $3.8$ essentially goes back to the 
unpublished manuscript \cite{rez2}, though it was
originally defined as left exact Bousfield localizations of model 
category of simplicial presheaves on some
usual category, which is not enough as we have seen. 
A different, but similar, version of our Giraud's theorem
\ref{tgiraud} appeared in \cite{rez2} as conjecture. 
The notion of $S$-topos introduced in Section $3.8$ seems 
new, though more or less equivalent to the notion of model topos. 
However, we think that both theories of $S$-categories and  
of model categories reach here their limits, as it seems quite difficult 
to define a reasonable notion of geometric morphisms between model topoi or between $S$-topoi. 
This problem can be solved by using Segal categories of \cite{sh,p} in order to introduce
a notion of \textit{Segal topos} as explained in \cite{msri}. 

A notion relatively closed to the notion of
Segal topos 
can also be found in \cite{s3} where \textit{Segal pre-topoi} are  
investigated and the question of the existence of a theory of Segal topoi is clearly addressed. 

Also closely related to our approach to model topoi is the
notion of \textit{$\infty$-topos} appeared in the recent preprint \cite{lu} by J. Lurie. The results of \cite{lu} are exposed in a rather different
context, and are essentially disjoints from ours. For example
the notion of topology is not considered in \cite{lu}
and results of the type \ref{t4}, \ref{c7} or \ref{cgiraud}
do not appear in it. Also, the notion of stack used by J. Lurie
is slightly different from ours (however the  
differences are quite subtle). An exception is Giraud's theorem
which first appeared in \cite{lu} in the context of
$\infty$-categories, and only later on in the last 
version of this work (February 2004) for model categories. 
These two works have been done independently,
though we must mention that the first version of the present paper has been
publicly available since July 2002 (an important
part of it was announced in \cite{web} which appeared
on the web during October 2001), whereas \cite{lu} appeared
in June 2003. 

Let us also mention that  
A. Joyal (see \cite{jo2}) has developed a theory of \textit{quasi-categories}, which is expected to be   
equivalent to the theory of $S$-categories and of Segal categories, and for which he has  
defined a notion of \textit{quasi-topos} very similar to   
the notion of Segal topos in \cite{msri}. The two approaches are expected to be  
equivalent.  
Also, the recent work of D-C. Cisinski (\cite{cis}) seems to be 
closely related to a notion of \textit{hypertopology} we discuss in Remark \ref{lurie} (3).  
  
Our definition of the \'etale topology for $\mathbb{S}$-algebras was strongly influenced  
by the content of \cite{min,min2}, and the definition of \'etale $K$-theory in the context of 
$\mathbb{S}$-algebras given in $\S 5$ was motivated by the note \cite{ro}.\\  
  
\medskip  
  

\noindent \textbf{Acknowledgments.} First, we would like to thank very warmly Markus Spitzweck for a very  
exciting discussion we had with him in Toulouse July $2000$, which  
turned out to be the starting point of our work. We wish  
especially to thank Carlos Simpson for precious conversations and  
friendly encouragement: the debt we owe to his huge amount of work on higher  
categories and higher stacks will be clear throughout this work.  
We are very thankful to Yuri Manin for his interest   
and in particular for his letter \cite{ma}. We thank Charles Rezk for a   
stimulating e-mail correspondence and for sharing with us   
his notion of \textit{model topos} (see Definition \ref{d18}).  
  
For many comments and discussions, we also thank Kai  
Behrend, Jan Gorsky, Vladimir Hinich, Rick Jardine, Andr\'e Joyal, Mikhail Kapranov, Ludmil Katzarkov,   
Maxim Kontsevich, Andrey Lazarev, Jacob Lurie, Michael Mandell, Peter May, Vahagn   
Minasian, Tony Pantev and John Rognes. It was Paul-Arne Ostv\ae r who pointed out to us the possible relevance of   
defining \'{e}tale $K$-theory of commutative ring spectra. We thank Stefan Schwede for pointing out to us the argument that led to the proof of Proposition \ref{stefan}. We are grateful to the referee for very useful suggestions concerning our exposition.
  
We are thankful to MSRI for support and for providing excellent working conditions   
during the Program \textit{Stacks, Intersection Theory and Nonabelian Hodge Theory}, January-May $2002$.  
The second author wishes to thank the  
Max Planck Institut f\"{u}r Mathematik in Bonn and the Laboratoire  
J. A. Dieudonn\'e of the Universit\'e de Nice Sophia-Antipolis for providing a  
particularly stimulating atmosphere during his visits when part of  
this work was conceived, written and tested in a seminar.  
In particular, Andr\'{e} Hirschowitz's enthusiasm was positive and contagious.  \\
During the preparation of this paper, the second-named author was partially supported by the University of Bologna, funds for selected research topics.

\medskip    

\noindent \textbf{Notations and conventions.}  
We will use the word \textit{universe} in the sense of \cite[Exp. I, Appendice]{sga4}.   
Universes will be denoted by $\mathbb{U} \in \mathbb{V} \in \mathbb{W} \dots$. For any universe  
$\mathbb{U}$ we will assume that $\mathbb{N} \in \mathbb{U}$. The category of sets (resp. simplicial sets, resp. \dots)   
belonging to a universe $\mathbb{U}$ will be denoted by $Set_{\mathbb{U}}$ (resp. $SSet_{\mathbb{U}}$, resp. \dots).   
The objects of $Set_{\mathbb{U}}$ (resp. $SSet_{\mathbb{U}}$, resp. \dots) will be called    
$\mathbb{U}$-sets (resp. $\mathbb{U}$-simplicial sets, resp. \dots). We will use the expression   
\textit{$\mathbb{U}$-small set}  
(resp. \textit{$\mathbb{U}$-small simplicial set}, resp. \dots) to mean \textit{a set isomorphic   
to a set in $\mathbb{U}$}  
(resp. \textit{a simplicial set isomorphic to a simplicial set in $\mathbb{U}$}, resp. \dots).   
  
Our references for model categories are \cite{ho} and \cite{hi}.   
By definition, our model categories will always be \textit{closed} model categories, will have all \textit{small}   
limits and colimits and the functorial factorization property.  
The word \textit{equivalence} will always mean \textit{weak equivalence} and will refer to a model category structure.  
  
The homotopy category of a model category $M$ is $W^{-1}M$ (see \cite[Def. $1.2.1$]{ho}),   
where $W$ is the subcategory of equivalences in $M$, and it   
will be denoted as $\mathrm{Ho}(M)$. The sets of morphisms in $\mathrm{Ho}(M)$   
will be denoted by $[-,-]_{M}$, or simply by  
$[-,-]$ when the reference to the model category $M$ is clear. We will say that two objects in a model category   
$M$ are equivalent if they are isomorphic in $\mathrm{Ho}(M)$.   
We say that two model categories are \textit{Quillen equivalent} if they can be connected by a finite string of   
Quillen adjunctions each one being a  
Quillen equivalence.   
  
The homotopy fibered product (see \cite[\S $11$]{hi} or \cite[Ch. XIV]{dkh}) of a diagram   
$\xymatrix{x \ar[r] & z & \ar[l] y}$ in a model category $M$ will be denoted  
by $x\times_{z}^{h}y$. In the same way, the homotopy push-out of a diagram $\xymatrix{x & \ar[l]  z  \ar[r] & y}$  
will be denoted by $x\coprod_{z}^{h}y$.   
When the model category $M$ is a simplicial model category, its simplicial sets of morphisms will be denoted by  
$\underline{Hom}(-,-)$, and their derived functors by $\mathbb{R}\underline{Hom}$ (see \cite[1.3.2]{ho}).   
  
For the notions of $\mathbb{U}$-cofibrantly generated, $\mathbb{U}$-combinatorial   
and $\mathbb{U}$-cellular model category, we refer to \cite{ho, hi, du2} or   
to Appendix B, where the basic definitions and crucial properties are   
recalled in a way that is suitable for our needs.  
  
As usual, the standard simplicial category will be denoted by $\Delta$. For any simplicial object  
$F \in C^{\Delta^{op}}$ in a category $C$, we will use the notation $F_{n}:=F([n])$. Similarly,  
for any co-simplicial object $F \in C^{\Delta}$, we will use the notation $F_{n}:=F([n])$.

For a Grothendieck site $(C,\tau)$ in a universe $\mathbb{U}$, we will denote by $Pr(C)$ the category of   
presheaves of $\mathbb{U}$-sets  
on $C$, $Pr(C):=C^{Set_{\mathbb{U}}^{op}}$. The subcategory of sheaves on $(C,\tau)$ will be denoted by   
$Sh_{\tau}(C)$, or  
simply by $Sh(C)$ if the topology $\tau$ is unambiguous.   
  
\end{section}  
  
\begin{section}{Review of $S$-categories}  
  
In this first section we recall some facts concerning $S$-categories. The main references on the subject are  
\cite{dk1,dk2,dkh}, except for the material covered in the two final   
subsections for which it does not seem to exist any reference.   
The notion of $S$-category will be of fundamental importance in all this work, as it will  
replace the notion of usual category in our higher sheaf theory.   
In Section 3, we will define what an $S$-topology on an $S$-category is, and study the   
associated notion of stack.   
  
We start by reviewing the definition of $S$\textit{-category}   
and the Dwyer-Kan \textit{simplicial localization} technique.   
We recall the existence of model categories of \textit{diagrams} over $S$-categories,   
as  well as their relations with the model categories of \textit{restricted diagrams}.   
The new materials are presented in  
the last two subsections: here, we first prove a \textit{Yoneda-like lemma}   
for $S$-categories and then introduce and study the notion of \textit{comma $S$-category}.  
  
\begin{subsection}{The homotopy theory of $S$-categories}  
  
We refer to \cite{ke} for the basic notions of enriched category theory.   
We will be especially interested in the case where the enrichement takes place   
in the cartesian closed category $SSet$ of simplicial sets.  
  
\begin{df}\label{d1}  
An $S$\emph{-category} $T$ is a category enriched in $SSet$. A   
\emph{morphism} of $S$-categories $T\rightarrow T'$ is a $SSet$-enriched   
functor.  
\end{df}  
  
\noindent More explicitly, an $S$\textit{-category} $T$ consists of the following data.  
\begin{itemize}  
\item A set $Ob(T)$ (whose elements are called the \textit{objects of $T$}).  
  
\item For any pair of objects $(x,y)$ of $Ob(T)$, a simplicial set   
$\underline{Hom}_{T}(x,y)$ (called the simplicial set of   
\textit{morphisms from $x$ to $y$}). A $0$-simplex in $\underline{Hom}_{T}(x,y)$  
will simply be called a \textit{morphism} from $x$ to $y$ in $T$.   
The $1$-simplices in $\underline{Hom}_{T}(x,y)$ will be called \textit{homotopies}.  
  
\item For any triple of objects $(x,y,z)$ in $Ob(T)$, a morphism  
of simplicial sets (called the \textit{composition} morphism)  
$$\underline{Hom}_{T}(x,y)\times \underline{Hom}_{T}(y,z) \longrightarrow \underline{Hom}_{T}(x,z).$$  
  
\item For any object $x \in Ob(T)$, a $0$-simplex $\mathrm{Id}_{x} \in \underline{Hom}_{T}(x,x)_{0}$  
(called the \textit{identity} morphism at $x$).  
  
\end{itemize}  
  
\noindent These data are required to satisfy the usual associativity and unit axioms. \\  
  
\noindent A \textit{morphism} between $S$-categories $f : T \longrightarrow T'$ consists of the following data.  
  
\begin{itemize}  
\item A map of sets $Ob(T) \longrightarrow Ob(T')$.  
  
\item For any two objects $x$ and $y$ in $Ob(T)$, a morphism of simplicial sets   
$$\underline{Hom}_{T}(x,y) \longrightarrow \underline{Hom}_{T'}(f(x),f(y)),$$  
compatible with the composition and unit in an abvious way.  
  
\end{itemize}  
  
\noindent Morphisms of $S$-categories can be composed in the obvious way,   
thus giving rise to the category of $S$-categories.  
  
\begin{df}\label{d2}  
The category of $S$-categories belonging to a universe $\mathbb{U}$, will be denoted by   
$S-Cat_{\mathbb{U}}$, or simply by $S-Cat$ if the universe $\mathbb{U}$ is clear from the context or irrelevant.  
\end{df}  

The natural inclusion functor $j : Set \longrightarrow SSet$, sending a   
set to the corresponding constant simplicial set, allows us to construct a natural inclusion   
$j : Cat \longrightarrow S-Cat$, and therefore to see any category as an $S$-category.   
Precisely, for a category $C$, $j(C)$ is  
the $S$-category with the same objects as $T$ and whose simplicial set of morphism from $x$ to $y$  
is just the constant simplicial set associated to the set $Hom_{C}(x,y)$.   
In the following we will simply write $C$ for $j(C)$.   
  
Any $S$-category $T$ has an \textit{underlying category of $0$-simplices} $T_{0}$; its set of objects is the same  
as that of $T$ while the set of morphisms from $x$ to $y$ in $T_{0}$ is the set of $0$-simplices of the simplicial set   
$\underline{Hom}_{T}(x,y)$. The construction $T \mapsto T_{0}$ defines a functor  
$S-Cat \longrightarrow Cat$ which is easily checked to be \textit{right adjoint} to the   
inclusion $j : Cat \longrightarrow S-Cat$ mentioned above. This is completely   
analogous to (and actually, induced by) the adjunction between the constant   
simplicial set functor $\mathrm{c}: Set \longrightarrow SSet$ and the $0$-th   
level set functor $(-)_{0} : SSet \longrightarrow Set$. 
  
Any $S$-category $T$ also has a \textit{homotopy category}, denoted by $\mathrm{Ho}(T)$; its set of  
objects is the same as that of $T$, and the set of morphisms from $x$ to $y$ in $\mathrm{Ho}(T)$ is given by    
$\pi_{0}(\underline{Hom}_{T}(x,y))$, the set of connected components of the simplicial set of morphisms from   
$x$ to $y$ in $T$. The construction $T \mapsto \mathrm{Ho}(T)$ defines a functor  
$S-Cat \longrightarrow Cat$ which is easily checked to be \textit{left adjoint} to the inclusion   
$j : Cat \longrightarrow S-Cat$. Again, this is completely analogous to (and actually,   
induced by) the adjunction between the constant simplicial set functor   
$\mathrm{c}: Set \longrightarrow SSet$ and the connected components' functor   
$\pi_{0}: SSet \longrightarrow Set$.\\  
  
Summarizing, we have the following two adjunction pairs (always ordered by writing the left adjoint on the left):  
  
$$j : Cat \longrightarrow S-Cat \qquad Cat \longleftarrow S-Cat : (-)_{0}$$  
$$Ob(T_{0}):=Ob(T) \qquad Hom_{T_{0}}(x,y):=\underline{Hom}_{T}(x,y)_{0}$$  
\medskip  
$$\mathrm{Ho}(-) : S-Cat \longrightarrow Cat \qquad S-Cat \longleftarrow Cat : j $$   
$$Ob(\mathrm{Ho}(T)):=Ob(T) \qquad Hom_{\mathrm{Ho}(T)}(x,y):=\pi_{0}(\underline{Hom}_{T}(x,y)).$$  
  
For an $S$-category $T$, the two associated categories $T_{0}$ and $\mathrm{Ho}(T)$  
are related in the following way. There exist natural morphisms of $S$-categories  
$$\xymatrix{T_{0} \ar[r]^{i} &  T \ar[r]^{p} & \mathrm{Ho}(T)},$$  
which induce a functor $q : T_{0} \longrightarrow \mathrm{Ho}(T)$.   
Being the underlying category of an $S$-category, the category $T_{0}$ has a natural notion  
of \textit{homotopy} between morphisms. This induces an equivalence relation on  
the set of morphisms of $T_{0}$, by declaring two morphisms equivalent if there is a string of homotopies   
between them.  
This equivalence relation is furthermore compatible with composition. The category obtained  
from $T_{0}$ by passing to the quotient with respect to this equivalence relation is  
precisely $\mathrm{Ho}(T)$.    
  
\begin{df}\label{d3}  
Let $f : T \longrightarrow T'$ be a morphism of $S$-categories.  
\begin{enumerate}  
\item The morphism $f$ is \emph{essentially surjective} if the induced functor   
$\mathrm{Ho}(f) : \mathrm{Ho}(T) \longrightarrow \mathrm{Ho}(T')$ is an essentially surjective functor of categories.  
  
\item The \emph{essential image} of $f$ is the inverse image by the natural projection  
$T' \longrightarrow \mathrm{Ho}(T')$ of the essential image of   
$\;\mathrm{Ho}(f) : \mathrm{Ho}(T) \longrightarrow \mathrm{Ho}(T')$.  
  
\item The morphism $f$ is \emph{fully faithful} if for any pair of objects $x$ and $y$ in $T$, the induced  
morphism $f_{x,y} : \underline{Hom}_{T}(x,y) \longrightarrow \underline{Hom}_{T'}(f(x),f(y))$  
is an equivalence of simplicial sets.  
  
\item The morphism $f$ is an \emph{equivalence} if it is essentially surjective and fully faithful.   
  
\end{enumerate}  
The category obtained from    
$S-Cat$ by formally inverting the equivalences will be denoted by $\mathrm{Ho}(S-Cat)$.   
The set of morphisms in $\mathrm{Ho}(S-Cat)$ between two objects $T$ and $T'$ will simply be   
denoted by $[T,T']$.   
\end{df}  
\begin{rmk}\emph{\cite[\S XII-48]{dkh} contains the sketch of a proof that the category $S-Cat$ admits   
a model structure whose equivalences are exactly the ones defined   
above. It seems however that this proof is not complete, as the generating   
trivial cofibrations of \cite[$48.5$]{dkh}  
fail to be equivalences. In his note \cite[Thm. 1.9]{may2},
P. May informed us that he knows an alternative proof, but
the reader will notice that the notion of fibrations
used in \cite{may2} is different from the one used
in \cite{dkh} and does not seem to be correct. 
We think however that the model structure described in 
\cite{dkh} exists\footnote{Progresses in this direction have been recently made by J. Bergner (private communication).}, as we have the feeling that one
could simply replace the wrong set of generating trivial cofibrations
by the set of all trivial cofibrations between 
countable $S$-categories. The existence of this model structure
would of course simplify some of our constructions, but 
it does not seem to be really unavoidable, and
because of the lack of clear references we have decided not 
to use it at all. This will cause a 
``lower degree'' of functoriality in some constructions, but
will be enough for all our purposes.} 
\end{rmk}

Since the natural localization functor $SSet \longrightarrow \mathrm{Ho}(SSet)$   
commutes with finite products, any category  
enriched in $SSet$ gives rise to a category enriched in $\mathrm{Ho}(SSet)$.   
The $\mathrm{Ho}(SSet)$-enriched category  
associated to an $S$-category $T$ will be denoted
by $\underline{\mathrm{Ho}(T)}$, and
has $\mathrm{Ho}(T)$ as underlying category.   
Furthermore, for any pair of objects $x$ and $y$   
in $\mathrm{Ho}(T)$, one has  
$\underline{Hom}_{\underline{\mathrm{Ho}(T)}}(x,y)=\underline{Hom}_{T}(x,y)$ considered  as objects in $\mathrm{Ho}(SSet)$. 
Clearly, $T \mapsto \underline{\mathrm{Ho}}(T)$
defines a functor from $S-Cat$ to the category $\mathrm{Ho}(SSet)-Cat$ 
of $\mathrm{Ho}(SSet)$-enriched categories, and
a morphism of $S$-categories is an equivalence if and only if
the induced $\mathrm{Ho}(SSet)$-enriched functor is an 
$\mathrm{Ho}(SSet)$-enriched equivalence. Therefore, this
construction induces a well defined functor
$$\begin{array}{ccc}
\mathrm{Ho}(S-Cat) & \longrightarrow & Ho(\mathrm{Ho}(SSet)-Cat) \\
T &\mapsto & \underline{\mathrm{Ho}(T)},
\end{array}$$
where $\mathrm{Ho}(\mathrm{Ho}(SSet)-Cat)$ is the 
localization of the category of $\mathrm{Ho}(SSet)$-enriched
categories along $\mathrm{Ho}(SSet)$-enriched equivalences. 

The previous construction allows one to define the notions of   
essentially surjective and fully faithful morphisms  
in $\mathrm{Ho}(S-Cat)$. Precisely, a morphism
$f : T \longrightarrow T'$ in $\mathrm{Ho}(S-Cat)$   
will be called essentially surjective (resp. fully faithful) if  
the corresponding $\mathrm{Ho}(SSet)$-enriched functor
$\underline{\mathrm{Ho}(f)} : \underline{\mathrm{Ho}(T)} 
\longrightarrow \underline{\mathrm{Ho}(T')}$
is essentially surjective (resp. fully faithful)   
in the $\mathrm{Ho}(SSet)$-enriched sense. 

Finally, for an $S$-category $T$ and a property $\mathbf{P}$ of morphisms in $\mathrm{Ho}(T)$, we will often say that   
\textit{a morphism $f$ in $T$ satisfies the property $\mathbf{P}$} to mean that    
\textit{the image of f in $\mathrm{Ho}(T)$ through the natural projection $T \longrightarrow \mathrm{Ho}(T)$,   
satisfies the property $\mathbf{P}$}. Recall that a \textit{morphism} $f$ in an   
$S$-category  $T$ is just an element in the zero simplex set of $\underline{Hom}_{T}(x,y)$  
for some $x$ and $y$ in $Ob(T)$.  
  
\end{subsection}  
  
\begin{subsection}{Simplicial localization}  
  
Starting from a category $C$ together with a subcategory $S \subset C$, W. Dwyer and D. Kan have defined  
in \cite{dk1} an $S$-category $\mathrm{L}(C,S)$, which is an enhanced version of the localized category $S^{-1}C$.   
It is an $S$-category with a diagram of morphisms in $S-Cat$ (viewing, according to our   
general conventions, any category as an $S$-category via the embedding $j:Cat \rightarrow S-Cat$)  
$$\xymatrix{C & \ar[l]_-{p}F_{*}C \ar[r]^-{L} & \mathrm{L}(C,S)}$$  
where $F_{*}C$ is the so-called \textit{standard simplicial free resolution   
of the category $C$}, and in particular, the projection   
$p$ is an equivalence of $S$-categories. Therefore, there exists a well   
defined \textit{localization morphism} in $\mathrm{Ho}(S-Cat)$  
$$L : C \longrightarrow \mathrm{L}(C,S).$$  
The construction $(C,S) \mapsto \mathrm{L}(C,S)$ is functorial in the pair $(C,S)$ and it also   
extends naturally to the case where $S$ is a sub-$S$-category  
of an $S$-category $C$ (see \cite[\S $6$]{dk1}). Note also that by construction, if  
$C$ belongs to a universe $\mathbb{U}$ so does $L(C,S)$. \\  
  
\begin{rmk}  
\emph{  
\begin{enumerate}  
\item   
One can also check that the localization morphism $L$ satisfies the following universal   
property.  
For each $S$-category $T$, let us denote by $[C,T]^{S}$ the subset of $[C,T]=Hom_{\mathrm{Ho}(S-Cat)}(C,T)$ consisting  
of morphisms for which the induced morphism $C \longrightarrow \mathrm{Ho}(T)$  
sends morphisms of $S$ into isomorphisms in $\mathrm{Ho}(T)$ (the reader will easily check   
that this property is well defined). Then the localization morphism $L$ is such that for any $S$-category $T$   
the induced map  
$$L^{*} : [\mathrm{L}(C,S),T] \longrightarrow [C,T]$$  
is injective and its image is $[C,T]^{S}$. This property characterizes the  
$S$-category $\mathrm{L}(C,S)$ as an object in the comma category  
$C/\mathrm{Ho}(S-Cat)$.  This universal property will not be used in   
the rest of the paper, but we believe it makes the meaning of the simplicial localization more transparent.  
\item It is important to mention the fact that any $S$-category $T$ is equivalent to some  
$L(C,S)$, for a category $C$ with a subcategory $S \subset C$ (this is the   
\textit{delocalization theorem}  of \cite{dk2}). Furthermore, it is clear by the construction given in \cite{dk1}   
that, if $T$ is $\mathbb{U}$-small, then so are $C$, $S$ and $L(C,S)$.  
\end{enumerate}  
}  
\end{rmk}  
  
\medskip  
  
\noindent Two fundamental properties of the functor $L:(C,S) \mapsto \mathrm{L}(C,S)$ are the following  
  
\begin{enumerate}  
  
\item The localization morphism $L$ induces a well defined (up to a unique isomorphism) functor   
$$\mathrm{Ho}(L) : C \simeq \mathrm{Ho}(F_{*}C) \longrightarrow \mathrm{Ho}(\mathrm{L}(C,S)),$$  
that identifies $\mathrm{Ho}(\mathrm{L}(C,S))$ with the (usual Gabriel-Zisman) localization $S^{-1}C$.  
  
\item Let $M$ be a simplicial model category, $W\subset M$ its subcategory of equivalences and  
let $\mathrm{Int}(M)$ be the $S$-category of fibrant and cofibrant objects in $M$ together with their  
simplicial sets of morphisms. The full (not simplicial) subcategory $M^{cf}\subset M$ of fibrant and cofibrant objects  
in $M$ has two natural morphisms in $S-Cat$  
$$\xymatrix{M & \ar[l] M^{cf} \ar[r] & \mathrm{Int}(M),}$$  
which induce isomorphisms in $\mathrm{Ho}(S-Cat)$  
$$\mathrm{L}(M,W) \simeq \mathrm{L}(M^{cf},W\cap M^{cf}) \simeq   
\mathrm{L}(\mathrm{Int}(M),W\cap M^{cf})\simeq \mathrm{Int}(M).$$  
  
In the same way, if $M^{f}$ (resp. $M^{c}$) is the full subcategory   
of fibrant (resp. cofibrant) objects  
in $M$, the natural morphisms $M^{f}\longrightarrow M$,   
$M^{c}\longrightarrow M$ induce isomorphisms in $\mathrm{Ho}(S-Cat)$  
$$\mathrm{L}(M^{f},W\cap M^{f}) \simeq \mathrm{L}(M,W) \qquad   
\mathrm{L}(M^{c},W\cap M^{c})\simeq \mathrm{L}(M,W).$$  
  
\end{enumerate}  
  
\begin{df}\label{d4}  
If $M$ is any model category, we set $LM:=\mathrm{L}(M,W)$, where   
$W \subset M$ is the subcategory of equivalences in $M$.   
\end{df}  
  
The construction $M \mapsto LM$ is functorial, up to equivalences,   
for Quillen functors between model categories. To see this, let $f : M \longrightarrow N$ be  
a right Quillen functor. Then, the restriction to the category of fibrant objects   
$f : M^{f} \longrightarrow N^{f}$ preserves equivalences, and   
therefore induces a morphism of $S$-categories  
$$Lf : LM^{f} \longrightarrow LN^{f}.$$  
Using the natural isomorphisms $LM^{f}\simeq LM$ and $LN^{f}\simeq LN$   
in $\mathrm{Ho}(S-Cat)$, one gets a well defined  
morphism $Lf : LM \longrightarrow LN$. This is a morphism in the   
homotopy category $\mathrm{Ho}(S-Cat)$, and one checks immediately  
that $M \mapsto LM$ is a functor from the category of model categories   
(belonging to a fixed universe $\mathbb{U}$) with    
right Quillen functors, to $\mathrm{Ho}(S-Cat_{\mathbb{U}})$.   
The dual construction gives rise to a functor $M \mapsto LM$ from the   
category of model categories which belongs to a universe $\mathbb{U}$   
and left Quillen functors to $\mathrm{Ho}(S-Cat_{\mathbb{U}})$.   
  
The reader will check easily that if  
$$f : M \longrightarrow N \qquad M \longleftarrow N : g$$  
is a Quillen adjunction which is a Quillen equivalence, then the morphisms $Lf : LM \longrightarrow LN$ and  
$Lg : LN \longrightarrow LM$ are isomorphisms inverse to each others in $\mathrm{Ho}(S-Cat)$. \\  
  
\end{subsection}  
  
\begin{subsection}{Model categories of diagrams}  
  
In this paragraph we discuss the notion of \textit{pre-stack over an $S$-category} which is a generalization of  
the notion of presheaf of sets on a usual category.  
  
\begin{subsubsection}{Diagrams} \label{diag}  
  
Let $T$ be any $S$-category in a universe $\mathbb{U}$, and $M$ a simplicial model category which is  
$\mathbb{U}$-cofibrantly generated (see \cite[13.2]{hi} and Appendix A). Since $M$ is simplicial, we may view it as  
an $S$-category, with the same set of objects as $M$ and whose simplicial sets of morphisms are provided by   
the simplicial structure.   
Therefore, we may consider the category $M^{T}$, of morphisms of $S$-categories  
$F : T \longrightarrow M$. To be more precise, an \textit{object}   
$F : T \longrightarrow M$ in $M^{T}$ consists of the following data  
  
\begin{itemize}  
  
\item A map $F : Ob(T) \longrightarrow Ob(M)$.   
  
\item For any pair of objects $(x,y)\in Ob(T) \times Ob(T)$, a morphism of simplicial sets  
$$F_{x,y} : \underline{Hom}_{T}(x,y) \longrightarrow \underline{Hom}(F(x),F(y)),$$  
(or equivalently, morphisms $F_{x,y} : \underline{Hom}_{T}(x,y)\otimes F(x) \longrightarrow F(y)$ in $M$)  
satisfying the obvious associativity and unit axioms.   
  
\end{itemize}  
  
A \textit{morphism} from $F$ to $G$ in $M^{T}$ consists of morphisms   
$H_{x} : F(x) \longrightarrow G(x)$ in $M$, for all $x \in Ob(T)$,   
such that the following diagram commutes in $M$  
$$\xymatrix{  
\underline{Hom}_{T}(x,y)\otimes F(x) \ar[r]^-{F_{x,y}} \ar[d]_-{\mathrm{Id}\otimes H_{x}} & F(y) \ar[d]^-{H_{y}} \\  
\underline{Hom}_{T}(x,y)\otimes G(x) \ar[r]_-{G_{x,y}} & G(y).}$$  
  
One defines a model structure on $M^{T}$, by defining a morphism $H$ to be a fibration (resp. an equivalence) if  
for all $x \in Ob(T)$, the induced morphism $H_{x}$ is a fibration (resp. an equivalence) in $M$. It is known that  
these definitions make $M^{T}$ into a simplicial model category which is   
again $\mathbb{U}$-cofibrantly generated (see \cite[$13.10.17$]{hi} and Appendix A). This model structure will be  
called the \textit{projective model structure on $M^{T}$}. Equivalences and fibrations in $M^{T}$ will be called   
\textit{objectwise equivalences} and \textit{objectwise fibrations}. \\  
  
Let us suppose now that $M$ is an \textit{internal} model category (i.e. a symmetric monoidal model category  
for the direct product, in the sense of \cite[Ch. 4]{ho}).  
The category $M^{T}$ is then naturally tensored and co-tensored over $M$.   
Indeed, the external product $A\otimes F \in M^{T}$ of $A \in M$   
and $F\in M^{T}$, is simply defined by the formula $(A\otimes F)(x):=A\times F(x)$ for any $x \in Ob(T)$. For any  
$x$ and $y$ in $Ob(T)$, the transition morphisms of $A\otimes F$ are defined by  
$$(A\otimes F)_{x,y}:=A\times F_{x,y} : \underline{Hom}_{T}(x,y)\times   
A\times F(x)\simeq A\times \underline{Hom}_{T}(x,y)\times F(x)  
\longrightarrow A\times F(y).$$  
In the same way,   
the exponential $F^{A}\in M^{T}$ of $F$ by $A$, is defined by $(F^{A})(x):=F(x)^{A}$ for any $x$ in $Ob(T)$.  
  
With these definitions the model category $M^{T}$ becomes a $M$-model category in the sense of  
\cite[Def. $4.2.18$]{ho}. When $M$ is the model category of simplicial sets, this implies that $SSet^{T}$ has  
a natural structure of simplicial model category where exponential   
and external products are defined levelwise. In particular, for any $x\in Ob(T)$, the  
evaluation functor  
$$\begin{array}{cccc}j_{x}^{*} :  & M^{T} & \longrightarrow & M \\  
& F & \mapsto & F(x)  
\end{array}$$  
commutes with the geometric realization and total space functors of   
\cite[\S $19.5$]{hi}. As fibrant (resp. cofibrant) objects  
in $M^{T}$ are also objectwise fibrant (resp. objectwise cofibrant), this easily implies that $j_{x}^{*}$  
commutes, up to an equivalence, with homotopy limits and homotopy colimits. One may also directly shows that   
$j_{x}^{*}$ is indeed a left and right Quillen functor.   
Finally, if $M$ is a proper model category, then so is $M^{T}$. \\  
  
Let $f : T \longrightarrow T'$ be a morphism in $S-Cat_{\mathbb{U}}$. It gives rise to an adjunction  
$$f_{!} : M^{T} \longrightarrow M^{T'} \qquad M^{T} \longleftarrow M^{T'} : f^{*},$$  
where $f^{*}$ is defined by the usual formula $f^{*}(F)(x):=F(f(x))$, for any $F \in M^{T'}$ and any $x \in Ob(T)$,   
and $f_{!}$ is its left adjoint. The functor $f^{*}$ is clearly a right Quillen functor, and therefore  
$(f_{!},f^{*})$ is a Quillen adjunction. \\  
  
The following theorem is proved in \cite{dk2} when $M$ is the category of simplicial sets; its proof generalizes   
immediately to our situation. As above, $M$ is a simplicial $\mathbb{U}$-cofibrantly generated model category.  
  
\begin{thm}\label{t1}  
If $f:T \rightarrow T'$ is an equivalence of $S$-categories, then   
$(f_{!},f^{*})$ is a Quillen equivalence of model categories.  
\end{thm}  
  
\begin{df}\label{d5}  
Let   
$T \in S-Cat_{\mathbb{U}}$ be an $S$-category in $\mathbb{U}$, and  
$M$ a $\mathbb{U}$-cofibrantly generated simplicial model category.   
The model category $Pr(T,M)$ of \emph{pre-stacks on} $T$ \emph{with values in} $M$ is defined as   
$$Pr(T,M):=M^{T^{op}}.$$  
We will simply write $SPr(T)$ for $Pr(T,SSet_{\mathbb{U}})$, and call it  
the model category of \emph{pre-stacks on} $T$.   
\end{df}  
  
Theorem \ref{t1} implies that the model category $Pr(T,M)$, for a fixed $M$, is an invariant, up to Quillen equivalence,  
of the isomorphism class of $T$ in $\mathrm{Ho}(S-Cat_{\mathbb{U}})$.  
In the same way, if $f : T \longrightarrow T'$  
is a morphism in $\mathrm{Ho}(S-Cat_{\mathbb{U}})$, one can represent   
$f$ by a string of morphisms in $S-Cat_{\mathbb{U}}$  
$$\xymatrix{T & \ar[l]_-{p_{1}} T_{1}  \ar[r]^-{f_{1}} & 
T_{2} & \ar[l]_-{p_{3}} T_{3} \ar[r]^-{f_{3}} & T_{4} & \cdots & &  
 \ar[l]_-{p_{2n-1}} T_{2n-1}\ar[r]^-{f_{2n-1}} & T',}$$  
where each $p_{i}$ is an equivalence of $S$-categories. 
We deduce a diagram  
of right Quillen functors  
$$\xymatrix{ Pr(T,M) \ar[r]^-{p_{1}^{*}} & Pr(T_{1},M) & \ar[l]_-{f_{1}^{*}} Pr(T_{2},M) \ar[r]^-{p_{3}^{*}} & Pr(T_{3},M) & \cdots & \ar[r]^-{p_{2n-1}^{*}} & Pr(T_{2n-1},M) & \ar[l]_-{f_{2n-1}^{*}} Pr(T',M),}$$  
such that each $p_{i}^{*}$ is a right adjoint of a Quillen equivalence. By definition, this diagram   
gives a 
\textit{Quillen adjunction between $Pr(T,M)$ and $Pr(T',M)$, up to Quillen equivalences}, which can also be interpreted as a morphism
in the category of model categories localized along
Quillen equivalences.   
In particular, we obtain a well defined morphism
in $\mathrm{Ho}(S-Cat)$  
$$\mathbb{R}f^{*}:=
(p_{1}^{*})^{-1} \circ (f_{1}^{*})\circ \dots \circ (p_{2n-1}^{*})^{-1}
\circ (f_{2n-1}^{*}): 
LPr(T',M) \longrightarrow LPr(T,M).$$
Using direct images (i.e. functors $(-)_{!}$) instead of inverse images, one also gets
a morphism in the other direction 
$$\mathbb{L}f_{!}:=
(f_{2n-1})_{!}\circ (p_{2n-1})_{!}^{-1} \circ \dots \circ (f_{1})_{!} \circ(p_{1})_{!}^{-1}: 
LPr(T,M) \longrightarrow LPr(T',M)$$ (again well defined
in $\mathrm{Ho}(S-Cat)$).
Passing to the associated $\mathrm{Ho}(SSet)$-enriched
categories, one obtains a
$\mathrm{Ho}(SSet)$-enriched adjunction
$$\mathbb{L}f_{!} : \underline{\mathrm{Ho}(Pr(T,M))}
\longrightarrow \underline{\mathrm{Ho}(Pr(T',M))} \qquad
\underline{\mathrm{Ho}(Pr(T,M)})
\longleftarrow \underline{\mathrm{Ho}(Pr(T',M))} : \mathbb{R}f^{*}.$$
The two $\mathrm{Ho}(SSet)$-enriched functors are well
defined up to a unique isomorphism. 
When $M$ is fixed, the construction above defines
a well defined functor from 
the category $\mathrm{Ho}(S-Cat)$ to the 
homotopy category of $\mathrm{Ho}(SSet)$-enriched
adjunctions. 
  
\end{subsubsection}  
  
\begin{subsubsection}{Restricted diagrams}\label{resdia}  
  
Let $C$ be a $\mathbb{U}$-small $S$-category, $S \subset C$ a sub-$S$-category, and $M$ a  
simplicial model category which is $\mathbb{U}$-cofibrantly generated. We will assume also that $M$ is a $\mathbb{U}$-combinatorial or $\mathbb{U}$-cellular model category so that   
the left Bousfield localization techniques of \cite[Ch. $4$]{hi} can be applied to homotopically invert any $\mathbb{U}$-set of morphisms (see Appendix A).  
  
We consider the model category $M^{C}$, of simplicial functors from $C$ to $M$, endowed with its projective model   
structure. For any object $x\in C$, the evaluation functor $i_{x}^{*} : M^{C} \longrightarrow M$,   
defined by $i_{x}^{*}(F):=F(x)$, has a left adjoint $(i_{x})_{!} : M \longrightarrow M^{C}$ which  
is a left Quillen functor. Let $I$ be a $\mathbb{U}$-set of generating cofibrations in $M$.   
For any $f : A \longrightarrow B$ in $I$ and any morphism   
$u : x \longrightarrow y$ in $S\subset C$, one consider the natural morphism  
in $M^{C}$  
$$f\square u : (i_{y})_{!}(A)\coprod_{(i_{x})_{!}(A)}(i_{x})_{!}(B) \longrightarrow (i_{y})_{!}(B).$$  
  
Since $M$ is a $\mathbb{U}$-combinatorial (or $\mathbb{U}$-cellular)   
model category, then so is $M^{C}$ (see \cite{du2,hi} and Appendix A). As the set  
of all $f\square u$, for $f \in I$ and $u$ a morphism in $S$, belongs to $\mathbb{U}$, the following definition   
is well posed.  
  
\begin{df}\label{d7}  
The model category $M^{C,S}$ is the left Bousfield localization of $M^{C}$ along the set   
of all morphisms $f\square u$, where $f \in I$ and $u$ is a morphism in  $S$.   
  
The model category $M^{C,S}$ will be called the \emph{model category of restricted diagrams} from $(C,S)$ to $M$.   
\end{df}  
  
\begin{rmk} \label{earth}  
\emph{If $M=SSet_{\mathbb{U}}$, we may take $I$ to be the usual set of generating cofibrations  
$$I=\left\{ f_{n}:\partial\Delta[n] \hookrightarrow \Delta[n]\;|\;n\geq 0\right\}.$$   
Since as it is easily checked, we have a canonical isomorphism   
$(i_{x})_{!}(*=\Delta[0])\simeq \underline{h}_{x}$ in $SSet^{(C,S)^{op}}$, for any $x \in C$,   
where $\underline{h}_{x}$ denotes the simplicial diagrams defined by $\underline{h}_{x}(y):=\underline{Hom}_{T}(y,x)$.   
Then, for any $u:x\rightarrow y$ in $S$, we have that the set of morphisms   
$f_{n}\square u$ is exactly   
the set of augmented horns on the set of morphisms $\underline{h}_{x} \rightarrow \underline{h}_{y}$ (see \cite[\S $4.3$]{hi}).   
This implies that $SSet^{C,S}$ is simply the left Bousfield localization of $SSet^{C}$ along the  
set of morphisms $\underline{h}_{x} \rightarrow \underline{h}_{y}$ for any $x \rightarrow y$ in $S$.}  
\end{rmk}  
  
By the general theory of left Bousfield localization of \cite{hi},   
the fibrant objects in the model category $M^{C,S}$ are the functors $F : C \longrightarrow M$  
satisfying the following two conditions:  
  
\begin{enumerate}  
  
\item For any $x \in C$, $F(x)$ is a fibrant object in $M$ (i.e.   
$F$ is fibrant in $M^{C}$ for the projective model structure).  
  
\item For any morphism $u : x\longrightarrow y$ in $S$, the induced   
morphism $F_{x,y}(u) : F(x) \longrightarrow F(y)$ is  
an equivalence in $M$.  
\end{enumerate}  
  
Now, let $(F_{*}C,F_{*}S)$ be the canonical free resolution of   
$(C,S)$ in $S-Cat_{\mathbb{U}}$ (see \cite{dk1}). Then, one  
has a diagram of pairs of $S$-categories  
$$\xymatrix{ (C,S) & \ar[l]_-{p} (F_{*}C,F_{*}S) \ar[r]^-{l} & (F_{*}S)^{-1}(F_{*}C)=\mathrm{L}(C,S),}$$  
inducing a diagram of right Quillen functors  
$$\xymatrix{\ar[r]^-{p^{*}} M^{C,S} &  M^{F_{*}C,F_{*}S} & \ar[l]_-{l^{*}} M^{\mathrm{L}(C,S)}}.$$  
The following result is proved  
in \cite{dk2} in the case where $M=SSet_{\mathbb{U}}$, and its proof generalizes easily to our situation.   
  
\begin{thm}\label{t2}  
The previously defined right Quillen functors $p^{*}$ and $l^{*}$ are Quillen equivalences. In particular  
the two model categories $M^{\mathrm{L}(C,S)}$ and $M^{C,S}$ are Quillen equivalent.   
\end{thm}  
  
The model categories of restricted diagrams are functorial in the following sense. Let $f : C \longrightarrow D$  
be a functor between two $\mathbb{U}$-small $S$-categories, and let $S \subset C$ and $T\subset D$ be two sub-$S$-categories  
such that $f(S)\subset T$. The functor $f$ induces the usual adjunction on the categories of diagrams in $M$  
$$f_{!} : M^{C,S} \longrightarrow M^{D,T} \qquad M^{C,S} \longleftarrow M^{D,T} : f^{*}.$$  
The adjunction $(f_{!},f^{*})$ is a Quillen adjunction for the objectwise model structures. Furthermore,   
using the description of fibrant objects given above, it is clear   
that $f^{*}$ sends fibrant objects in $M^{D,T}$ to fibrant objects in $M^{C,S}$. By   
the general formalism of left Bousfield localizations (see \cite[\S 3]{hi}), this implies that $(f_{!},f^{*})$ is also a   
Quillen adjunction for the restricted model structures.   
  
\begin{cor}\label{c0}  
Let $f : (C,S) \longrightarrow (D,T)$ be as above. If the induced morphism of   
$S$-categories $Lf : \mathrm{L}(C,S) \longrightarrow \mathrm{L}(D,T)$  
is an equivalence, then the Quillen adjunction $(f_{!},f^{*})$ is a Quillen equivalence between $M^{C,S}$ and $M^{D,T}$.   
\end{cor}  
  
\textit{Proof:} This is a consequence of theorems \ref{t1} and \ref{t2}. \hfill \textbf{$\Box$} 
  
\end{subsubsection}  
  
\end{subsection}

\begin{subsection}{The Yoneda embedding}\label{SYoneda}  
  
In this paragraph we define a Yoneda embedding for $S$-categories. To be   
precise it will be only defined as a morphism in $S-Cat$  
for fibrant $S$-categories, i.e. for $S$-categories whose simplicial sets   
of morphisms are all fibrant; for arbitrary $S$-categories,  
the Yoneda embedding will only be defined as a morphism in the homotopy category $\mathrm{Ho}(S-Cat)$. 
  
We fix $T$, a $\mathbb{U}$-small $S$-category. The category $SPr(T)$   
(see Definition \ref{d5}) is naturally enriched over $SSet$ and the corresponding    
$S$-category will be denoted by $SPr(T)_{s}$. Note that $\mathrm{Int}(SPr(T))$ is a full  
sub-$S$-category of $SPr(T)_{s}$ (recall that $\mathrm{Int}(SPr(T))$ is the $S$-category of  
fibrant and cofibrant objects in the simplicial model category $SPr(T)$).  
  
\noindent Recall the following $SSet$-enriched version of Yoneda lemma (e.g., \cite[IX Lemma $1.2$]{gj})  
  
\begin{prop}\label{non-h-Yoneda}  
Let $T$ be an $S$-category. For any object $x$ in $T$,    
let us denote by $\underline{h}_{x}$ the object in $SPr(T)_{s}$ defined by   
$\underline{h}_{x}(y):=\underline{Hom}_{T}(y,x)$. Then, for any simplicial functor   
$F:T \rightarrow SSet$, there is a canonical isomorphism   
of simplicial sets $$F(x)\simeq \underline{Hom}_{SPr(T)_{s}}(\underline{h}_{x},F)$$  
which is functorial in the pair $(F,x)$.  
\end{prop}  
  
Then, for any $T \in S-Cat_{\mathbb{U}}$, one defines a morphism of $S$-categories  
$\underline{h} : T\longrightarrow SPr(T)_{s}$, by setting for $x \in Ob(T)$  
$$\begin{array}{cccc}  
\underline{h}_{x} : & T^{op} & \longrightarrow & SSet_{\mathbb{U}} \\  
& y & \mapsto & \underline{Hom}_{T}(y,x).  
\end{array}$$  
  
Note that Proposition \ref{non-h-Yoneda} defines immediately $\underline{h}$ at   
the level of morphisms between simplicial $Hom$'s and shows that   
$\underline{h}$ is fully-faithful (in a strong sense) as a morphism in $S-Cat_{\mathbb{V}}$.  
Now, the morphism $\underline{h}$ induces a functor between the associated   
homotopy categories that we will still denote by  
$$\underline{h} : \mathrm{Ho}(T) \longrightarrow \mathrm{Ho}(SPr(T)_{s}).$$  
Now, we want to compare $\mathrm{Ho}(SPr(T)_{s})$ to $\mathrm{Ho}(SPr(T))$;   
note that the two $\mathrm{Ho}(-)$'s here have different meanings, as the   
first one refers to the homotopy category of an $S$-category while the   
second one to the homotopy category of a model category. By definition,   
in the set of morphisms between $F$ and $G$ in $\mathrm{Ho}(SPr(T)_{s})$,   
simplicially homotopic maps in $Hom_{SPr(T)}(F,G)=\underline{Hom}_{SPr(T)_{s}}(F,G)_{0}$,   
give rise to the same element. Then, since simplicially homotopic maps in $Hom_{SPr(T)}(F,G)$   
have the same image in $\mathrm{Ho}(SPr(T))$ (see, for example, \cite[Cor. 10.4.5]{hi}),   
the identity functor induces a well defined localization morphism  
$$\mathrm{Ho}(SPr(T)_{s}) \longrightarrow \mathrm{Ho}(SPr(T)).$$  
Composing this with the functor $\underline{h}$, one deduces a well defined functor (denoted with the same symbol)  
$$\underline{h} : \mathrm{Ho}(T) \longrightarrow \mathrm{Ho}(SPr(T)).$$  
  
The following is a homotopy version of the enriched Yoneda lemma   
(i.e. a homotopy variation of Proposition \ref{non-h-Yoneda})  
  
\begin{prop}\label{l1}  
For any object $F \in SPr(T)$ and any $x \in Ob(T)$,   
there exists an isomorphism in $\mathrm{Ho}(SSet_{\mathbb{U}})$  
$$F(x)\simeq \mathbb{R}\underline{Hom}(\underline{h}_{x},F)$$  
which is functorial in the pair $(F,x) \in \mathrm{Ho}(SPr(T))\times \mathrm{Ho}(T)$.  
In particular, the functor $\underline{h} : \mathrm{Ho}(T) \longrightarrow \mathrm{Ho}(SPr(T))$ is fully faithful.  
\end{prop}  
  
\textit{Proof:} Using Proposition \ref{non-h-Yoneda}, since equivalences in $SPr(T)$   
are defined objectwise, by taking a fibrant replacement of $F$, we may suppose that $F$ is fibrant.   
Moreover, again by Proposition \ref{non-h-Yoneda}, the unique morphism $* \rightarrow \underline{h}_{x}$ has  
the right lifting property with respect to all trivial fibrations, hence $\underline{h}_{x}$ is a   
cofibrant object in $SPr(T)$. Therefore, for any   
fibrant object $F \in SPr(T)$, one has  
natural isomorphisms in $\mathrm{Ho}(SSet_{\mathbb{U}})$  
$$F(x)\simeq \underline{Hom}(\underline{h}_{x},F)\simeq\mathbb{R}\underline{Hom}(\underline{h}_{x},F).$$\hfill \textbf{$\Box$}   
  
The following corollary is a refined version of Proposition \ref{l1}.   
  
\begin{cor}\label{c1}  
Let $T$ be an $S$-category in $\mathbb{U}$ with fibrant
simplicial Hom-sets. 
Then, the morphism $\underline{h} : T \longrightarrow SPr(T)_{s}$  
factors through $\mathrm{Int}(SPr(T))$ and the induced morphism    
$\underline{h} : T \longrightarrow \mathrm{Int}(SPr(T))$ in $S-Cat$  
is fully faithful.   
\end{cor}  
   
\textit{Proof:} The assumption on $T$ implies that $\underline{h}_{x}$   
is fibrant and cofibrant in $SPr(T)$, for any $x \in Ob(T)$ and therefore that $\underline{h}$ factors through  
$\mathrm{Int}(SPr(T))\subset SPr(T)_{s}$. Finally, Proposition \ref{l1} immediately implies   
that $\underline{h}$ is fully faithful. Actually, this is already true for   
$\underline{h} : T \longrightarrow SPr(T)_{s}$, by Proposition \ref{non-h-Yoneda},   
and hence this is true for our factorization since $\mathrm{Int}(SPr(T))$ is a full  
sub-$S$-category of $SPr(T)_{s}$. \hfill \textbf{$\Box$} \\  
  
In case $T$ is an \textit{arbitrary} $S$-category in $\mathbb{U}$ (possibly with non-fibrant simplicial Hom sets), one can consider  
a fibrant replacement $j : T \longrightarrow RT$, defined
by applying the Kan $Ex^{\infty}$-construction (\cite{gj}, III.4) to 
each simplicial set of morphisms in $T$, together with  
its Yoneda embedding  
$$\xymatrix{T \ar[r]^-{j} & RT \ar[r]^-{\underline{h}} & \mathrm{Int}(SPr(RT)).}$$  
When viewed in $\mathrm{Ho}(S-Cat_{\mathbb{V}})$, this induces a well defined morphism  
$$\xymatrix{T \ar[r]^-{j} & RT \ar[r]^-{\underline{h}} & \mathrm{Int}(SPr(RT))\simeq LSPr(RT).}$$  
Finally, composing with the isomorphism $Lj_{!}=(j^{*})^{-1} : LSPr(RT) \simeq LSPr(T)$ of Theorem \ref{t1},  
one gets a morphism   
$$L\underline{h} : T \longrightarrow LSPr(T).$$   
   
This is a morphism in $\mathrm{Ho}(S-Cat_{\mathbb{V}})$, called the   
$S$-\textit{Yoneda embedding of} $T$; when no confusion is possible,   
we will simply call it the Yoneda embedding of $T$. Now, Corollary \ref{c1} immediately  
implies that $L\underline{h}$ is fully faithful, and is indeed isomorphic   
to the morphism $\underline{h}$ defined above when $T$ has  
fibrant simplicial Hom-sets.   
   
\begin{df}\label{d6}  
Let $T$ be an $S$-category.  
An object in $\mathrm{Ho}(SPr(T))$ is called   
\emph{representable} if it belongs to the essential image (see Definition \ref{d3}, 2.) of the functor  
$L\underline{h} : T \longrightarrow LSPr(T)$.  
\end{df}  
  
For any $T \in \mathrm{Ho}(S-Cat_{\mathbb{U}})$, the Yoneda embedding  
$L\underline{h} : T \longrightarrow LSPr(T)$  
induces an isomorphism in $\mathrm{Ho}(S-Cat_{\mathbb{V}})$ between   
$T$ and the full sub-$S$-category of $LSPr(T)$ consisting  
of representable objects.   
  
Note that the functor induced on the level of homotopy categories  
$$L\underline{h} : \mathrm{Ho}(T) \longrightarrow \mathrm{Ho}(LSPr(T))=\mathrm{Ho}(SPr(T))$$  
simply sends $x\in Ob(T)$ to the simplicial presheaf $\underline{h}_{x} \in \mathrm{Ho}(SPr(T))$.    
  
\end{subsection}  
  
\begin{subsection}{Comma $S$-categories}  
  
In this subsection we will use the Yoneda embedding defined above,   
in order to define, for an $S$-category $T$ and an object $x \in T$,   
the comma $S$-category $T/x$ in a meaningful way.\\  
  
Let $T$ be an $S$-category in $\mathbb{U}$, and let us consider its (usual, enriched) Yoneda  
embedding  
$$\underline{h} : T \longrightarrow SPr(T):=SSet_{\mathbb{U}}^{T^{op}}.$$  
For any object $x \in Ob(T)$, we consider the comma category $SPr(T)/\underline{h}_{x}$,   
together with its natural   
induced model structure (i.e. the one created by the forgetful functor   
$SPr(T)/\underline{h}_{x} \rightarrow SPr(T)$, see \cite[p. $5$]{ho}).   
For any object $y \in Ob(T)$, and any morphism    
$u : \underline{h}_{y} \longrightarrow \underline{h}_{x}$, let   
$F_{u}\in SPr(T)/\underline{h}_{x}$ be a fibrant replacement   
of $u$. Since $u$ is already a cofibrant object in $SPr(T)/\underline{h}_{x}$   
(as we already observed in the proof of Proposition \ref{l1}), the object $F_{u}$ is then actually fibrant and cofibrant.   
  
\begin{df}\label{d8}  
The \emph{comma $S$-category $T/x$} is defined to be the full   
sub-$S$-category of $L(SPr(T)/\underline{h}_{x})$  
consisting of all objects $F_{u}$, for all $u$ of the form   
$u: \underline{h}_{y} \rightarrow \underline{h}_{x}$, $y \in Ob(T)$.   
\end{df}  
  
\medskip  
  
\noindent Note that since $T$ belongs to $\mathbb{U}$, so does the $S$-category $T/x$, for any object $x \in Ob(T)$.   
  
There exists a natural morphism in $\mathrm{Ho}(S-Cat_{\mathbb{V}})$  
$$T/x \longrightarrow L(SPr(T)/\underline{h}_{x}) \longrightarrow LSPr(T),$$  
where the morphism on the right is induced by the forgetful functor $SPr(T)/\underline{h}_{x} \longrightarrow SPr(T)$.   
One checks immediately that the essential image of this morphism is contained in the essential image of the  
Yoneda embedding $L\underline{h} : T \longrightarrow LSPr(T)$.   
Therefore, there exists a natural   
factorization in $\mathrm{Ho}(S-Cat_{\mathbb{V}})$  
$$\xymatrix{  
T/x \ar[rd]_-{j_{x}} \ar[rr] & &  LSPr(T) \\  
 & T \ar[ru]_-{L\underline{h}} &}$$  
As the inclusion functor $\mathrm{Ho}(S-Cat_{\mathbb{U}}) \longrightarrow   
\mathrm{Ho}(S-Cat_{\mathbb{V}})$ is fully faithful (see Appendix A), this  
gives a well defined morphism in $\mathrm{Ho}(S-Cat_{\mathbb{U}})$  
$$j_{x} : T/x \longrightarrow T.$$  
It is important to observe that the functor $\mathbb{R}(j_{x})_{!} :   
\mathrm{Ho}(SPr(T/x)) \longrightarrow \mathrm{Ho}(SPr(T))$, induced by $j_{x}$  
is such that $\mathbb{R}(j_{x})_{!}(*)\simeq \underline{h}_{x}$.  
  
Up to a natural equivalence of categories, the homotopy category   
$\mathrm{Ho}(T/x)$ has the following explicit description.   
For the sake of simplicity we will assume that $T$ is a fibrant $S$-category   
(i.e. all the simplicial sets $\underline{Hom}_{T}(x,y)$ of morphisms are fibrant).   
The objects of $\mathrm{Ho}(T/x)$ are simply pairs $(y,u)$, consisting of an object $y \in Ob(T)$ and a $0$-simplex  
$u \in \underline{Hom}_{T}(y,x)_{0}$ (i.e. a morphism $y \rightarrow x$ in the category $T_{0}$).   
  
Let us consider two objects $(y,u)$ and $(z,v)$, and a pair $(f,h)$, consisting of   
a $0$-simplex $f \in \underline{Hom}_{T}(y,z)$ and a $1$-simplex $h\in \underline{Hom}_{T}(y,x)^{\Delta^{1}}$  
such that   
$$\partial_{0}(h)=u \qquad \partial_{1}(h)=v\circ f.$$  
We may represent diagramatically this situation as:   
$$\xymatrix{  
y \ar[rr]^-{f} \ar[rdd]_-{u} & & z \ar[ldd]^-{v} \\  
 & h \Rightarrow & \\  
 & x & }$$  
Two such pairs $(f,h)$ and $(g,k)$ are defined to be equivalent if there exist  
a $1$-simplex $H \in \underline{Hom}_{T}(y,z)^{\Delta^{1}}$ and a $2$-simplex   
$G \in \underline{Hom}_{T}(y,x)^{\Delta^{2}}$ such that  
$$\partial_{0}(H)=f \qquad \partial_{1}(H)=g \qquad \partial_{0}(G)=h \qquad   
\partial_{1}(G)=k \qquad \partial_{2}(G)=v\circ H.$$  
The set of morphisms in $\mathrm{Ho}(T/x)$ from $(y,u)$ to $(z,v)$ is then the set of equivalences classes  
of such pairs $(f,h)$.   
In other words, the set of morphisms from $(y,u)$ to $(z,v)$ is the set of connected components  
of the homotopy fiber of $\;-\circ v : \underline{Hom}_{T}(y,z) \longrightarrow \underline{Hom}_{T}(y,x)$ at  
the point $u$.  
  
Let $(f,h) : (y,u) \longrightarrow (z,v)$ and $(g,k) : (z,v) \longrightarrow (t,w)$ be two morphisms  
in $\mathrm{Ho}(T/x)$.  The composition of $(f,h)$ and $(g,k)$ in   
$\mathrm{Ho}(T/x)$ is the class of $(g\circ f,k\dot h)$, where $k\dot h$  
is the concatenation of the $1$-simplices $h$ and $k\circ f$ in $\underline{Hom}_{T}(y,x)$.   
Pictorially, one composes the triangles as  
$$
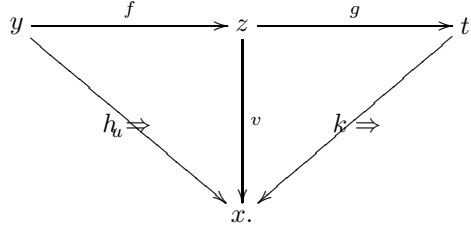
\xymatrix{  
y \ar[rr]^-{f} \ar[rrdd]_-{u} & & z \ar[dd]^-{v} \ar[rr]^-{g} & & t \ar[lldd]\\  
 & h \Rightarrow &  & k  \Rightarrow & \\  
 & & x. & & }$$  
As the  
concatenation of $1$-simplices is well defined, associative and unital up to homotopy, this gives a   
well defined, associative and unital composition of morphisms in $\mathrm{Ho}(T/x)$. 
  
Note that there is a natural projection $\mathrm{Ho}(T/x) \longrightarrow \mathrm{Ho}(T)/x$,   
which sends an object $(y,u)$ to   
the object $y$ together with the image of $u$ in   
$\pi_{0}(\underline{Hom}_{T}(y,x))=Hom_{\mathrm{Ho}(T)}(y,x)$. This   
functor is not an equivalence but is always full and essentially surjective. The composition functor  
$\mathrm{Ho}(T/x) \longrightarrow \mathrm{Ho}(T)/x \longrightarrow \mathrm{Ho}(T)$   
is isomorphic to the functor induced by the  
natural morphism $T/x \longrightarrow T$.   
  
\end{subsection}  
  
\end{section}  
  
\begin{section}{Stacks over $S$-sites}  
  
This section is devoted to the definition of the notions of $S$-topologies, $S$-sites and stacks over them. We start by  
defining \textit{$S$-topologies} on $S$-categories, generalizing the notion of Grothendieck topologies on   
usual categories and inducing an obvious notion of \textit{$S$-site}.   
For an $S$-site $T$, we define a notion of \textit{local equivalence}  
in the model category of pre-stacks $SPr(T)$, analogous to the notion   
of local isomorphism between presheaves on a given Grothendieck site.   
The first main result of this section is the existence of a model structure on $SPr(T)$, the  
\textit{local model structure}, whose equivalences are exactly the local equivalences. This model structure   
is called the \textit{model category of stacks}. To motivate this   
terminology we prove a criterion characterizing fibrant  
objects in the model category of stacks as objects satisfying a   
\textit{hyperdescent} property with respect to the given $S$-topology,   
which is a homotopy analog of the usual descent or sheaf condition.   
We also investigate functoriality properties (i.e. inverse and direct   
images functors) of the model categories of stacks, as well   
as the very useful notion of \textit{stack of morphisms} (i.e. internal $Hom$'s).   
  
The second main result of this section is a correspondence between $S$-topologies on an $S$-category $T$ and   
$t$-complete left Bousfield localizations of the model category of pre-stacks $SPr(T)$.   
Finally, we relate our definition of stacks  
over $S$-sites to the notion of \textit{model topos} due to C. Rezk, and we conclude from our previous results   
that almost all model topoi are equivalent to a model category of stacks over an $S$-site.   
  
\begin{subsection}{$S$-topologies and  $S$-sites} \label{stop}  
  
We refer to \cite[Exp. II]{sga4} or \cite{mm} for the   
definition of a Grothendieck topology and for the associated sheaf theory. \\  
  
\begin{df}\label{d9}  
An \emph{$S$-topology} on an $S$-category $T$ is a Grothendieck topology on the category $\mathrm{Ho}(T)$. An \emph{$S$-site} $(T,\tau)$ is the datum of an $S$-category $T$ together with an $S$-topology $\tau$ on $T$.  
\end{df}  
  
\begin{rmk}\label{infty}  
\begin{enumerate}  
\item \emph{It is important to remark that the notion of an $S$-topology on an $S$-category $T$ only depends  
on the isomorphism class of $T \in Ho(S-Cat)$, since equivalent $S$-categories have equivalent homotopy categories.}   
  
\item \emph{From the point of view of higher category theory, $S$-categories are models   
for $\infty$-categories in which all  
$i$-arrows are \textit{invertible} for all $i>1$. Therefore, if one tries to define the notion of a topology on   
this kind of higher categories, the stability axiom will imply that all   
$i$-morphisms should be automatically coverings for $i>1$.   
The datum of the topology should therefore only depends on isomorphism   
classes of $1$-morphisms, or, in other words, on the  
homotopy category. This might give a more conceptual explanation of Definition \ref{d9}.   
See also Remark \ref{lurie} for more on topologies on higher categories.}  
\end{enumerate}  
\end{rmk}  
  
\medskip  
  
Let $T \in S-Cat_{\mathbb{U}}$ be a $\mathbb{U}$-small $S$-category and $SPr(T)$ its  
model category of pre-stacks. Given any pre-stack $F\in SPr(T)$,   
one can consider its associated presheaf of  
connected components  
$$\begin{array}{ccc}  
T^{op} & \longrightarrow & Set_{\mathbb{U}} \\  
x & \mapsto & \pi_{0}(F(x)).  
\end{array}$$  
  
The universal property of the homotopy category of $T^{op}$ implies that there exists a unique  
factorization  
$$\xymatrix{  
T^{op} \ar[r] \ar[d] & Set_{\mathbb{U}} \\  
\mathrm{Ho}(T)^{op} \ar[ru]_-{\pi_{0}^{pr}(F)} & }$$  
  
The construction $F \mapsto \pi_{0}^{pr}(F)$, being obviously functorial in $F$, induces a well defined functor  
$SPr(T) \longrightarrow Set_{\mathbb{U}}^{\mathrm{Ho}(T)^{op}}$; but,   
since equivalences in $SPr(T)$ are defined objectwise, this induces a functor    
$$\pi_{0}^{pr}(-) : \mathrm{Ho}(SPr(T)) \longrightarrow Set_{\mathbb{U}}^{\mathrm{Ho}(T)^{op}}.$$  
  
\begin{df}\label{d10}  
Let $(T,\tau)$ be a $\mathbb{U}$-small $S$-site.   
\begin{enumerate}  
\item For any object $F \in SPr(T)$, the sheaf associated to the presheaf $\pi_{0}^{pr}(F)$  
is denoted by $\pi_{0}^{\tau}(F)$ (or $\pi_{0}(F)$ if the topology $\tau$ is   
clear from the context). It is a sheaf on the site $(\mathrm{Ho}(T),\tau)$, and is called the  
\emph{sheaf of connected components of} $F$.   
\item   
A morphism $F \longrightarrow G$ in $\mathrm{Ho}(SPr(T))$  
is called a \emph{$\tau$-covering} (or just a \emph{covering}   
if the topology $\tau$ is clear from the context) if the induced morphism   
$\pi_{0}^{\tau}(F) \longrightarrow \pi_{0}^{\tau}(G)$  
is an epimorphism of sheaves.  
  
\item A morphism $F \longrightarrow G$ in $SPr(T)$  
is called a \emph{$\tau$-covering} (or just a \emph{covering} if the topology $\tau$ is unambiguous) if its image  
by the natural functor $SPr(T) \longrightarrow \mathrm{Ho}(SPr(T))$ is a $\tau$-covering as defined in the previous item.  
\end{enumerate}  
\end{df}  
  
Clearly, for two objects $x$ and $y$ in $T$,   
any morphism $x \longrightarrow y$ such that the sieve generated by   
its image in $\mathrm{Ho}(T)$ is a covering sieve of $y$, induces  
a covering $\underline{h}_{x} \longrightarrow \underline{h}_{y}$.  
  
More generally, one has the following characterization of coverings as   
\textit{homotopy locally surjective} morphisms. This is the homotopy analog of the notion of   
epimorphism of stacks (see for example \cite[\S 1]{lm}), where one requires that any object in the target is   
locally isomorphic to the image of an object in the source.   
  
\begin{prop} \label{hls}   
A morphism $f : F \longrightarrow G$ in   
$SPr(T)$ is a covering if it has the following \emph{homotopy local surjectivity} property. For any object  
$x\in Ob(T)$, and any morphism in $\mathrm{Ho}(SPr(T))$,  
$\underline{h}_{x} \longrightarrow G$, there exists a covering sieve $R$ of $x$ in $\mathrm{Ho}(T)$, such that for   
any morphism $u \rightarrow x$ in $R$ there is a commutative diagram in $\mathrm{Ho}(SPr(T))$  
$$\xymatrix{  
F \ar[r] & G \\  
\underline{h}_{u} \ar[u] \ar[r] & \underline{h}_{x}. \ar[u]}$$  
In other words, $f$ is a covering if and only if \emph{any object   
of $G$ over $x$ lifts locally and up to homotopy to an object of $F$}.   
\end{prop}  
  
\textit{Proof:} First of all, let us observe that both the definition of a   
covering and the \textit{homotopy local surjectivity property} hold true for the   
given $f : F \rightarrow G$ if and only if they hold true for   
$RF \rightarrow RG$, where $R(-)$ is a fibrant replacement functor in $SPr(T)$.   
Therefore, we may suppose both $F$ and $G$ fibrant.  
Now, by \cite[III.$7$, Cor. $6$]{mm}, $f : F \rightarrow G$ is a covering iff the   
induced map of presheaves $\pi_{0}^{pr}(F) \rightarrow \pi_{0}^{pr}(G)$ is locally   
surjective. But, by Yoneda $\pi_{0}^{pr}(H)(y) \simeq \pi_{0}(\underline{Hom}_{SPr(T)}(\underline{h}_{y},H))$,   
for any $H \in SPr(T)$ and any object $y$ in $T$. Since $F$ and $G$ are fibrant,   
we then have $\pi_{0}^{pr}(F)(y) \simeq Hom_{\mathrm{Ho}(SPr(T))}(\underline{h}_{y},F)$   
and $\pi_{0}^{pr}(G)(y) \simeq Hom_{\mathrm{Ho}(SPr(T))}(\underline{h}_{y},G)$,   
for any object $y$ in $T$. But then, the local surjectivity of $\pi_{0}^{pr}(F)   
\rightarrow \pi_{0}^{pr}(G)$ exactly translates to the \textit{homotopy local surjectivity property}  
in the proposition and we conclude.  
\hfill \textbf{$\Box$}  
  
\begin{rmk}\label{locstric}  
\emph{If the morphism $f$ of Proposition \ref{hls} is an objectwise fibration (i.e.   
for any $x \in T$, the morphism $F(x) \longrightarrow G(x)$  
is a fibration of simplicial sets), then the homotopy local surjectivity property  
implies the local surjectivity property. This means that the diagrams  
$$\xymatrix{  
F \ar[r] & G \\  
\underline{h}_{u} \ar[u] \ar[r] & \underline{h}_{x} \ar[u]}$$  
of Proposition \ref{hls} can be chosen to be commutative in $SPr(T)$, and not only in $\mathrm{Ho}(SPr(T))$. }  
\end{rmk}  
  
From this characterization one concludes easily that coverings have the following stability properties.   
  
\begin{cor}\label{p1}  
\begin{enumerate}  
  
\item A morphism in $SPr(T)$ which is a composition of coverings is a covering.  
  
\item Let   
$$\xymatrix{  
F' \ar[r]^{f'} \ar[d] & G' \ar[d] \\  
F \ar[r]_{f} & G }$$  
be a homotopy cartesian diagram in $SPr(T)$. If $f$ is a covering so is $f'$.  
  
\item Let $\xymatrix{F \ar[r]^-{u} & G \ar[r]^-{v} & H}$ be two morphisms in $SPr(T)$. If   
the morphism $v\circ u$ is a covering then so is $v$.  
  
\item Let   
$$\xymatrix{  
F' \ar[r]^{f'} \ar[d] & G' \ar[d]^{p} \\  
F \ar[r]_{f} & G }$$  
be a homotopy cartesian diagram in $SPr(T)$. If $p$ and $f'$ are coverings then so is $f$.   
  
\end{enumerate}  
\end{cor}  
  
\textit{Proof:}   
Properties $(1)$ and $(3)$ follows immediately from Proposition \ref{hls},   
and $(4)$ follows from $(3)$. It remains to prove $(2)$.  
Let us $f$ and $f'$ be as in $(2)$ and let us consider a diagram  
$$\xymatrix{  
& \underline{h}_{x} \ar[d] \\  
F' \ar[d] \ar[r]^-{f'} & G' \ar[d] \\  
F \ar[r]^-{f} & G.}$$  
As $f$ is a covering, there exists a covering sieve $R$   
over $x \in \mathrm{Ho}(T)$, such that for any $u\rightarrow x$ in $R$,  
one has a commutative diagram  
$$\xymatrix{  
\underline{h}_{u} \ar[d] \ar[r] & \underline{h}_{x} \ar[d]  \\  
F \ar[r]^-{f} & G.}$$  
By the universal property of homotopy fibered products, the morphisms  
$\underline{h}_{u} \longrightarrow F$ and $\underline{h}_{u} \longrightarrow \underline{h}_{x} \longrightarrow G'$  
are the two projections of a (non unique) morphism $\underline{h}_{u} \longrightarrow F'$.   
This gives, for all $u \rightarrow x$, the required liftings  
$$\xymatrix{  
\underline{h}_{u} \ar[d] \ar[r] & \underline{h}_{x} \ar[d] \\  
F' \ar[r]^-{f'} & G'.}$$  

\hfill \textbf{$\Box$} 
  
\end{subsection}  
  
\begin{subsection}{Simplicial objects and hypercovers}\label{sobj}  
  
Let us now consider $sSPr(T):=SPr(T)^{\Delta^{op}}$, the category of simplicial objects in $SPr(T)$. Its objects  
will be denoted as   
$$\begin{array}{cccc}  
F_{*} : & \Delta^{op} & \longrightarrow & SPr(T) \\  
 & [m] & \mapsto & F_{m}.  
\end{array}$$  
  
As the category $SPr(T)$ has all kind of limits and colimits indexed in $\mathbb{U}$, the category $sSPr(T)$   
has a natural structure of tensored and co-tensored category over $SSet_{\mathbb{U}}$ (see \cite[Ch. II, Thm. 2.5]{gj}).   
The external product of $F_{*} \in sSPr(T)$ by $A \in SSet_{\mathbb{U}}$, denoted by   
$\underline{A}\otimes F_{*}$, is the simplicial object in $SPr(T)$ defined by  
$$\begin{array}{cccc}  
\underline{A}\otimes F_{*} : & \Delta^{op} & \longrightarrow & SPr(T) \\  
& [n] & \mapsto & \coprod_{A_{n}}F_{n}.  
\end{array}$$  
The exponential (or co-tensor) of $F_{*}$ by $A$, is denoted by $F_{*}^{\underline{A}}$ and is determined by the usual   
adjunction isomorphism  
$$Hom(\underline{A}\otimes F_{*},G_{*}) \simeq Hom(F_{*},G_{*}^{\underline{A}}).$$  
  
\noindent \textbf{Notation:} We will denote by $F_{*}^{A} \in SPr(T)$   
the $0$-th level of the simplicial object $F_{*}^{\underline{A}} \in sSPr(T)$. \\  
  
\noindent Explicitly, the object $F_{*}^{A}$ is the \textit{end} of the functor  
$$\begin{array}{ccc}  
\Delta^{op}\times \Delta & \longrightarrow & SPr(T) \\  
([n],[m]) & \mapsto & \prod_{A_{m}}F_{n}.  
\end{array}$$  
  
One checks immediately that for any $F_{*} \in sSPr(T)$, one has a natural   
isomorphism $F_{*}^{\Delta^{n}}\simeq F_{n}$. \\  
  
We endow the category $sSPr(T)$ with its Reedy model structure (see \cite[Thm. 5.2.5]{ho}).   
The equivalences in $sSPr(T)$ are  
the morphisms $F_{*} \longrightarrow G_{*}$ such that, for any $n$, the induced morphism  
$F_{n} \longrightarrow G_{n}$ is an equivalence in $SPr(T)$. The fibrations are the   
morphisms $F_{*} \longrightarrow G_{*}$  
such that, for any $[n] \in \Delta$, the induced morphism  
$$F_{*}^{\Delta^{n}} \simeq F_{n} \longrightarrow F_{*}^{\partial \Delta^{n}}  
\times_{G_{*}^{\partial \Delta^{n}}}G_{*}^{\Delta^{n}}$$  
is a fibration in $SPr(T)$.   
  
Given any simplicial set $A \in SSet_{\mathbb{U}}$, the functor  
$$\begin{array}{ccc}  
sSPr(T) & \longrightarrow & SPr(T) \\  
 F_{*} & \mapsto & F_{*}^{A}  
\end{array}$$  
is a right Quillen functor for the Reedy model structure on $sSPr(T)$ (\cite[Prop. 5.4.1]{ho}). Its  
right derived functor will be denoted by  
$$\begin{array}{ccc}  
\mathrm{Ho}(sSPr(T)) & \longrightarrow & \mathrm{Ho}(SPr(T)) \\  
F_{*} & \mapsto & F_{*}^{\mathbb{R}A}.  
\end{array}$$  
  
For any object $F \in SPr(T)$, one can consider the constant simplicial object $c(F)_{*} \in sSPr(T)$  
defined by $c(F)_{n}:=F$ for all $n$. One the other hand, one can consider  
$$\begin{array}{cccc}  
(RF)^{\Delta^{*}}: & \Delta^{op} & \longrightarrow & SPr(T) \\  
& [n] & \mapsto & (RF)^{\Delta^{n}},  
\end{array}$$  
where $RF$ is a fibrant model for $F$ in $SPr(T)$, and  
$(RF)^{\Delta^{n}}$ is the exponential object defined using the simplicial structure on $SPr(T)$.   
The object $(RF)^{\Delta^{*}}$ is a fibrant replacement of $c_{*}(F)$ in $sSPr(T)$.   
Furthermore, for any object $G \in SPr(T)$ and $A \in SSet_{\mathbb{U}}$, there exists a natural isomorphism in $SPr(T)$  
$$(G^{\Delta^{*}})^{A}\simeq G^{A}.$$  
This induces a natural isomorphism in $\mathrm{Ho}(SPr(T))$  
$$(c(F)_{*})^{\mathbb{R}A}\simeq ((RF)^{\Delta^{*}})^{A}\simeq (RF)^{A}.$$  
However, we remark that   
$c(F)_{*}^{A}$ is not isomorphic to $F^{A}$ as an object in $SPr(T)$. \\  
  
\noindent \textbf{Notation:} For any $F \in SPr(T)$ and $A \in SSet_{\mathbb{U}}$, we will simply denote by  
$F^{\mathbb{R}A} \in \mathrm{Ho}(SPr(T))$ the object $c(F)_{*}^{\mathbb{R}A}\simeq (RF)^{A}$. \\  
  
We let $\Delta_{\leq n}$ be the full subcategory of $\Delta$ consisting of objects  
$[p]$ with $p\leq n$, and denote by $s_{n}SPr(T)$ the category of functors  
$\Delta^{op}_{\leq n} \longrightarrow SPr(T)$. The natural inclusion   
$i_{n} : \Delta_{\leq n} \rightarrow \Delta$ induces  
a restriction functor  
$$i_{n}^{*} : sSPr(T) \longrightarrow s_{n}SPr(T)$$  
which has a right adjoint $(i_{n})_{*} : s_{n}SPr(T) \longrightarrow sSPr(T)$,   
as well as a left adjoint   
$(i_{n})_{!} : s_{n}SPr(T) \longrightarrow sSPr(T).$  
The two adjunction morphisms induce isomorphisms   
$i_{n}^{*}(i_{n})_{*}\simeq \mathrm{Id}$ and $i_{n}^{*}(i_{n})_{!}\simeq \mathrm{Id}$:  
therefore both functors $(i_{n})_{*}$ and $(i_{n})_{!}$ are fully faithful.   
  
\begin{df}\label{d11}  
Let $F_{*} \in sSPr(T)$ and $n\geq 0$.   
\begin{enumerate}  
\item One defines the $n$\emph{-th skeleton} and $n$\emph{-th coskeleton} of $F_{*}$ as  
$$Sk_{n}(F_{*}):=(i_{n})_{!}i_{n}^{*}(F_{*}) \qquad Cosk_{n}(F_{*}):=(i_{n})_{*}i_{n}^{*}(F_{*}).$$  
  
\item The simplicial object $F_{*}$ is called $n$\emph{-bounded} if the adjunction morphism   
$F_{*}\longrightarrow Cosk_{n}(F_{*})$ is an isomorphism.  
  
\end{enumerate}  
\end{df}  
  
\medskip   
  
It is important to note that $F_{*}$, $Cosk_{n}(F_{*})$ and $Sk_{n}(F_{*})$ all coincide in degrees $\leq n$  
$$i_{n}^{*}(F_{*})\simeq i_{n}^{*}(Cosk_{n}F_{*})\simeq i_{n}^{*}(Sk_{n}F_{*}).$$  
The adjunctions $(i_{n}^{*}, (i_{n})_{*})$ and $((i_{n})_{!},i_{n}^{*})$ induce a natural adjunction isomorphism  
$$Hom(Sk_{n}(F_{*}),G_{*})\simeq Hom(F_{*},Cosk_{n}(G_{*})),$$  
for any $F_{*}$ and $G_{*}$ in $sSPr(T)$ and any $n \geq 0$. As a special case, for any $A \in SSet_{\mathbb{U}}$, one  
has an isomorphism in $SPr(T)$  
$$F_{*}^{Sk_{n}A}\simeq (Cosk_{n}F_{*})^{A}.$$  
As $Sk_{n}\Delta^{n+1}=\partial \Delta^{n+1}$, one gets natural isomorphisms  
  
\begin{equation}\label{num}  
F_{*}^{\partial \Delta^{n+1}}\simeq Cosk_{n}(F_{*})_{n+1}.   
\end{equation}  
  
\begin{lem}\label{l2}  
The functor $Cosk_{n} : sSPr(T) \longrightarrow sSPr(T)$ is a right Quillen functor for the Reedy model structure  
on $sSPr(T)$.  
\end{lem}  
  
\textit{Proof:} By adjunction, for any integer $p$ with $p\leq n$, one has  
$$(Cosk_{n}(F_{*}))^{\partial\Delta^{p}}\simeq F_{*}^{\partial \Delta^{p}}   
\qquad (Cosk_{n}(F_{*}))^{\Delta^{p}}\simeq F_{*}^{\Delta^{p}},$$  
while, for $p>n+1$, one has  
$$(Cosk_{n}(F_{*}))^{\partial \Delta^{p}}\simeq (Cosk_{n}(F_{*}))^{\Delta^{p}}.$$  
Finally, for $p=n+1$ one has  
$$(Cosk_{n}(F_{*}))^{\partial\Delta^{n+1}}\simeq F_{*}^{\partial \Delta^{n+1}}   
\qquad (Cosk_{n}(F_{*}))^{\Delta^{n+1}}\simeq  
F_{*}^{\partial \Delta^{n+1}}.$$  
Using these formulas and the definition of Reedy fibrations in $sSPr(T)$ one checks immediately that the functor  
$Cosk_{n}$ preserves fibrations and trivial fibrations. As it is a right adjoint (its left adjoint being $Sk_{n}$),  
this implies that $Cosk_{n}$ is a right Quillen functor. \hfill $\Box$ \\  
  
The previous lemma allows us to consider the right derived version of the coskeleton functor  
$$\mathbb{R}Cosk_{n}  : \mathrm{Ho}(sSPr(T)) \longrightarrow \mathrm{Ho}(sSPr(T)).$$  
It comes with a natural morphism $\mathrm{Id}_{\mathrm{Ho}(sSPr(T))}   
\longrightarrow \mathbb{R}Cosk_{n}(F)$, induced by the adjunction   
morphism $\mathrm{Id}_{sSPr(T)} \longrightarrow (i_{n})_{*}i_{n}^{*}$.   
There exist obvious \textit{relative} notions of the functors $Sk_{n}$   
and $Cosk_{n}$ whose formulations are left to the reader. Let us  
only mention that   
the relative derived coskeleton of a morphism $F_{*} \longrightarrow G_{*}$   
in $sSPr(T)$ may be defined by the following   
homotopy cartesian square in $SPr(T)$  
$$\xymatrix{  
\mathbb{R}Cosk_{n}(F_{*}/G_{*}) \ar[r] \ar[d] & G_{*} \ar[d] \\  
\mathbb{R}Cosk_{n}(F_{*}) \ar[r] & \mathbb{R}Cosk_{n}(G_{*}).  
}$$  
  
\medskip  
  
The functor $\mathbb{R}Cosk_{0}(-/c(G)_{*})$, relative to a constant diagram   
$c(G)_{*}$, where $G \in SPr(T)$, has the following interpretation in terms of \textit{derived nerves}.   
For any morphism $F_{*} \longrightarrow c_{*}(G)$ in $sSPr(T)$, with $c_{*}(G)$ the constant simplicial   
diagram with value $G$, we consider the induced morphism $f : F_{0}   
\longrightarrow G$ in $\mathrm{Ho}(SPr(T))$. Let us represent this morphism   
by a fibration in $SPr(T)$, and let us consider its usual nerve $N(f)$:    
$$\begin{array}{cccc}  
N(f) : & \Delta^{op} & \longrightarrow & SPr(T) \\  
& [n] & \mapsto & \underbrace{F_{0}\times_{G}F_{0}\times_{G} \dots \times_{G}F_{0}}_{\textrm{n times}}.  
\end{array}  
$$  
The nerve $N(f)$ is naturally augmented over $G$, and therefore is an object  
of $sSPr(T)/c_{*}(G)$.  
Then, there is a natural isomorphism in $\mathrm{Ho}(sSPr(T)/c_{*}(G))$  
$$\mathbb{R}Cosk_{0}(F_{*}/c_{*}(G))\simeq N(f).$$  
   
\begin{df}\label{d12}  
Let $(T,\tau)$ be a $\mathbb{U}$-small $S$-site.  
\begin{enumerate}  
\item A morphism in $sSPr(T)$  
$$F_{*} \longrightarrow G_{*}$$  
is called a $\tau$\emph{-hypercover} (or just a \emph{hypercover} if   
the topology $\tau$ is unambiguous) if for any $n$, the induced morphism   
$$F_{*}^{\mathbb{R}\Delta^{n}}\simeq F_{n} \longrightarrow  
F_{*}^{\mathbb{R}\partial \Delta^{n}}\times^{h}_{G_{*}^{\mathbb{R}  
\partial\Delta^{n}}}G_{*}^{\mathbb{R}\Delta^{n}}$$  
is a covering in $\mathrm{Ho}(SPr(T))$ (see Definition \ref{d10} $(2)$).  
  
\item A morphism in $\mathrm{Ho}(sSPr(T))$  
$$F_{*} \longrightarrow G_{*}$$  
is called a $\tau$\emph{-hypercover} (or just a \emph{hypercover} if the topology   
$\tau$ is unambiguous) if one of its representatives in $sSPr(T)$ is a $\tau$-hypercover.  
\end{enumerate}  
\end{df}  
  
Using the isomorphisms (\ref{num}), Definition \ref{d12} may also be stated  
as follows. A morphism $f : F_{*} \longrightarrow G_{*}$ is a $\tau$-hypercover \textit{if and only if}  
for any $n\geq 0$ the induced morphism  
$$F_{n} \longrightarrow \mathbb{R}Cosk_{n-1}(F_{*}/G_{*})_{n}$$  
is a covering in $\mathrm{Ho}(SPr(T))$.   
   
Note also that in Definition \ref{d12} $(2)$, if one of the representatives of $f$ is a hypercover, then so is  
any representative. Being a hypercover is therefore a property of morphisms in $\mathrm{Ho}(sSPr(T))$.  
   
\end{subsection}  
  
\begin{subsection}{Local equivalences}\label{loceq}  
  
Throughout this subsection, we fix a $\mathbb{U}$-small $S$-site $(T,\tau)$.  
  
Let $x$ be an object in $T$. The topology on $\mathrm{Ho}(T)$ induces a natural   
topology on the comma category $\mathrm{Ho}(T)/x$.  
We define a Grothendieck topology on $\mathrm{Ho}(T/x)$ by   
pulling back the topology of $\mathrm{Ho}(T)/x$ through the   
natural projection $\mathrm{Ho}(T/x) \longrightarrow \mathrm{Ho}(T)/x$.   
By this, we mean that a sieve $S$ over an object $y \in Ho(T/x)$, is   
defined to be a covering sieve if and only if (the sieve generated by) its image   
in $Ho(T)$ is a $\tau$-covering sieve of the object $y \in Ho(T)/x$. The reader   
will check easily that this indeed defines a topology on $Ho(T/x)$, and therefore an $S$-topology on $T/x$.  
This topology will still be denoted by $\tau$.  
  
\begin{df}\label{d13}  
The $S$-site $(T/x,\tau)$ will be called the \emph{comma $S$-site} of $(T,\tau)$ over $x$.  
\end{df}  
  
\medskip  
  
Let $F \in SPr(T)$, $x \in Ob(T)$ and $s \in \pi_{0}(F(x))$   
be represented by a morphism $s : \underline{h}_{x} \longrightarrow F$ in $\mathrm{Ho}(SPr(T))$ (see \ref{l1}).    
By pulling-back this morphism along the natural morphism $j_{x} :   
T/x \longrightarrow T$, one gets a morphism in $\mathrm{Ho}(SPr(T/x))$  
$$s : j_{x}^{*}(\underline{h}_{x}) \longrightarrow j_{x}^{*}(F).$$  
By definition of the comma category $T/x$, it is immediate that   
$j_{x}^{*}(\underline{h}_{x})$ has a natural global point  
$* \longrightarrow j_{x}^{*}(\underline{h}_{x})$ in $\mathrm{Ho}(SPr(T/x))$.   
Note that the morphism $* \longrightarrow j_{x}^{*}(\underline{h}_{x})$  
is also induced by adjunction from the identity of $\underline{h}_{x}  
\simeq \mathbb{R}(j_{x})_{!}(*)$.  
Therefore we obtain a global point of $j_{x}^{*}(F)$   
$$s : * \longrightarrow j_{x}^{*}(\underline{h}_{x})\longrightarrow j_{x}^{*}(F).$$   
  
\begin{df}\label{d14}  
Let $F \in SPr(T)$ and $x \in Ob(T)$.  
\begin{enumerate}  
\item   
For any integer $n>0$, the sheaf $\pi_{n}(F,s)$ is defined as   
$$\pi_{n}(F,s):=\pi_{0}(j_{x}^{*}(F)^{\mathbb{R}\Delta^{n}}  
\times_{j_{x}^{*}(F)^{\mathbb{R}\partial \Delta^{n}}}\;*).$$  
It is a sheaf on the site $(\mathrm{Ho}(T/x),\tau)$ called the   
$n$\emph{-th homotopy sheaf} of $F$ \emph{pointed} at $s$.   
  
\item A morphism $f : F \longrightarrow G$ in $SPr(T)$ is called a   
$\pi_{*}$\emph{-equivalence} or, equivalently, a  
\emph{local equivalence}, if the following two conditions are satisfied:  
\begin{enumerate}  
\item The induced morphism   
$\pi_{0}(F) \longrightarrow \pi_{0}(G)$  
is an isomorphism of sheaves on $\mathrm{Ho}(T)$.  
\item For any object $x \in Ob(T)$, any section $s \in \pi_{0}(F(x))$ and any integer $n>0$, the induced morphism  
$\pi_{n}(F,s) \longrightarrow \pi_{n}(G,f(s))$  
is an isomorphism of sheaves on $\mathrm{Ho}(T/x)$.   
\end{enumerate}  
  
\item A morphism in $\mathrm{Ho}(SPr(T))$ is a $\pi_{*}$-equivalence if one of its representatives in $SPr(T)$  
is a $\pi_{*}$-equivalence.  
  
\end{enumerate}  
\end{df}  
  
\medskip  
  
Obviously an equivalence in the model category $SPr(T)$ is always a $\pi_{*}$-equivalence for any topology $\tau$ on $T$.   
Indeed, an equivalence in $SPr(T)$ induces isomorphisms between the homotopy presheaves which are the homotopy   
sheaves for the trivial topology.   
  
Note also that in Definition \ref{d14} $(3)$, if a representative of $f$ is a $\pi_{*}$-equivalence then so is  
any of its representatives. Therefore, being a $\pi_{*}$-equivalence is actually a property of morphisms  
in $\mathrm{Ho}(SPr(T))$. 
   
The following characterization of $\pi_{*}$-equivalences is interesting as it does not involve any base point.   
  
\begin{lem}\label{l3}  
A morphism $f : F \longrightarrow G$ in $SPr(T)$ is a $\pi_{*}$-equivalence if and only if  
for any $n\geq 0$ the induced morphism   
$$F^{\mathbb{R}\Delta^{n}} \longrightarrow F^{\mathbb{R}\partial \Delta^{n}}  
\times^{h}_{G^{\mathbb{R}\partial \Delta^{n}}}G^{\mathbb{R}\Delta^{n}}$$  
is a covering.   
  
In other words, $f : F \longrightarrow G$ is a $\pi_{*}$-equivalence if and only if it is a $\tau$-hypercover   
when considered as a morphism of constant simplicial objects in $SPr(T)$.  
\end{lem}  
  
\textit{Proof:} Without loss of generality we can assume that $f$ is a fibration   
between fibrant objects in the model category $SPr(T)$. This means that for any $x \in Ob(T)$,   
the induced morphism $f : F(x) \longrightarrow G(x)$  
is a fibration between fibrant simplicial sets. In particular, the morphism  
$$F^{\mathbb{R}\Delta^{n}} \longrightarrow   
F^{\mathbb{R}\partial \Delta^{n}}\times^{h}_{G^{\mathbb{R}\partial \Delta^{n}}}G^{\mathbb{R}\Delta^{n}}$$  
in $\mathrm{Ho}(SPr(T))$ is represented by the morphism in $SPr(T)$  
$$F^{\Delta^{n}} \longrightarrow F^{\partial \Delta^{n}}\times_{G^{\partial \Delta^{n}}}G^{\Delta^{n}}.$$  
This morphism is furthermore an objectwise fibration, and therefore the local lifting property of  
$\tau$-coverings (see Prop. \ref{hls}) holds not only in   
$\mathrm{Ho}(SPr(T))$ but in $SPr(T)$ (see Remark \ref{locstric}).   
Hence, $f$ is a hypercover if and only if it  
satisfies the following local lifting property.   
  
For any $x \in \mathrm{Ho}(T)$, and any morphism in $SPr(T)$  
$$\underline{h}_{x} \longrightarrow   
F^{\partial \Delta^{n}}\times_{G^{\partial \Delta^{n}}}G^{\Delta^{n}},$$  
there exists a covering sieve $R$ of $x$ and, for any $u \rightarrow x$ in $R$, a commutative  
diagram in $SPr(T)$  
$$\xymatrix{  
F^{\Delta^{n}} \ar[r] & F^{\partial \Delta^{n}}  
\times_{G^{\partial \Delta^{n}}}G^{\Delta^{n}} \\  
\underline{h}_{u} \ar[u] \ar[r] & \underline{h}_{x}. \ar[u] }$$  
  
By adjunction, this is equivalent to the following condition.  
For any object $x \in Ob(T)$ and any commutative diagram in $SSet_{\mathbb{U}}$  
$$\xymatrix{  
F(x) \ar[r] & G(x) \\  
\partial \Delta^{n} \ar[u] \ar[r] & \Delta^{n} \ar[u]}$$  
there exists a covering sieve $R$ of $x$ in $\mathrm{Ho}(T)$   
such that for any morphism $u \rightarrow x$ in $T$, whose image  
belongs to $R$, there is a commutative diagram in $SSet_{\mathbb{U}}$  
$$\xymatrix{  
F(u) \ar[r] & G(u) \\  
F(x) \ar[r] \ar[u] & G(x) \ar[u] \\  
\partial \Delta^{n} \ar[u] \ar[r] & \ar[luu] \Delta^{n} \ar[u].}$$  
By definition of the homotopy sheaves, this   
last condition is easily seen to be equivalent to being a   
$\pi_{*}$-equivalence (the details are left to the reader, who might  
also wish to consult \cite[Thm. $1.12$]{ja}). \hfill \textbf{$\Box$} 
  
\begin{cor}\label{c2}  
Let $f : F \longrightarrow G$ be a morphism in $SPr(T)$ and   
$G' \longrightarrow G$ be a covering. Then, if the induced morphism  
$$f' : F\times_{G}^{h}G' \longrightarrow G'$$  
is a $\pi_{*}$-equivalence, then so is $f$.  
\end{cor}  
  
\textit{Proof:} Apply Lemma \ref{l3} and Proposition \ref{p1} $(2)$. \hfill \textbf{$\Box$} 
  
\begin{cor}\label{c3}  
Let $f : F \longrightarrow G$ be a $\pi_{*}$-equivalence in $SPr(T)$ and   
$G' \longrightarrow G$ be an objectwise fibration. Then, the induced morphism  
$$f' : F\times_{G}G' \longrightarrow G'$$  
is a $\pi_{*}$-equivalence.  
\end{cor}  
  
\textit{Proof:} This follows from Corollary \ref{c2} since $SPr(T)$   
is a proper model category. \hfill \textbf{$\Box$} \\  
  
Let $x$ be an object in $T$ and $f:F \rightarrow G$ be a morphism in $\mathrm{Ho}(SPr(T))$.   
For any morphism $s : \underline{h}_{x} \longrightarrow G$ in $\mathrm{Ho}(SPr(T))$, let   
us define $F_{s} \in \mathrm{Ho}(SPr(T/x))$ by the following homotopy cartesian square in $SPr(T/x)$  
$$\xymatrix{  
j_{x}^{*}(F) \ar[r]^-{j_{x}^{*}(f)} & j_{x}^{*}(G) \\  
F_{s} \ar[u] \ar[r] & \mathrm{*} \ar[u]}$$  
where the morphism $* \longrightarrow j_{x}^{*}(G)$ is adjoint to the   
morphism $s : \mathbb{R}(j_{x})_{!}(*)\simeq \underline{h}_{x} \longrightarrow G$.   
  
\begin{cor}\label{c3''}  
Let $f : F \longrightarrow G$ be a morphism in $SPr(T)$. With the same notations as above, the morphism $f$ is a  
$\pi_{*}$-equivalence if and only for any $s : \underline{h}_{x} \longrightarrow G$ in $\mathrm{Ho}(SPr(T))$,   
the induced morphism $F_{s} \longrightarrow *$ is a $\pi_{*}$-equivalence in $\mathrm{Ho}(SPr(T/x))$.  
\end{cor}  
  
\textit{Proof:} By Lemma \ref{l3} it is enough to show that the morphism $f$ is a covering if and only if  
all the $F_{s} \longrightarrow *$ are coverings in $\mathrm{Ho}(SPr(T/x))$. The \textit{only if} part follows  
from Proposition \ref{p1} $(2)$, so it is enough to show that if all the $F_{s} \longrightarrow *$ are  
coverings then $f$ is a covering.   
  
Given $s : \underline{h}_{x} \longrightarrow G$ in $\mathrm{Ho}(SPr(T))$, let us prove that it lifts locally  
to $F$. By adjunction, $s$ corresponds to a morphism $* \longrightarrow j_{x}^{*}(G)$. As the corresponding  
morphism $F_{s} \longrightarrow *$ is a covering, there exists a covering sieve $R$ of $*$ in $\mathrm{Ho}(T/x)$  
and, for each $u \rightarrow *$ in $R$, a commutative diagram in $\mathrm{Ho}(SPr(T/x))$  
$$\xymatrix{  
j_{x}^{*}(F) \ar[r] & j_{x}^{*}(G) \\  
\underline{h}_{u} \ar[r] \ar[u] & \textrm{\large{*}} \ar[u]}$$  
By adjunction, this commutative diagram induces a commutative diagram in $\mathrm{Ho}(SPr(T))$  
$$\xymatrix{  
F \ar[r] & G \\  
\mathbb{R}(j_{x})_{!}(\underline{h}_{u}) \ar[r] \ar[u] & \underline{h}_{x} \ar[u]}$$  
But $\mathbb{R}(j_{x})_{!}(\underline{h}_{u})\simeq \underline{h}_{j_{x}(u)}$,   
and by definition of the induced topology on   
$\mathrm{Ho}(T/x)$, the morphisms in $(j_{x})(R)$ form a covering sieve of $x$.   
Therefore, the commutative diagram above   
shows that the morphism $s$ lifts locally to $F$.  \hfill \textbf{$\Box$} \\  
  
We end this paragraph by describing the behaviour of $\pi_{*}$-equivalences under homotopy push-outs.  
  
\begin{prop}\label{c3'}  
Let $f : F \longrightarrow G$ be a $\pi_{*}$-equivalence in $SPr(T)$ and   
$F \longrightarrow F'$ be an objectwise cofibration (i.e. a monomorphism). Then, the induced morphism  
$$f' : F' \longrightarrow F'\coprod_{F}G$$  
is a $\pi_{*}$-equivalence.  
\end{prop}  
  
\textit{Proof:} It is essentially the same proof as that of \cite[Prop. $2.2$]{ja}.  
\hfill \textbf{$\Box$}  
  
\end{subsection}  
  
\begin{subsection}{The local model structure}  
  
Throughout this subsection, we fix a $\mathbb{U}$-small $S$-site $(T,\tau)$. 
  
The main purpose of this paragraph is to prove the  
following theorem which is a generalization of the existence of the local projective model structure  
on the category of simplicial presheaves on a Grothendieck site   
(see for example \cite{bl} and \cite[\S $5$]{sh}). The proof   
we present here is based on some arguments found in \cite{s2},   
\cite{sh} and \cite{dhi} (as well as on some  
hints from V. Hinich) and uses the Bousfield localization   
techniques of \cite{hi}, but does not assume the results of \cite{bl,ja}.   
  
\begin{thm}\label{t3}  
Let $(T,\tau)$ be an $S$-site. There exists a closed model   
structure on $SPr(T)$, called the \emph{local projective model structure},   
for which the equivalences are the  
$\pi_{*}$-equivalences and the cofibrations are the cofibrations for the projective model  
structure on $SPr(T)$. Furthermore the local projective model structure is  
$\mathbb{U}$-cofibrantly generated and proper. The category $SPr(T)$ together with its local projective model structure will be denoted by $SPr_{\tau}(T)$.   
\end{thm}  
  
\textit{Proof:} We are going to apply the existence theorem for left Bousfield localizations \cite[Thm. 4.1.1]{hi}  
to the objectwise model structure $SPr(T)$ along a certain   
$\mathbb{U}$-small set $H$ of morphisms. The main point will be  
to check that equivalences in this localized model structure are exactly $\pi_{*}$-equivalences. \\  
  
\begin{center} \textit{Definition of the set $H$} \end{center}   
  
As the $S$-category $T$ is $\mathbb{U}$-small, the set   
$$E(T):=\coprod_{n \in \mathbb{N}}\;\coprod_{(x,y) \in Ob(T)^{2}}\underline{Hom}_{T}(x,y)_{n},$$  
of all simplices in all simplicial set of morphisms of $T$ is also $\mathbb{U}$-small.  
We denote by $\alpha$ a $\mathbb{U}$-small cardinal bigger than the cardinal of $E(T)$  
and than $\aleph_{0}$.  Finally, we let  
$\beta$ be a $\mathbb{U}$-small cardinal with $\beta > 2^{\alpha}$.   
  
The size of a simplicial presheaf $F \in SPr(T)$ is by definition the cardinality of the set  
$$\coprod_{n \in \mathbb{N}}\;\coprod_{x\in Ob(T)}F_{n}(x).$$  
We will denote it by $\mathrm{Card}(F)$.   
  
For an object $x \in Ob(T)$ we consider a fibrant replacement $\underline{h}_{x} \hookrightarrow R(\underline{h}_{x})$   
as well as the simplicial object it   
defines $R(\underline{h}_{x})^{\Delta^{*}} \in sSPr(T)$. Note that as $\underline{h}_{x}$ is a cofibrant object, so  
is $R(\underline{h}_{x})$.   
We define a subset $\mathcal{H}_{\beta}(x)$ of objects in $sSPr(T)/R(\underline{h}_{x})^{\Delta^{*}}$  
in the following way.   
We consider the following two conditions.  
  
\begin{enumerate}  
  
\item The morphism  $F_{*} \longrightarrow R(\underline{h}_{x})^{\Delta^{*}} \in \mathrm{Ho}(sSPr(T))$ is a hypercover.  
  
\item For all $n\geq 0$, one has $\mathrm{Card}(F_{n})<\beta$. Furthermore, for each $n\geq 0$, $F_{n}$ is isomorphic  
in $\mathrm{Ho}(SPr(T))$ to a coproduct of representable objects  
$$F_{n}\simeq \coprod_{u\in I_{n}}\underline{h}_{u}.$$  
  
\end{enumerate}  
  
We define $\mathcal{H}_{\beta}(x)$ to be a set of representatives in $sSPr(T)/R(\underline{h}_{x})^{\Delta^{*}}$,   
for the isomorphism classes of objects $F_{*} \in sSPr(T)/R(\underline{h}_{x})^{\Delta^{*}}$  
which satisfy conditions $(1)$ and $(2)$ above. Note that condition   
$(2)$ insures that $\mathcal{H}_{\beta}(x)$ is a $\mathbb{U}$-small set   
for any $x \in Ob(T)$.   
  
Now, for any $x \in Ob(T)$, any $F_{*} \in \mathcal{H}_{\beta}(x)$   
we consider its geometric realization  
$|F_{*}|$ in $SPr(T)$, together with its natural adjunction   
morphism $|F_{*}| \longrightarrow R(\underline{h}_{x})$ (see \cite[$19.5.1$]{hi}).   
Note that $|F_{*}|$ is naturally equivalent to the  
homotopy colimit of the diagram $[n] \mapsto F_{n}$. Indeed,   
for any $y\in Ob(T)$, $|F_{*}|(y)$ is naturally isomorphic  
to the diagonal of the the bi-simplicial set $F_{*}(y)$ (see \cite[$16.10.6$]{hi}). We define the set $H$ to be  
the union of all the $\mathcal{H}_{\beta}(x)$'s when $x$ varies   
in $Ob(T)$. In other words, $H$ consists of all morphisms  
$$|F_{*}| \longrightarrow R(\underline{h}_{x}),$$   
for all $x \in Ob(T)$ and all $F_{*} \in \mathcal{H}_{\beta}(x)$.   
Clearly, the set $H$ is $\mathbb{U}$-small, so one can apply  
Theorem \ref{tb1} or \ref{tb2} to the objectwise model category $SPr(T)$ and the set of morphisms $H$.   
We let $\mathrm{L}_{H}SPr(T)$ be the left Bousfield localization of   
$SPr(T)$ along the set of morphisms $H$. We are  
going to show that equivalences in $\mathrm{L}_{H}SPr(T)$ are exactly   
$\pi_{*}$-equivalences. This will clearly implies the  
existence of the local model structure of \ref{t3}. \\  
  
\begin{center}  \textit{The morphisms in $H$ are $\pi_{*}$-equivalences} \end{center}   
  
The main point in the proof is the following lemma.  
  
\begin{lem}\label{l4}  
For any object $x \in Ob(T)$ and any hypercover $F_{*} \longrightarrow   
R(\underline{h}_{x})^{\Delta^{*}}$, the natural morphism   
in $\mathrm{Ho}(SPr(T))$  
$$\mathrm{hocolim}_{[n] \in \Delta^{n}}(F_{n}) \longrightarrow   
R(\underline{h}_{x})\simeq \underline{h}_{x}$$  
is a $\pi_{*}$-equivalence.  
\end{lem}  
  
\textit{Proof:} By applying the base change functor $j_{x}^{*} :   
\mathrm{Ho}(SPr(T)) \longrightarrow \mathrm{Ho}(SPr(T/x))$   
one gets a natural morphism  
$j_{x}^{*}(\mathrm{hocolim}_{[n] \in \Delta^{n}}(F_{n})) \longrightarrow  j_{x}^{*}(\underline{h}_{x})$.  
By definition of the homotopy sheaves one sees that it is enough to show that the homotopy fiber of this  
morphism at the natural point $* \longrightarrow j_{x}^{*}\underline{h}_{x}$  
is $\pi_{*}$-contractible (see Corollary \ref{c3''}).  In other words, one can always  
assume that $x$ is a final object in $T$, or in other words that $\underline{h}_{x}\simeq *$  
(this reduction is not necessary but simplifies notations). We can also   
assume that $F_{*}$ is fibrant as an object in $sSPr(T)$, so   
$Cosk_{n}(F_{*})\simeq \mathbb{R}Cosk_{n}(F_{*})$. We will simply denote by   
$|G_{*}|$ the homotopy colimit of a simplicial diagram $G_{*}$ in $SPr(T)$. \\  
  
\textsf{Step $1$:}  
Let us first assume that $F_{*}$ is a $0$-bounded hypercover. Recall   
that this means that for any $n>0$ one has $F_{n}\simeq   
F_{*}^{\mathbb{R}\partial \Delta^{n}}$, or in other words that $F_{*}$ is the nerve of the  
covering $F_{0} \longrightarrow *$. Therefore, we can suppose that $F_{0}$ is fibrant in $SPr(T)$, and  
that $F_{n}=F_{0}^{n}$ (the face and degeneracy morphisms being induced by the various projections and  
diagonals). As $F_{0} \longrightarrow *$ is a covering, one can find a covering sieve $R$ of $*$ such that  
for any object $u\rightarrow *$ in $S$, there exists a commutative diagram  
$$\xymatrix{F_{0} \ar[r] & \textrm{\large{*}} \\  
\underline{h}_{u} \ar[u] \ar[ru] }$$  
Furthermore, as $\pi_{*}$-equivalences are local for the topology $\tau$ (see Corollary \ref{c2}), it is enough to  
prove that for any such $u$, the nerve of the morphism  
$$F_{0}\times \underline{h}_{u} \longrightarrow \underline{h}_{u}$$  
is a $\pi_{*}$-equivalence. We can therefore assume that the morphism $F_{0}\longrightarrow *$ admits a section.   
But then, for any object $x \in Ob(T)$, $|F_{*}|(x) \in \mathrm{Ho}(SSet_{\mathbb{U}})$ is the geometric  
realization of the nerve of a morphism of simplicial sets which has a section, and therefore is contractible.   
This proves Lemma \ref{l4} for $0$-bounded hypercovers. \\  
  
\textsf{Step $2$:}  
Let us now assume that $F_{*}$ is $(n+1)$-bounded for some integer   
$n>0$ (see Definition \ref{d11}), and let us consider  
the morphism   
$$F_{*} \longrightarrow Cosk_{n}F_{*}.$$  
For any integer $p$, and any simplicial set $K \in SSet_{\mathbb{U}}$,   
there is a co-cartesian square of simplicial sets  
$$\xymatrix{  
Sk_{p}K \ar[r] & Sk_{p+1}K \\  
\coprod_{K^{\partial \Delta^{p+1}}}\partial \Delta^{p+1} \ar[r] \ar[u] & \coprod_{K_{p+1}}\Delta^{p+1} \ar[u].}$$  
This induces a cartesian square in $SPr(T)$  
$$\xymatrix{  
F_{*}^{Sk_{p+1}K} \ar[r] \ar[d] & F_{*}^{Sk_{p}K} \ar[d] \\  
\ar[r] \prod_{K_{p+1}}F_{p+1}& \prod_{K^{\partial \Delta^{p+1}}}F_{*}^{\partial \Delta^{p+1}}.}$$  
As $F_{*}$ is fibrant for the Reedy structure and a hypercover, each bottom horizontal morphism   
is a fibration which is again a covering. This shows by induction and by Proposition \ref{p1} $(1)$, that   
$F_{*}^{Sk_{p+i}K} \longrightarrow F_{*}^{Sk_{p}K}$  
is a covering and a fibration for any $i>0$. But, since we have  
$$(Cosk_{n}F_{*})^{K}\simeq F_{*}^{Sk_{n}K},$$  
we easily conclude that for any $K \in SSet_{\mathbb{U}}$ such that $K=Sk_{p}K$ for some $p$, the natural morphism  
$$F_{*}^{K} \longrightarrow (Cosk_{n}F_{*})^{K}$$  
is again a fibration and a covering. In particular, taking $K=\Delta^{p}$, one finds that the natural morphism   
$$F_{p} \longrightarrow (Cosk_{n}F_{*})_{p}.$$  
is a fibration and a covering.  
  
Let $U_{*,*}$ be the bi-simplicial object such that $U_{p,*}$ is the nerve of the morphism   
$F_{p} \longrightarrow (Cosk_{n}F_{*})_{p}$. It fits into a commutative diagram of  
bi-simplicial objects  
$$\xymatrix{  
F_{*} \ar[r] \ar[d] & Cosk_{n}F_{*} \\  
U_{*,*}, \ar[ru]}$$  
where $F_{*}$ and $Cosk_{n}F_{*}$ are considered as constant in the second spot. Furthermore, for any $p$,   
$U_{p,*} \longrightarrow (Cosk_{n}F_{*})_{p}$ is a $0$-truncated hypercover. Therefore, by   
\textsf{Step $1$}, we deduce that   
$$|diag(U_{*,*})|\simeq \mathrm{hocolim}_{p}\mathrm{hocolim}_{q}(U_{p,q}) \longrightarrow  
|Cosk_{n}F_{*}|$$  
is a $\pi_{*}$-equivalence.   
  
Now, let $U_{*}:=diag(U_{*,*})$ be the diagonal of $U_{*,*}$. It fits into a commutative diagram  
$$\xymatrix{  
F_{*} \ar[r]^-{\pi} \ar[d]_{f} & Cosk_{n}F_{*} \\  
U_{*}. \ar[ru]_-{\phi}}$$  
We are going to construct a morphism $U_{*} \longrightarrow F_{*}$ that will be  
a retract of $f$ compatible with the two projections $\pi$ and $\phi$ (i.e. construct a retraction  
of $\phi$ on $\pi$).  
  
The above diagram consists clearly of isomorphisms in degrees $p\leq n$, showing that $\pi$  
is a retract of $\phi$ is degrees $p\leq n$.  
As $F_{*}$ is $(n+1)$-bounded, to extend this retraction to the whole $\phi$, it is  
enough to define a morphism $U_{n+1} \longrightarrow F_{n+1}$ which is equalized by all the  
face morphisms $F_{n+1} \longrightarrow F_{n}$. But, by definition   
$$U_{n+1}=\underbrace{F_{*}^{\Delta^{n+1}}\times_{F_{*}^{\partial \Delta^{n+1}}}F_{*}^{\Delta^{n+1}}\times \dots  
\times_{F_{*}^{\partial \Delta^{n+1}}}F_{*}^{\Delta^{n+1}}}_{\textrm{(n+1)  times}},$$  
and so any of the natural projections $U_{n+1} \longrightarrow F_{n+1}$   
to one of these factors will produce the required extension.   
  
In conclusion, the morphism $F_{*} \longrightarrow Cosk_{n}F_{*}$ is a retract  
of $U_{*} \longrightarrow Cosk_{n}F_{*}$, which itself induces a $\pi_{*}$-equivalence on the  
homotopy colimits. As $\pi_{*}$-equivalences are stable by retracts, this shows that the  
induced morphism   
$|F_{*}| \longrightarrow |Cosk_{n}F_{*}|$  
is also a $\pi_{*}$-equivalence. Therefore, by induction on $n$ and   
\textsf{Step $1$}, this implies that   
$|F_{*}| \longrightarrow *$  
is a $\pi_{*}$-equivalence. \\  
  
\textsf{Step $3$:}  
Finally, for a general hypercover $F_{*}$, the $i$-th homotopy presheaf  
of $|F_{*}|$ only depends on the $n$-th coskeleton of $F_{*}$ for $i<n$  
(as the $(n-1)$-skeleton of $|F_{*}|$ and $|Cosk_{n}F_{*}|$ coincide). In particular, the $i$-th homotopy sheaf of  
$|F_{*}|$ only depends on $\mathbb{R}Cosk_{n}(F_{*})$ for $i<n$. Therefore  
one can always suppose that $F_{*}=Cosk_{n}F_{*}$ for some integer $n$ and apply \textsf{Step $2$}.  
  
Lemma \ref{l4} is proved.   
\hfill \textbf{$\Box$} \\  
  
Now, let $f : F \longrightarrow G$ be any $H$-local equivalence (i.e. an   
equivalence in $\mathrm{L}_{H}SPr(T)$), and let us   
prove that it is a $\pi_{*}$-equivalence. By definition of $H$-local equivalences,    
the induced morphism on the $H$-local models  
$$\mathrm{L}_{H}f : \mathrm{L}_{H}F \longrightarrow \mathrm{L}_{H}G$$  
is an objectwise equivalence, and in particular a $\pi_{*}$-equivalence. By considering the commutative diagram  
$$\xymatrix{F \ar[r]^-{f} \ar[d] & G \ar[d] \\  
\mathrm{L}_{H}F \ar[r]_-{\mathrm{L}_{H}f} & \mathrm{L}_{H}G,}$$  
one sees that it is enough to show that the localization morphisms $F \longrightarrow \mathrm{L}_{H}F$ and   
$G\longrightarrow \mathrm{L}_{H}G$  
are $\pi_{*}$-equivalences. But the functor $\mathrm{L}_{H}$ can be defined via the small object argument applied to  
the set of augmented horns on $H$, $\overline{\Lambda(H)}$   
(see \cite[\S $4.3$]{hi}). In the present situation, the morphisms  
in $\overline{\Lambda(H)}$ are either trivial cofibrations in $SPr(T)$   
or projective cofibrations which are isomorphic in   
$\mathrm{Ho}(SPr(T))$ to  
$$\Delta^{n}\otimes |F_{*}|\coprod^{h}_{\partial \Delta^{n}\otimes |F_{*}|}  
\partial \Delta^{n}\otimes R(\underline{h}_{x}) \longrightarrow \Delta^{n}\otimes R(\underline{h}_{x}).$$  
By Proposition \ref{c3'} and Lemma \ref{l4}, these morphisms are $\pi_{*}$-equivalences, and  
therefore all morphisms in $\overline{\Lambda(H)}$ are projective cofibrations and $\pi_{*}$-equivalences.  
As $\pi_{*}$-equivalences are also stable  
by filtered colimits, another application of Proposition \ref{c3'} shows that   
relative cell complexes on $\overline{\Lambda(H)}$ are again $\pi_{*}$-equivalences. This shows that   
the localization morphisms $F \longrightarrow \mathrm{L}_{H}F$ are always   
$\pi_{*}$-equivalences, and finish the proof that   
$H$-local equivalences are $\pi_{*}$-equivalences. \\  
  
\begin{center} \textit{$\pi_{*}$-Equivalences are $H$-local equivalences} \end{center}   
  
To conclude the proof of Theorem \ref{t3}, we are left to show that $\pi_{*}$-equivalences are $H$-local equivalences.   
  
Recall that we denoted by $\alpha$ a $\mathbb{U}$-small cardinal bigger than $\aleph_{0}$ and than   
the cardinality of the set $E(T)$ of  
all simplicies in all simplicial set of morphisms  
in $T$. Recall also that  
$\beta$ is a $\mathbb{U}$-small cardinal with $\beta > 2^{\alpha}$.   
  
\begin{lem}\label{l5}  
Let $f : F \longrightarrow G$ be a morphism in $SPr(T)$ which is a $\pi_{*}$-equivalence and   
an objectwise fibration between fibrant objects. Then, for any object   
$x \in Ob(T)$ and any morphism $R(\underline{h}_{x}) \longrightarrow G$,   
there exists an $F_{*} \in \mathcal{H}_{\beta}(x)$ and a commutative diagram in $SPr(T)$  
$$\xymatrix{  
F \ar[r] & G \\  
|F_{*}| \ar[r] \ar[u] & R(\underline{h}_{x}). \ar[u]}$$  
\end{lem}  
  
\textit{Proof:} By adjunction, it is equivalent to find a commutative diagram in $sSPr(T)$  
$$\xymatrix{  
F^{\Delta^{*}} \ar[r] & G^{\Delta^{*}} \\  
F_{*} \ar[r] \ar[u] & R(\underline{h}_{x})^{\Delta^{*}}, \ar[u]}$$  
with $F_{*} \in \mathcal{H}_{\beta}(x)$. We will define $F_{*}$ inductively. Let us suppose   
we have constructed $F(n)_{*} \in sSPr(T)/R(\underline{h}_{x})^{\Delta^{*}}$, with  
a commutative diagram  
$$\xymatrix{  
F^{\Delta^{*}} \ar[r] & G^{\Delta^{*}} \\  
F(n)_{*} \ar[r]_{p_{n}} \ar[u] & R(\underline{h}_{x})^{\Delta^{*}}, \ar[u]}$$  
such that $Sk_{n}F(n)_{*}=F(n)_{*}$, and $p_{n}$ is a Reedy fibration and   
a hypercover in degrees $i\leq n$. By the latter condition we mean that  
$$F(n)_{i} \longrightarrow F(n)^{\partial \Delta^{i}}\times_{R(\underline{h}_{x})^{\partial \Delta^{i}}}  
R(\underline{h}_{x})^{\Delta^{i}}$$  
is an objectwise fibration and a covering for any $i\leq n$ (we do not require $p_{n}$ to be a Reedy fibration).  
We also assume that $\mathrm{Card}(F(n)_{m})<\beta$ for any $m$.  
We need the following (technical) factorization result with control on the cardinality.  
  
\begin{lem}\label{tech}  
Let $f : F \longrightarrow G$ be a morphism in $SPr(T)$ such that   
$\mathrm{Card}(F)$ and $\mathrm{Card}(G)$ are both stricly smaller than $\beta$. Then, there   
exists a factorization in $SPr(T)$  
$$\xymatrix{F \ar[r]^-{i}  & RF \ar[r]^-{p} & G},$$  
with $i$ a trivial cofibration, $p$ a fibration, and $\mathrm{Card}(RF)<\beta$.   
\end{lem}  
  
\textit{Proof:} We use the standard small object argument in order to   
produce such a factorization (see \cite[\S $2.1.2$]{ho}).  
The trivial cofibrations in $SPr(T)$ are generated by the set of morphisms  
$$\Lambda^{n,k}\otimes \underline{h}_{x} \longrightarrow \Delta^{n}\otimes \underline{h}_{x},$$  
for all $x \in Ob(T)$ and all $n \in \mathbb{N}$, $0\leq k\leq n$. This set is clearly of cardinality smaller than  
$\aleph_{0}.\alpha$, and therefore is stricly smaller than $\beta$. Furthermore, for any of these generating   
trivial cofibrations, the set of all commutative diagrams  
$$\xymatrix{  
F \ar[r] & G \\  
\Lambda^{n,k}\otimes \underline{h}_{x} \ar[r] \ar[u] & \Delta^{n}\otimes \underline{h}_{x}, \ar[u]}$$  
is in bijective correspondence with the set of all commutative diagrams  
$$\xymatrix{  
F(x) \ar[r] & G(x) \\  
\Lambda^{n,k} \ar[r] \ar[u] & \Delta^{n}. \ar[u]}$$  
By the assumptions made on $F$ and $G$, this set is therefore of   
cardinality stricly smaller than $\beta$. Furthermore,   
by the choice of $\beta$, it is clear that $\mathrm{Card}(A\otimes \underline{h}_{x})\leq \alpha<\beta$ for  
any finite simplicial set $A$. Therefore, the push-out  
$$\xymatrix{  
F \ar[r] & F_{1} \\  
\coprod_{I}\Lambda^{n,k}\otimes \underline{h}_{x} \ar[r] \ar[u] &   
\coprod_{I}\Delta^{n}\otimes \underline{h}_{x} \ar[u],}$$  
where $I$ consists of all objects $x \in Ob(T)$ and   
commutative diagrams  
$$\xymatrix{  
F \ar[r] & G \\  
\Lambda^{n,k}\otimes \underline{h}_{x} \ar[r] \ar[u] & \Delta^{n}\otimes \underline{h}_{x}, \ar[u]}$$  
is such that   
$$\mathrm{Card}(F_{1})\leq \mathrm{Card}(F)+\mathrm{Card}(\coprod_{I}\Delta^{n}\otimes \underline{h}_{x})<  
\beta+\mathrm{Card}(I).\alpha.$$   
But $\mathrm{Card}(I)< \alpha.\beta$, and therefore one has $\mathrm{Card}(F_{1})<\beta$.  
As the factorization $\xymatrix{F \ar[r] & RF \ar[r] & G}$  
is obtained after a numerable number of such push-outs constructions (see \cite[Thm. $2.1.14$]{ho})  
$$\xymatrix{F \ar[r] & F_{1} \ar[r] & \dots \ar[r] & F_{n} \ar[r] & \dots \ar[r] & RF=colim_{i}F_{i},}$$  
we conclude that   
$\mathrm{Card}(RF)<\beta$. The proof of Lemma \ref{tech} is achieved. \hfill \textbf{$\Box$} \\  
  
Let us come back to the proof of Lemma \ref{l5}.  
We consider the following diagram  
$$\xymatrix{  
F^{\Delta^{n+1}} \ar[r] & F^{\partial \Delta^{n+1}}\times_{G^{\partial \Delta^{n+1}}}G^{\Delta^{n+1}} \\  
 &  F(n)_{*}^{\partial \Delta^{n+1}}\times_{R(\underline{h}_{x})^{\partial \Delta^{n+1}}}  
R(\underline{h}_{x})^{\Delta^{n+1}}. \ar[u] }$$  
  
By Lemma \ref{tech}, we can suppose that $\mathrm{Card}(R(\underline{h}_{x}))<\beta$. Therefore, by induction on $n$  
$$\mathrm{Card}(F(n)_{*}^{\partial \Delta^{n+1}}\times_{R(\underline{h}_{x})^{\partial \Delta^{n+1}}}  
R(\underline{h}_{x})^{\Delta^{n+1}})<\beta.$$  
This implies that there exists a $\mathbb{U}$-small set $J$ of objects in $T$, with $\mathrm{Card}(J)<\beta$, and a covering  
$$\coprod_{z \in J}\underline{h}_{z} \longrightarrow   
F(n)_{*}^{\partial \Delta^{n+1}}\times_{R(\underline{h}_{x})^{\partial \Delta^{n+1}}}  
R(\underline{h}_{x})^{\Delta^{n+1}}.$$  
Now, by considering the induced diagram  
$$\xymatrix{  
F^{\Delta^{n+1}} \ar[r] & F^{\partial \Delta^{n+1}}\times_{G^{\partial \Delta^{n+1}}}G^{\Delta^{n+1}} \\  
 &  \coprod_{z \in J}\underline{h}_{z}, \ar[u] }$$  
and using the fact that the top horizontal morphism is a covering, one sees that there  
exists, for all $z\in J$, a covering sieve $S_{z}$ of $z \in \mathrm{Ho}(T)$, and a commutative diagram  
$$\xymatrix{  
F^{\Delta^{n+1}} \ar[r] & F^{\partial \Delta^{n+1}}\times_{G^{\partial \Delta^{n+1}}}G^{\Delta^{n+1}} \\  
\coprod_{z \in J, (u\rightarrow z) \in S_{z}}\underline{h}_{u} \ar[u] \ar[r]  
 &  \coprod_{z \in J}\underline{h}_{z}. \ar[u] }$$  
Clearly, one has  
$$\mathrm{Card}(\coprod_{z \in J, (u\rightarrow z) \in S_{z}}\underline{h}_{u})  
\leq \mathrm{Card}(J).2^{\alpha}.\alpha < \beta.$$  
We now consider the commutative diagram  
$$\xymatrix{  
F^{\Delta^{n+1}} \ar[r] & F^{\partial \Delta^{n+1}}\times_{G^{\partial \Delta^{n+1}}}G^{\Delta^{n+1}} \\  
\coprod_{z \in J, (u\rightarrow z) \in S_{z}}\underline{h}_{u} \ar[r]  \ar[u] &    
F(n)_{*}^{\partial \Delta^{n+1}}\times_{R(\underline{h}_{x})^{\partial \Delta^{n+1}}}  
R(\underline{h}_{x})^{\Delta^{n+1}}. \ar[u] }$$  
Lemma \ref{tech} implies the existence of an object $H(n+1) \in SPr(T)$, with $\mathrm{Card}(H(n+1))<\beta$, and a factorization  
$$\xymatrix{  
\coprod_{z \in J, (u\rightarrow z) \in S_{z}}\underline{h}_{u} \ar[r] &   
H(n+1) \ar[r] & F(n)_{*}^{\partial \Delta^{n+1}}\times_{R(\underline{h}_{x})^{\partial \Delta^{n+1}}}  
R(\underline{h}_{x})^{\Delta^{n+1}},}$$  
into an objectwise trivial cofibration followed by a fibration in $SPr(T)$. Since the morphism  
$$F^{\Delta^{n+1}} \longrightarrow F^{\partial \Delta^{n+1}}\times_{G^{\partial \Delta^{n+1}}}G^{\Delta^{n+1}}$$  
is an objectwise fibration, there exists a commutative diagram in $SPr(T)$  
$$\xymatrix{  
F^{\Delta^{n+1}} \ar[rr] & & F^{\partial \Delta^{n+1}}\times_{G^{\partial \Delta^{n+1}}}G^{\Delta^{n+1}} \\  
\coprod_{z \in J, (u\rightarrow z) \in S_{z}}\underline{h}_{u} \ar[r] \ar[u] &  
H(n+1) \ar[ul] \ar[r] & F(n)_{*}^{\partial \Delta^{n+1}}\times_{R(\underline{h}_{x})^{\partial \Delta^{n+1}}}  
R(\underline{h}_{x})^{\Delta^{n+1}}. \ar[u]}$$  
  
We define $F(n+1)_{p}:=F(n)_{p}$ for any $p<n+1$, and  
$F(n+1)_{n+1}$ to be the coproduct of $H(n+1)$ together with $L_{n+1}F$,
the $(n+1)$-th latching space of $F(n)$. The face morphisms   
$F(n+1)_{n+1} \longrightarrow F(n)_{n}$ are defined as the identity on 
$L_{n+1}F(n)$ and  
via the $(n+1)$ natural projections (corresponding to the face inclusions $\Delta^{n} \subset \partial \Delta^{n+1}$)  
$$F(n)^{\partial \Delta^{n+1}} \longrightarrow F(n)^{\Delta^{n}}=F(n)_{n}$$  
on the factor $H(n+1)$. Then, by adjunction, one has a natural   
commutative diagram in $s_{n+1}SPr(T)$  
$$\xymatrix{  
F^{\Delta^{*}} \ar[r] & G^{\Delta^{*}} \\  
F(n+1)_{*} \ar[r]_-{p_{n+1}} \ar[u] & R(\underline{h}_{x})^{\Delta^{*}}, \ar[u]}$$  
which extends via the functor $(i_{n+1})_{!}$ to the required diagram in $sSPr(T)$.   
It is clear by construction, that $p_{n+1}$ is a Reedy fibration and a   
hypercover in degrees $i\leq n+1$ and that its $n$-th skeleton is $p_{n}$.   
Therefore, by defining $F_{*}$ to be the limit of the $F(n)$'s,   
the natural morphism $F_{*} \longrightarrow R(\underline{h}_{x})^{\Delta^{*}}$  
is a hypercover. It is also clear by construction that $F_{*}$ satisfies   
condition $(2)$ defining the set $\mathcal{H}_{\beta}(x)$.   
\hfill \textbf{$\Box$} \\  
  
We are now ready to finish the proof that $\pi_{*}$-equivalences are $H$-local equivalences.   
Let $f : F \longrightarrow G$ be a $\pi_{*}$-equivalence; we can clearly assume $f$ to be an objectwise fibration  
between fibrant objects. Furthermore, as $H$-local equivalences are already known to be $\pi_{*}$-equivalences,   
we can also suppose that $f$ is a $H$-local fibration between $H$-local objects.   
We are going to prove that $f$ is in fact   
an objectwise equivalence.  
  
Let   
$$\xymatrix{  
F \ar[r]^{f} & G \\  
\partial \Delta^{n}\otimes \underline{h}_{x} \ar[r] \ar[u] & \Delta^{n}\otimes \underline{h}_{x} \ar[u]}$$  
be a commutative diagram in $SPr(T)$. We need to show that there exist a   
lifting $\Delta^{n}\otimes \underline{h}_{x} \longrightarrow F$.   
By adjunction, this is equivalent to showing that the natural morphism   
$$\underline{h}_{x} \longrightarrow F^{\partial \Delta^{n}}\times_{G^{\partial \Delta^{n}}}G^{\Delta^{n}}$$  
lifts to a morphism $\underline{h}_{x} \longrightarrow F^{\Delta^{n}}$.   
  
As $F$ and $G$ are objectwise fibrant, the previous morphism factors through  
$$\underline{h}_{x} \longrightarrow R(\underline{h}_{x})\longrightarrow   
F^{\partial \Delta^{n}}\times_{G^{\partial \Delta^{n}}}G^{\Delta^{n}}.$$  
An application of Lemma \ref{l5} to the morphism   
$$F^{\Delta^{n}} \longrightarrow F^{\partial \Delta^{n}}\times_{G^{\partial \Delta^{n}}}G^{\Delta^{n}},$$  
which satisfies the required hypothesis,  
shows that there exists an $F_{*} \in \mathcal{H}_{\beta}(x)$ and a commutative diagram  
$$\xymatrix{  
F^{\Delta^{n}} \ar[r] & F^{\partial \Delta^{n}}\times_{G^{\partial \Delta^{n}}}G^{\Delta^{n}} \\  
|F_{*}| \ar[r] \ar[u] & R(\underline{h}_{x}). \ar[u]}$$  
By adjunction, this commutative diagram yields a commutative diagram  
$$\xymatrix{  
F \ar[r]^{f} & G \\  
\Delta^{n}\otimes |F_{*}|\coprod_{\partial \Delta^{n}\otimes |F_{*}|}\partial \Delta^{n}\otimes R(\underline{h}_{x})   
\ar[r] \ar[u] & \Delta^{n}\otimes R(\underline{h}_{x}).\ar[u]}$$  
The horizontal bottom morphism is an $H$-local equivalence by definition, and therefore a lifting   
$\Delta^{n}\otimes R(\underline{h}_{x}) \longrightarrow F$ exists in the   
homotopy category $\mathrm{Ho}(\mathrm{L}_{H}SPr(T))$.   
But, as $f$ is a $H$-local fibration, $F$ and $G$ are $H$-local objects and $R(\underline{h}_{x})$   
is cofibrant, this lifting can be represented in $SPr(T)$ by a commutative diagram  
$$\xymatrix{  
F \ar[r]^{f} & G \\  
& \ar[ul]\Delta^{n}\otimes R(\underline{h}_{x}).\ar[u]}$$  
Composing with $\underline{h}_{x} \longrightarrow R(\underline{h}_{x})$, we obtain the   
required lifting. This implies that $\pi_{*}$-equivalences are $H$-local equivalences, and completes the proof of the  
existence of the local model structure.  
  
By construction, $SPr_{\tau}(T)$ is the left Bousfield localization of   
$SPr(T)$ along the set of morphisms $H$: this implies   
that it is a $\mathbb{U}$-cellular and $\mathbb{U}$-combinatorial model category. In particular, it is  
$\mathbb{U}$-cofibrantly generated.  
Finally, properness of $SPr_{\tau}(T)$ follows from Corollary \ref{c3} and Proposition \ref{c3'}.   
  
This concludes the proof of Theorem \ref{t3}. \hfill \textbf{$\Box$} \\  
  
Let us keep the notations introduced in the proof of Theorem \ref{t3}.   
We choose a $\mathbb{U}$-small cardinal $\beta$ as in the proof and consider,   
for any object $x \in Ob(T)$, the subset of hypercovers $\mathcal{H}_{\beta}(x)$.   
  
\begin{cor}\label{c4}  
The model category $SPr_{\tau}(T)$ is the left Bousfield localization of  
$SPr(T)$ with respect to the set of morphisms $$\left\{|F_{*}| \longrightarrow \underline{h}_{x} \;|\;  
x \in Ob(T), \;  F_{*}\in \mathcal{H}_{\beta}(x) \right\}.$$  
\end{cor}  
  
\textit{Proof:} This is exactly the way we proved Theorem \ref{t3}. \hfill \textbf{$\Box$} 
  
\begin{rmk} \label{twoinone} \emph{It is worthwile emphasizing that the proof of Theorem   
\ref{t3} shows actually a bit more than what's in its statement. In fact,   
the argument proves both Theorem \ref{t3} and Corollary \ref{c4}, in that   
it gives \textit{two descriptions} of the same model category $SPr_{\tau}(T)$:   
one as the left Bousfield localization of $SPr(T)$ with respect to \textit{local equivalences}   
and the other as the left Bousfield localization of the same $SPr(T)$ but this time   
with respect to \textit{hypercovers} (more precisely, with respect to the set of   
morphisms defined in the statement of Corollary \ref{c4}).}  
\end{rmk}  
  
In the special case where $(T,\tau)$ is a usual Grothendieck site (i.e. when $T$ is a category),   
the following corollary was announced in \cite{du} and proved in \cite{dhi}.   
  
\begin{cor}\label{c6}  
An object $F \in SPr_{\tau}(T)$ is fibrant if and only if it is objectwise fibrant and  
for any object $x\in Ob(T)$ and any $H_{*}\in \mathcal{H}_{\beta}(x)$,   
the natural morphism   
$$F(x)\simeq \mathbb{R}\underline{Hom}(\underline{h}_{x},F) \longrightarrow \mathbb{R}\underline{Hom}(|H_{*}|,F)$$  
is an isomorphism in $\mathrm{Ho}(SSet)$.  
\end{cor}  
  
\textit{Proof:} This follows from Thm. \ref{t3} and from the explicit description of fibrant objects in a   
left Bousfield localization  
(see \cite[Thm. 4.1.1]{hi}). \hfill \textbf{$\Box$} \\  
  
The previous corollary is more often described in the following way. For any  
$H_{*}\in \mathcal{H}_{\beta}(x)$ and any $n \geq 0$, $H_{n}$ is equivalent to a coproduct of representables  
$$H_{n}\simeq \coprod_{i \in I_{n}}\underline{h}_{u_{i}}.$$  
Therefore, for any $H_{*} \in \mathcal{H}_{\beta}(x)$ and any fibrant object  $F$ in $SPr(T)$,   
the simplicial set $\mathbb{R}\underline{Hom}(|H_{*}|,F)$ is naturally equivalent to the  
homotopy limit of the cosimplicial diagram in $SSet$  
$$[n] \mapsto \prod_{i\in I_{n}}F(u_{i}).$$  
Then, Corollary \ref{c6} states that an object $F \in SPr(T)$ is fibrant if and only if,   
for any $x\in Ob(T)$, $F(x)$ is fibrant, and the natural morphism  
$$F(x) \longrightarrow \mathrm{holim}_{[n] \in \Delta}\left(\prod_{i\in I_{n}}F(u_{i})\right)$$  
is an equivalence of simplicial sets, for any $H_{*} \in \mathcal{H}_{\beta}(x)$.  
  
\medskip  
  
\begin{df}\label{d15}  
\begin{enumerate}  
\item A hypercover $H_{*} \longrightarrow \underline{h}_{x}$ is said to be \emph{semi-representable} if  
for any $n \geq 0$, $H_{n}$ is isomorphic in $\mathrm{Ho}(SPr(T))$ to   
a coproduct of representable objects  
$$H_{n} \simeq \coprod_{u\in I_{n}}\underline{h}_{u}.$$  
  
\item An object $F \in SPr(T)$ is said to \emph{have hyperdescent} if,   
for any object $x \in Ob(T)$ and any semi-representable hypercover $H_{*} \longrightarrow \underline{h}_{x}$,   
the induced  
morphism  
$$F(x)\simeq \mathbb{R}\underline{Hom}(\underline{h}_{x},F)\longrightarrow    
\mathbb{R}\underline{Hom}(|H_{*}|,F)$$  
is an isomorphism in $\mathrm{Ho}(SSet_{\mathbb{U}})$.   
\end{enumerate}  
\end{df}  
  
An immediate consequence of the proof of Theorem \ref{t3} is that an object  
$F \in SPr(T)$ has hyperdescent with respect to all hypercover $H_{*} \in \mathcal{H}_{\beta}(x)$   
\textit{if and only} if it has hyperdescent with respect to all semi-representable hypercovers.   
  
From now on we will adopt the following terminology and notations.  
  
\begin{df}\label{d16}  
  
Let $(T,\tau)$ be an $S$-site in $\mathbb{U}$.  
\begin{enumerate}  
  
\item A \emph{stack} on the site $(T,\tau)$ is a pre-stack $F \in SPr(T)$ which satisfies the  
hyperdescent condition of Definition \ref{d15}.   
  
\item The model category $SPr_{\tau}(T)$ is also called the \emph{model category of stacks} on the $S$-site $(T,\tau)$. The category $\mathrm{Ho}(SPr(T))$ (resp. $\mathrm{Ho}(SPr_{\tau}(T))$) is called the \emph{homotopy category of pre-stacks}, and  
(resp. the \emph{homotopy category of stacks}). Objects of $\mathrm{Ho}(SPr(T))$ (resp.   
$\mathrm{Ho}(SPr_{\tau}(T))$) will simply be called \emph{pre-stacks} on $T$ (resp., \emph{stacks} on $(T,\tau)$). The functor $a : \mathrm{Ho}(SPr(T)) \longrightarrow   
\mathrm{Ho}(SPr_{\tau}(T))$ will be called the \emph{associated stack functor}.

  
\item The topology $\tau$ is said to be \emph{sub-canonical} if for any $x \in Ob(T)$, the 
pre-stack $\underline{h}_{x} \in \mathrm{Ho}(SPr(T))$ is a stack (in other words, if the Yoneda   
embedding $L\underline{h} : \mathrm{Ho}(T) \longrightarrow \mathrm{Ho}(SPr(T))$ factors through the subcategory  
of stacks).  
    
\item For pre-stacks $F$ and $G$ on $T$, we will denote by  
$\mathbb{R}\underline{Hom}(F,G) \in \mathrm{Ho}(SSet_{\mathbb{U}})$ (resp. by $\mathbb{R}_{\tau}\underline{Hom}(F,G) \in \mathrm{Ho}(SSet_{\mathbb{U}})$)   
the derived $Hom$-simplicial set computed in the simplicial model category $SPr(T)$ (resp. $SPr_{\tau}(T)$).  
  
\end{enumerate}  
\end{df}  
  
\medskip  
  
Let's explain why, given Definition \ref{d16} (1), we also call the objects in $\mathrm{Ho}(SPr_{\tau}(T))$ stacks ( Definition \ref{d16} (2)). As $SPr_{\tau}(T)$ is a left Bousfield localization of $SPr(T)$, the identity functor  
$SPr(T) \longrightarrow SPr_{\tau}(T)$ is left Quillen, and its right adjoint (which is still the identity functor) induces a fully faithful functor  
$$j : \mathrm{Ho}(SPr_{\tau}(T)) \longrightarrow \mathrm{Ho}(SPr(T)).$$  
Furthermore, the essential image of this inclusion functor is exactly the full subcategory consisting  
of objects having the hyperdescent property; in other words, the essential image of $j$ is the full subcategory of $\mathrm{Ho}(SPr(T))$ consisting of stacks. We will often identify $\mathrm{Ho}(SPr_{\tau}(T))$ with its essential image via $j$ (which is equivalent to $\mathrm{Ho}(SPr_{\tau}(T))$). The left adjoint    
$$a : \mathrm{Ho}(SPr(T)) \longrightarrow \mathrm{Ho}(SPr_{\tau}(T))$$  
to the inclusion $j$, is a left inverse to $j$. Note that $F\in \mathrm{Ho}(SPr(T))$ is a stack iff the canonical adjunction map $F\rightarrow ja(F)$ (which we will write as $F\rightarrow a(F)$ taking into account our identification) is an isomorphism in $\mathrm{Ho}(SPr(T))$.  
  
As explained in the Introduction, this situation is the analog for stacks   
over $S$-sites of the usual picture for sheaves over Grothendieck sites.   
In particular, this gives a \textit{sheaf-like} description of objects of   
$\mathrm{Ho}(SPr_{\tau}(T))$, via the hyperdescent property.   
However, this description is not as useful as one might at first think, though it allows to prove easily that some adjunctions are Quillen adjunctions (see for example,   
\cite[$7.1$]{dhi}, \cite[Prop. $2.2.2$]{to2} and \cite[Prop. $2.9$]{to3})  
or to check that an $S$-topology is sub-canonical.   
  
We will finish this paragraph with the following proposition.  
  
\begin{prop}\label{p3}  
\begin{enumerate}  
\item Let $F$ and $G$ be two pre-stacks on $T$. If $G$ is a stack, then the natural morphism  
$$\mathbb{R}\underline{Hom}(F,G) \longrightarrow \mathbb{R}_{\tau}\underline{Hom}(F,G)$$  
is an isomorphism in $\mathrm{Ho}(SSet)$.  
  
\item The functor $\mathrm{Id} : SPr(T) \longrightarrow SPr_{\tau}(T)$ preserves   
homotopy fibered products.  
  
\end{enumerate}  
\end{prop}  
  
\textit{Proof:} $(1)$ follows formally from Corollary \ref{c4}. To prove $(2)$ it is enough to show    
that $\pi_{*}$-equivalences are stable under pull-backs along    
objectwise fibrations, and this follows from Corollary \ref{c3}.   
\hfill \textbf{$\Box$} 
  
\begin{rmk}  
\emph{If $M$ is \textit{any} left proper $\mathbb{U}$-combinatorial or   
$\mathbb{U}$-cellular (see Appendix A) simplicial model category, one  
can also define the local projective model structure on $Pr(T,M):=M^{T^{op}}$ as the left Bousfield localization of  
the objectwise model structure, obtained by \textit{inverting hypercovers}.   
This allows one to consider the model category   
of stacks on the $S$-site $(T,\tau)$ \textit{with values in} $M$. Moreover, in many cases   
(e.g., symmetric spectra \cite{hss}, simplicial abelian groups, simplicial   
groups, etc.) the local equivalences also have a description in terms of some   
appropriately defined $\pi_{*}$-equivalences.   
We will not pursue this here as it is a purely formal exercise to adapt the   
proof of Theorem \ref{t3} to these situations. }  
  
\emph{In many cases these model categories of stacks with values in $M$ may   
also be described by performing the constructions defining $M$ \textit{directly} in the  
model category $SPr_{\tau}(T)$. More precisely, one can consider e.g. the   
categories of symmetric spectra, abelian group objects, group objects etc.,   
in $SPr(T)$, and use some general results to provide these categories with model   
structures. For reasonable model categories $M$ both approaches give Quillen equivalent  
model categories (e.g. for group objects in $SPr_{\tau}(T)$, and \textit{stacks of  
simplicial groups} on $(T,\tau)$). The reader might wish to consult   
\cite{bk} in which a very general approach to these considerations  
is proposed.}  
\end{rmk}  
  
\end{subsection}  
  
\begin{subsection}{Functoriality}\label{Sfunctoriality}  
  
Let $(T,\tau)$ and $(T',\tau')$ be two $\mathbb{U}$-small $S$-sites and $f : T \longrightarrow T'$ a morphism  
of $S$-categories. As we saw in Subsection \ref{diag} before Thm. \ref{t1}, the morphism $f$ induces a   
Quillen adjunction on the model categories of pre-stacks  
$$f_{!} : SPr(T) \longrightarrow SPr(T') \qquad SPr(T) \longleftarrow SPr(T') : f^{*}.$$  
  
\begin{df}\label{d17}  
We say that the morphism $f$ is \emph{continuous} (with respect to the topologies $\tau$ and $\tau'$) if  
the functor $f^{*} :  SPr(T') \longrightarrow SPr(T)$ preserves the subcategories  
of stacks.   
\end{df}  
  
As the model categories of stacks $SPr_{\tau}(T)$ and $SPr_{\tau'}(T)$ are left Bousfield localizations of  
$SPr(T)$ and $SPr(T')$, respectively, the general machinery of \cite{hi} implies that $f$ is continuous if and   
only if the adjunction $(f_{!},f^{*})$ induces a Quillen adjunction  
$$f_{!} : SPr_{\tau}(T) \longrightarrow SPr_{\tau}(T') \qquad SPr_{\tau}(T) \longleftarrow SPr_{\tau'}(T') : f^{*}$$  
between the model category of stacks.  
  
Recall from the proof of Theorem \ref{t3} that we have defined the sets of distinguished hypercovers  
$\mathcal{H}_{\beta}(x)$, for any object $x \in T$. These distinguished hypercovers   
detect continuous functors, as shown in the following proposition.  
  
\begin{prop}\label{p4}  
The morphism $f$ is continuous if and only if, for any $x\in Ob(T)$ and any $H_{*}\in \mathcal{H}_{\beta}(x)$,   
the induced morphism  
$$\mathbb{L}f_{!}(|H_{*}|) \longrightarrow \mathbb{L}f_{!}(\underline{h}_{x})\simeq \underline{h}_{f(x)}$$  
is an isomorphism in $\mathrm{Ho}(SPr_{\tau'}(T'))$.  
\end{prop}  
  
\textit{Proof:} This follows immediately by adjunction, from Corollary \ref{c6}.  
\hfill \textbf{$\Box$} 
  
\end{subsection}  
  
\begin{subsection}{Injective model structure and
stacks of morphisms}\label{morphisms}  
  
The goal of this paragraph is to present an injective version of
the local model structure on $SPr(T)$ for which 
cofibrations are monomorphisms, and to use it in order to 
construct \emph{stacks of morphisms}. Equivalently, we will
show that the injective model category of stacks over an
$S$-site possesses derived internal Hom's, and as a consequence
the homotopy category of stacks $\mathrm{Ho}(SPr_{\tau}(T))$
is \textit{cartesian closed} (in the usual sense of \cite[Ch. IV \S 10]{mcl}).  
These stacks of morphisms will be  
important especially for applications to Derived Algebraic Geometry (see \cite{camb, partII}),   
since many of the \textit{moduli stacks} are   
defined as stacks of morphisms to a certain \textit{classifying stack} (for example,   
the stack of vector bundles on a scheme).   
  
Before going into details, let us observe that in general, as explained in   
\cite[\S $11$]{sh}, the projective model structure on $SPr_{\tau}(T)$ is not an \textit{internal model category},    
i.e. is not a closed symmetric monoidal model category for the direct product (\cite[Def. $4.2.6$]{ho}),   
and therefore the internal $Hom$'s of the category $SPr_{\tau}(T)$ are not compatible  
with the model structure. This prevents one from defining   
derived internal $Hom$'s in the usual way (i.e. by applying the internal $Hom$'s of $SPr(T)$ to   
fibrant models for the targets and cofibrant models for the sources). One way to solve this problem is to work with another model category  
which is internal and Quillen equivalent to $SPr(T)$. 
The canonical choice is to use   
an \textit{injective model structure on $SPr(T)$}, analogous to the one described in \cite{ja}.   
  
\begin{prop}\label{pp1}  
Let $(T,\tau)$ be an $S$-site in $\mathbb{U}$. 
Then there exists a simplicial closed model structure  
on the category $SPr(T)$, called the \emph{local injective model structure} and denoted
by $SPr_{\mathrm{inj}, \tau}(T)$,   
where the cofibrations are the monomorphisms and the equivalences are the  
local equivalences. Moreover, the local injective model structure on   
$SPr(T)$ is proper and internal\footnote{  
Recall once again that a model category is said to be \textit{internal} if it is a monoidal model  
category (in the sense of \cite[Def. $4.2.6$]{ho}) for the monoidal structure given by the direct product.}.  
\end{prop}  
  
\textit{Proof:} The proof is essentially the same as
the proof of our Theorem \ref{t3}. The starting point is the objectwise
injective model structure $SPr_{\textrm{inj}}(T)$, for which 
equivalences and cofibrations are defined 
objectwise. The existence of this model structure 
can be proved by the same cardinality argument as in the case
where $T$ is a usual category (see \cite{ja}). The model 
category $SPr_{\mathrm{inj}}(T)$ is clearly proper, 
$\mathbb{U}$-cellular and $\mathbb{U}$-combinatorial, so one
can apply the localization techniques of \cite{hi}.
We define the model category $SPr_{\mathrm{inj},\tau}(T)$ as the left Bousfield
localization of $SPr_{\mathrm{inj}}(T)$ along the set of
hypercovers $H$ defined in the proof of Theorem \ref{t3}. 
Note that the identity functor 
$SPr_{\mathrm{inj},\tau}(T) \longrightarrow SPr_{\tau}(T)$
is the right adjoint of a Quillen equivalence. From this and
Theorem \ref{t3}
we deduce that equivalences in $SPr_{\mathrm{inj},\tau}(T)$ are exactly
the local equivalences of Definition \ref{d14}. This proves the
existence of the model category $SPr_{\mathrm{inj},\tau}(T)$.
The fact that it is proper follows  
easily from the fact the the model category $SSet$ is proper
and from the description of equivalences
in $SPr_{\mathrm{inj},\tau}(T)$ as $\pi_{*}$-equivalences. 
It only remains to show that   
$SPr_{\mathrm{inj}, \tau}(T)$ is internal. But, 
as cofibrations are the monomorphisms this follows easily from the fact that 
finite products preserves local equivalences. \hfill \textbf{$\Box$} \\  
  
As the equivalences in $SPr_{\mathrm{inj},\tau}(T)$ and $SPr_{\tau}(T)$ are the same,   
the corresponding homotopy categories coincide  
$$\mathrm{Ho}(SPr_{\textrm{inj}, \tau}(T))=\mathrm{Ho}(SPr_{\tau}(T)).$$  
Since the homotopy category of an internal model category is known to be cartesian closed,
 Proposition \ref{pp1} implies the following corollary.  
  
\begin{cor}\label{cc2}  
For any $S$-site $T$ in $\mathbb{U}$, the homotopy 
category of stacks $\mathrm{Ho}(SPr_{\tau}(T))$ is   
cartesian closed.  
\end{cor} 

\textit{Proof:} Apply  \cite[Thm. 4.3.2]{ho} to the symmetric
monoidal model category $SPr_{\mathrm{inj},\tau}(T)$, with 
the monoidal structure given by the direct product. \hfill $\Box$ \\
  
\begin{df}\label{dd1}  
\begin{enumerate}  
\item  
The internal $Hom$'s of the category $\mathrm{Ho}(SPr_{\tau}(T))$ will be denoted  
by   
$$\mathbb{R}_{\tau}\underline{\mathcal{H}om}(-,-) :   
\mathrm{Ho}(SPr_{\tau}(T)) \times \mathrm{Ho}(SPr_{\tau}(T)) \longrightarrow   
\mathrm{Ho}(SPr_{\tau}(T)).$$  
\item  
Let $(T,\tau)$ be an $S$-site in $\mathbb{U}$, and $F$, $G$ be stacks in $\mathrm{Ho}  
(SPr_{\tau}(T))$. The \emph{stack of morphisms}  
from $F$ to $G$ is defined to be the stack  
$$\mathbb{R}_{\tau}\underline{\mathcal{H}om}(F,G) \in \mathrm{Ho}(SPr_{\tau}(T)).$$  
\end{enumerate}  
\end{df}  

Explicitly, we have for any pair of stacks $F$ and $G$   
$$\mathbb{R}_{\tau}\underline{\mathcal{H}om}  
(F,G)\simeq \underline{\mathcal{H}om}(F,R_{\mathrm{inj}}G),$$  
where $R_{\mathrm{inj}}$ is the fibrant replacement functor in
the objectwise injective model category $SPr_{\mathrm{inj}}(T)$, and  
$\underline{\mathcal{H}om}$ is the internal $Hom$ functor of the category $SPr(T)$. In fact, if $G$ is a stack, then both $R_{\textrm{inj}}G$ and $\underline{\mathcal{H}om}(F,R_{\mathrm{inj}}G)$ are stacks.  \\

Actually, Proposition \ref{pp1} gives more than 
the cartesian closedness of $\mathrm{Ho}(SPr_{\tau}(T))$.
Indeed, one can consider the full sub-category
$SPr_{\mathrm{inj},\tau}(T)^{f}$ of fibrant objects
in $SPr_{\mathrm{inj},\tau}(T)$. As any object 
is cofibrant in $SPr_{\mathrm{inj},\tau}(T)$, 
for any two objects $F$ and $G$ in $SPr_{\mathrm{inj},\tau}(T)^{f}$
the internal Hom $\underline{\mathcal{H}om}(F,G)$ is
also a fibrant object and therefore lives in 
$SPr_{\mathrm{inj},\tau}(T)^{f}$. This shows in particular that
$SPr_{\mathrm{inj},\tau}(T)^{f}$ becomes cartesian 
closed for the direct product, and therefore one can 
associate to it a natural $SPr_{\mathrm{inj},\tau}(T)^{f}$-enriched
category $\underline{SPr_{\mathrm{inj},\tau}(T)^{f}}$. Precisely, 
the set of object of $\underline{SPr_{\mathrm{inj},\tau}(T)^{f}}$
is the set of fibrant objects in $SPr_{\mathrm{inj},\tau}(T)$, 
and for two such objects $F$ and $G$ the object of morphisms
is $\underline{\mathcal{H}om}(F,G)$. 

The $SPr_{\mathrm{inj},\tau}(T)^{f}$-enriched category
$\underline{SPr_{\mathrm{inj},\tau}(T)^{f}}$ yields in fact 
a \emph{up-to-equivalence $SPr_{\mathrm{inj},\tau}(T)^{f}$-enrichement of the $S$-category
$LSPr_{\tau}(T)$}. Indeed, as $SPr_{\tau}(T)$ and 
$SPr_{\mathrm{inj},\tau}(T)$ has the same simplicial localizations
(because they are the same categories with the same
notion of equivalence), one has
a natural equivalence of $S$-categories
$$LSPr_{\tau}(T)=LSPr_{\mathrm{inj},\tau}(T)\simeq 
Int(SPr_{\mathrm{inj},\tau}(T)).$$
Recall that the $S$-category $Int(SPr_{\mathrm{inj},\tau}(T))$
consists of fibrant objects in $SPr_{\mathrm{inj},\tau}(T)$
and their simplicial Hom-sets. In other words
the $SSet$-enriched category $Int(SPr_{\mathrm{inj},\tau}(T))$
is obtained from the $SPr_{\mathrm{inj},\tau}(T)^{f}$-enriched category
$\underline{SPr_{\mathrm{inj},\tau}(T)^{f}}$
by applying the global section functor
$\Gamma : SPr_{\mathrm{inj},\tau}(T) \longrightarrow SSet$.
In conclusion, one has a triple
$$(LSPr_{\tau}(T),\underline{SPr_{\mathrm{inj},\tau}(T)^{f}},\alpha)$$
where $\alpha$ is an isomorphism in $\mathrm{Ho}(S-Cat)$ between
$LSPr_{\tau}(T)$ and 
the underlying $S$-category of $\underline{SPr_{\mathrm{inj},\tau}(T)^{f}}$.
This triple is what we refer to as an \emph{up-to-equivalence
$SPr_{\mathrm{inj},\tau}(T)^{f}$-enrichement of $LSPr_{\tau}(T)$}. 
For example, the $SPr_{\mathrm{inj},\tau}(T)^{f}$-enriched functor
$$\underline{\mathcal{H}om} : 
(\underline{SPr_{\mathrm{inj},\tau}(T)^{f}})^{op}\times
\underline{SPr_{\mathrm{inj},\tau}(T)^{f}} \longrightarrow \underline{SPr_{\mathrm{inj},\tau}(T)^{f}}$$
gives rise to a well defined morphism in $\mathrm{Ho}(S-cat)$
$$\mathbb{R}_{\tau}\underline{\mathcal{H}om} : LSPr_{\tau}(T)^{op}\times LSPr_{\tau}(T) \longrightarrow
LSPr_{\tau}(T),$$
lifting the internal Hom-structure on the homotopy
category $\mathrm{Ho}(SPr_{\tau}(T))$.

\begin{rmk}\emph{This last structure is at first sight more subtle than 
the cartesian closedness of the homotopy category
$Ho(SPr_{\tau}(T))$, as $\underline{SPr_{\mathrm{inj},\tau}(T)^{f}}$
encodes strictly associative and unital
compositions between stacks of morphisms, which 
are only described by $\mathrm{Ho}(SPr_{\tau}(T))$ as 
up-to-homotopy associative and unital 
compositions. This looks like comparing the
notions of simplicial monoids (i.e. monoids
in $SSet$) and 
up-to-homotopy simplicial monoids (i.e. monoids
in $\mathrm{Ho}(SSet)$), and the former is well known
to be the \emph{right notion}. However, we would like to 
mention that we think that the $S$-category alone
$LSPr_{\tau}(T) \in \mathrm{Ho}(S-Cat)$, together with the
fact that $\mathrm{Ho}(SPr_{\tau}(T))$ is cartesian closed,
completely determines its up-to-equivalence 
$SPr_{\mathrm{inj},\tau}(T)^{f}$-enrichement. In other words, 
the structure 
$$(LSPr_{\tau}(T),\underline{SPr_{\mathrm{inj},\tau}(T)^{f}},\alpha)$$
only depends, up to an adequate notion of equivalence, on the $S$-category $LSPr_{\tau}(T)$. Unfortunately, 
investigating this question would drive us way too far
from our purpose, as we think the right context to treat
it is the general theory of \emph{symmetric monoidal
$S$-categories}, as briefly exposed in \cite[\S 5.1]{to4}.} 
\end{rmk}  

\end{subsection}  
  
\begin{subsection}{Truncated stacks and truncation functors}\label{truncations}  
  
We start by recalling some very general definition of 
truncated objects in model categories. 

\begin{df}\label{dtruncmodel}
\begin{enumerate}  
\item Let $n\geq 0$.  
An object $x \in \mathrm{Ho}(M)$ is called $n$-\emph{truncated} if for  
any $y\in \mathrm{Ho}(M)$, the mapping space $Map_{M}(y,x) \in \mathrm{Ho}(SSet)$ is $n$-truncated.   
  
\item An object $x \in \mathrm{Ho}(M)$ is called \emph{truncated} if it is $n$-truncated for some integer  
$n\geq 0$.  
\end{enumerate}
\end{df}

Clearly, a simplicial set $X$ is $n$-truncated in the
sense above if and only it is $n$-truncated in the classical
sense (i.e. if for any base point $x \in X$, $\pi_{i}(X,x)=0$
for all $i>n$). \\

We now fix an $S$-site $(T,\tau)$ in $\mathbb{U}$, and we consider the correponding model category of stacks $SPr_{\tau}(T)$.  
  
\begin{df}\label{dd2}  
Let $n\geq 0$ be an integer.  
A morphism $f : F \longrightarrow G$ in $SPr_{\tau}(T)$ is a $\pi_{\leq n}$\emph{-equivalence} (or   
a \emph{local $n$-equivalence}) if   
the following two conditions are satisfied.  
\begin{enumerate}  
\item The induced morphism   
$\pi_{0}(F) \longrightarrow \pi_{0}(G)$  
is an isomorphism of sheaves on $\mathrm{Ho}(T)$.  
\item For any object $x \in Ob(T)$, any section $s \in \pi_{0}(F(x))$ and   
any integer $i$ such that $n\geq i>0$, the induced morphism  
$\pi_{i}(F,s) \longrightarrow \pi_{i}(G,f(s))$  
is an isomorphism of sheaves on $\mathrm{Ho}(T/x)$.   
\end{enumerate}  
\end{df}  
  
\begin{thm}\label{tt1}  
There exists a closed model structure on $SPr(T)$, called the $n$\emph{-truncated local projective model structure},   
for which the equivalences are the  
$\pi_{\leq n}$-equivalences and the cofibrations are the cofibrations for the projective model  
structure on $SPr(T)$. Furthermore the $n$-local projective model structure is  
$\mathbb{U}$-cofibrantly generated and proper.   
  
The category $SPr(T)$ together with its $n$-truncated local projective model   
structure will be denoted by $SPr_{\tau}^{\leq n}(T)$.  
\end{thm}  
  
\textit{Proof:}  The proof is essentially a corollary of Theorem \ref{t3}.   
Let $J$ (resp., $I$) be a $\mathbb{U}$-small set of   
generating trivial cofibrations (resp., generating cofibrations) for the model category $SPr_{\tau}(T)$.  
Let $J'$ be the set of morphisms $\partial \Delta^{i}\otimes \underline{h}_{x} \longrightarrow  
\Delta^{i}\otimes \underline{h}_{x}$, for all $i>n$ and all $x \in Ob(T)$. We define $J(n)=J\cup J'$.   
Finally, let $W(n)$ be the set of $\pi_{\leq n}$-equivalences. It is easy (and left to the reader) to   
prove that \cite[Thm. 2.1.19]{ho} can be applied to the sets $W(n)$, $I$ and $J(n)$. \hfill \textbf{$\Box$}  
  
\begin{cor}\label{cc3}  
The model category $SPr_{\tau}^{\leq n}(T)$ is the left Bousfield localization of  
$SPr_{\tau}(T)$ with respect to the morphisms $\partial \Delta^{i}\otimes \underline{h}_{x} \longrightarrow  
\Delta^{i}\otimes \underline{h}_{x}$, for all $i>n$ and all $x \in Ob(T)$.  
\end{cor}  
  
\textit{Proof:} This follows immediately from the explicit description of the set $J(n)$ of generating  
cofibrations given in the proof of Theorem \ref{tt1} above. \hfill \textbf{$\Box$} \\  
  
Note that corollaries \ref{c4} and \ref{cc3} also imply that $SPr_{\tau}^{\leq n}(T)$ is a left  
Bousfield localization of $SPr(T)$.   
  
For the next corollary, an object $F \in SPr(T)$ is called \emph{objectwise $n$-truncated} if   
for any $x \in Ob(T)$, the simplicial set $F(x)$ is $n$-truncated (i.e. for any base point $s \in F(x)_{0}$,   
one has $\pi_{i}(F(x),s)=0$ for $i>n$).  
  
\begin{cor}\label{cc4}  
An object $F \in SPr_{\tau}^{\leq n}(T)$ is fibrant if and only if it is 
objectwise fibrant, satisfies the  
hyperdescent condition (see Def. \ref{d15}) and is objectwise $n$-truncated.  
\end{cor}  
  
\textit{Proof:} This again follows formally from the explicit description of the set $J(n)$ of  
generating cofibrations given in the proof of Theorem \ref{tt1}. \hfill \textbf{$\Box$} \\  
  
From the previous corollaries we deduce that the identity functor  
$\mathrm{Id} : SPr_{\tau}(T) \longrightarrow SPr_{\tau}^{\leq n}(T)$  
is a left Quillen functor, which then induces an adjunction on the homotopy categories  
$$t_{\leq n}:=\mathbb{L}\mathrm{Id} : \mathrm{Ho}(SPr_{\tau}(T))   
\longrightarrow \mathrm{Ho}(SPr_{\tau}^{\leq n}(T)) \qquad   
\mathrm{Ho}(SPr_{\tau}(T)) \longleftarrow \mathrm{Ho}(SPr_{\tau}^{\leq n}(T)) : j_{n}:=\mathbb{R}\mathrm{Id}.$$  
Note however that  
the functor  
$$t_{\leq n} : \mathbb{L}\mathrm{Id} : \mathrm{Ho}(SPr_{\tau}(T))   
\longrightarrow \mathrm{Ho}(SPr_{\tau}^{\leq n}(T))$$  
does not preserves homotopy fibered products in general.  
Finally, $j_{n}$ is fully faithful and a characterization of its   
essential image is given in the following lemma.  
  
\begin{lem}\label{ll1}  
Let $F \in SPr_{\tau}(T)$ and $n\geq 0$. The following conditions are equivalent.  
\begin{enumerate}  
\item $F$ is an $n$-truncated object in the mdoel 
category $SPr_{\tau}(T)$ (in the sense of Definition \ref{dtruncmodel}).
  
\item For any $x \in Ob(T)$ and any base point $s \in F(x)$, one has  
$\pi_{i}(F,s)=0$ for any $i>n$.  
  
\item The adjunction morphism $F \longrightarrow j_{n}t_{\leq n}(F)$ is an isomorphism in $\mathrm{Ho}(SPr_{\tau}(T))$.  
  
\end{enumerate}  
\end{lem}  
  
\textit{Proof:} The three conditions are invariant under isomorphisms in   
$\mathrm{Ho}(SPr_{\tau}(T))$; we can therefore always assume  
that $F$ is fibrant in $SPr_{\tau}(T)$.   
  
To prove that $(1) \Rightarrow (2)$, it is enough to observe  
that $\mathbb{R}_{\tau}\underline{Hom}(\underline{h}_{x},F)\simeq F(x)$. Conversely, let us suppose that   
$(2)$ holds and let    
$j : F \longrightarrow RF$ be a fibrant replacement in $SPr_{\tau}^{\leq n}(T)$. The hypothesis on $F$   
and Corollary \ref{cc4} imply that  
$j$ is a $\pi_{*}$-equivalence, thus showing that we can assume $F$ to be fibrant in $SPr_{\tau}^{\leq n}(T)$, and  
by Corollary \ref{cc4} again, that $F$ can be also assumed to be objectwise $n$-truncated. In particular, the natural  
morphism $F^{\Delta^{i}} \longrightarrow F^{\partial \Delta^{i}}$ is an objectwise trivial fibration  
for any $i>n$. Therefore, one has for any $i>n$,  
$$\mathbb{R}_{\tau}\underline{Hom}(G,F)^{\mathbb{R} \partial \Delta^{i}}\simeq  
\mathbb{R}\underline{Hom}_{\tau}(G,F^{\partial \Delta^{i}}) \simeq   
\mathbb{R}\underline{Hom}_{\tau}(G,F^{\Delta^{i}})\simeq  
\mathbb{R}_{\tau}\underline{Hom}(G,F)^{\mathbb{R} \Delta^{i}}.$$  
This implies that $\mathbb{R}_{\tau}\underline{Hom}(G,F)$ is $n$-truncated   
for any $G \in SPr_{\tau}(T)$. This proves the equivalence between $(1)$ and $(2)$.  
  
For any $F \in \mathrm{Ho}(SPr_{\tau}(T))$, the adjunction morphism   
$F \longrightarrow j_{n}t_{\leq n}(F)$ is represented  
in $SPr(T)$ by a fibrant resolution $j : F \longrightarrow RF$ in the model category $SPr_{\tau}^{\leq n}(T)$.   
If $F$ satisfies condition $(2)$, we have already seen that $j$ is a $\pi_{*}$-equivalence, and therefore   
that $(3)$ is satisfied. Conversely, by Corollary \ref{cc4}, $RF$ always satisfies condition $(2)$ and then   
$(3)\Rightarrow (2)$. \hfill \textbf{$\Box$} \\  
  
In the rest of the paper we will systematically use Lemma \ref{ll1} and the functor  
$j_{n}$ to identify the homotopy category $\mathrm{Ho}(SPr_{\tau}^{\leq n}(T))$   
with the full subcategory of $\mathrm{Ho}(SPr_{\tau}(T))$  
consisting of $n$-truncated objects. We will therefore never specify the functor $j_{n}$. With this  
convention, the functor $t_{\leq n}$ becomes an endofunctor  
$$t_{\leq n} : \mathrm{Ho}(SPr_{\tau}(T)) \longrightarrow \mathrm{Ho}(SPr_{\tau}(T)),$$  
called the $n$\textit{-th truncation functor}. There is an adjunction morphism $\mathrm{Id}   
\longrightarrow t_{\leq n}$, and  
for any $F \in \mathrm{Ho}(SPr_{\tau}(T))$, the morphism $F \longrightarrow t_{\leq n}(F)$   
is universal among morphisms from $F$ to an $n$-truncated object. More precisely, for any   
$n$-truncated object $G \in \mathrm{Ho}(SPr_{\tau}(T))$,   
the natural morphism  
$$\mathbb{R}_{\tau}\underline{Hom}(t_{\leq n}(F),G) \longrightarrow \mathbb{R}_{\tau}\underline{Hom}(F,G)$$  
is an isomorphism in $\mathrm{Ho}(SSet)$.   
  
\begin{df}\label{trunc}  
The \emph{$n$-th truncation functor} is the functor previously defined  
$$t_{\leq n} : \mathrm{Ho}(SPr_{\tau}(T)) \longrightarrow \mathrm{Ho}(SPr_{\tau}(T)).$$  
The essential image of $t_{\leq n}$ is called the subcategory of \emph{$n$-truncated stacks}.  
\end{df}  
  
Note that the essential image of $t_{\leq n}$ is by constuction equivalent to  
the category $\mathrm{Ho}(SPr_{\tau}^{\leq n}(T))$. \\  
  
The following proposition gives a complete characterization of the category of  
$0$-truncated stacks and of the $0$-th truncation functor $t_{\leq 0}$.  
  
\begin{prop}\label{pp2}  
The functor $\pi^{pr}_{0} : SPr(T) \longrightarrow Pr(\mathrm{Ho}(T))$ induces an equivalence of categories  
$$\mathrm{Ho}(SPr_{\tau}^{\leq 0}(T)) \simeq Sh_{\tau}(\mathrm{Ho}(T))$$  
\noindent where $Sh_{\tau}(\mathrm{Ho}(T))$ denotes the category of sheaves   
of sets on the usual Grothendieck site $(\mathrm{Ho}(T),\tau)$.  
\end{prop}  
  
\textit{Proof:} Let us first suppose that the topology $\tau$ is trivial. In this case, we define  
a quasi-inverse functor as follows. By considering sets as constant simplicial sets,  
we obtain an embedding $Pr(\mathrm{Ho}(T)) \subset SPr(\mathrm{Ho}(T))$ that we compose with the pullback  
$p^{*} : SPr(\mathrm{Ho}(T)) \longrightarrow SPr(T)$ along the natural projection $p : T \longrightarrow \mathrm{Ho}(T)$.   
It is quite clear that $F \mapsto \pi_{0}^{pr}(F)$ and $F \mapsto p^{*}(F)$ induce two functors, inverse of   
each others  
$$\pi_{0}^{pr} : \mathrm{Ho}(SPr^{\leq 0}(T)) \simeq Pr(\mathrm{Ho}(T)) : p^{*}.$$  
In the general case, we use Corollary \ref{c4}. We need to show that a   
presheaf $F \in Pr(\mathrm{Ho}(T))$ is a sheaf for the topology $\tau$ if and only   
if the corresponding object $p^{*}(F)$ has the hyperdescent property. This  
last step is left to the reader as an exercice. \hfill \textbf{$\Box$} 
  
\begin{rmk}\label{nella}  
\emph{  
\begin{enumerate}  
\item The previous Proposition implies in particular that the homotopy   
category of stacks $\mathrm{Ho}(SPr_{\tau}(T))$ always contains the category   
of sheaves on the site $(\mathrm{Ho}(T),\tau)$ as the full subcategory   
of $0$-truncated objects. Again, we will not mention explicitly the functor $p^{*} : Sh_{\tau}(\mathrm{Ho}(T))   
\longrightarrow \mathrm{Ho}(SPr_{\tau}(T))$ and identify  
$Sh_{\tau}(\mathrm{Ho}(T))$ with the full subcategory of   
$\mathrm{Ho}(SPr_{\tau}(T))$ consisting of $0$-truncated   
objects.  
\item Proposition \ref{pp2} is actually just the 0-th stage of a   
series of similar results involving higher truncations.   
In fact Proposition \ref{pp2} can be generalized to a Quillen equivalence  
between $SPr_{\tau}^{\leq n}(T)$ and a certain model category of presheaves of $n$-groupoids   
on the $(n+1)$-category $t_{\leq n}(T)$ obtained from $T$ by applying the  
$n$-th fundamental groupoid functor to its simplicial sets of morphisms (see \cite[\S 2, p. 28]{sh}).   
We will not investigate these results further   
in this paper.   
\end{enumerate}  }  
\end{rmk}  
  
\end{subsection}  
  
\begin{subsection}{Model topoi}\label{modeltopoi}  
  
Let $M$ be any $\mathbb{U}$-cellular (\cite[\S 14.1]{hi}) or  
$\mathbb{U}$-combinatorial (\cite{sm}, \cite[Def. 2.1]{du2}) left proper model category  
(see also Appendix A).   
Let us recall from Theorem \ref{tb1} and \ref{tb2} that for any $\mathbb{U}$-set of morphisms  
$S$ in $M$, the left Bousfield localization $\mathrm{L}_{S}M$ exists. It is a model category, whose underlying category  is still $M$, whose cofibrations are those of $M$ and whose equivalences   
are the so-called $S$-local equivalences (\cite[\S 3.4]{hi}).   
A left Bousfield localization of $M$ is any model   
category of the form $\mathrm{L}_{S}M$, for a $\mathbb{U}$-small set $S$ of morphisms  
in $M$.   
  
The following definition is a slight modification of the 
a notion communicated to us by   
C. Rezk (see \cite{rez2}). It is a model categorical analog of the notion of  
topos defined as a reflexive subcategory of the category of presheaves with an exact localization functor  
(see for example, \cite[Ch. 20]{schu}).   
  
\begin{df}\label{d18}  
\begin{enumerate}  
\item   
If $T$ is an $S$-category, a \emph{left exact Bousfield localization} of $SPr(T)$ is a   
left Bousfield localization $\mathrm{L}_{S}SPr(T)$ of $SPr(T)$, such that   
the identity functor $\mathrm{Id} : SPr(T) \longrightarrow \mathrm{L}_{S}SPr(T)$  
preserves homotopy fiber products.   
\item  
A $\mathbb{U}$\emph{-model topos} is a model category in $\mathbb{V}$ which is Quillen   
equivalent to a left exact Bousfield localization   
of $SPr(T)$ for some $T\in S-Cat_{\mathbb{U}}$.   
\end{enumerate}  
\end{df}  
  
For $2.$, recall our convention throughout the paper, according to which two   
model categories are Quillen equivalent if they can be connected by a finite \textit{chain} of   
Quillen equivalences, regardless of their direction.  
\noindent We will also need the 
following general definitions related to the notion of 
truncated objects in a model category 
(see Remark \ref{lurie} for some comments on it).   
  
\begin{df}\label{d19}  
Let $M$ be any model category.  
\begin{enumerate}
\item We say that $M$ is $t$\emph{-complete} if truncated objects detect   
isomorphisms in $\mathrm{Ho}(M)$ i.e. if a morphism $u : a \rightarrow b$  
in $\mathrm{Ho}(M)$ is an isomorphism if and only if, for any truncated object $x$ in $\mathrm{Ho}(M)$, the map   
$u^{*} : [b,x] \longrightarrow [a,x]$ is bijective.  
  
\item A $\mathbb{U}$-model topos is $t$\emph{-complete} if its underlying model category is $t$-complete.  
  
\end{enumerate}  
\end{df}  
   
The next Theorem shows that given an $S$-category $T$,     
$t$-complete left exact Bousfield localizations of $SPr(T)$ correspond exactly to simplicial topologies on $T$.  
It should be considered as a homotopy analog of the correspondence for usual Grothendieck   
topologies as described e.g. in \cite[Thm. 20.3.7]{schu}.   
  
\begin{thm}\label{t4}  
Let $T$ be a $\mathbb{U}$-small $S$-category. There exists a bijective correspondence between  
$S$-topologies on $T$ and left exact Bousfield localizations of $SPr(T)$ which are $t$-complete.  
\end{thm}  
  
\textit{Proof:} Let $\mathcal{T}(T)$ be the set of $S$-topologies on $T$, which by definition is also the set of  
Grothendieck topologies on $\mathrm{Ho}(T)$. Let $\mathcal{B}(T)$ be the set   
of left exact Bousfield localizations of $SPr(T)$, and   
$\mathcal{B}_{t}(T) \subset \mathcal{B}(T)$ the subset of those which are $t$-complete.   
We are first going to define maps $\phi:\mathcal{T}(T) \rightarrow \mathcal{B}_{t}(T)$   
and $\psi:\mathcal{B}_{t}(T) \rightarrow \mathcal{T}(T)$. \\ 
    
\begin{center} \textit{The map $\phi:\mathcal{T}(T) \rightarrow \mathcal{B}_{t}(T)$} \end{center}  
  
Let $\tau \in \mathcal{T}(T)$ be an $S$-topology on $T$. According to Corollary \ref{c4} and Proposition   
\ref{p3} $(2)$, $SPr_{\tau}(T)$ is a left exact Bousfield localization of $SPr(T)$. We are going to show that  
$SPr_{\tau}(T)$ is also $t$-complete.  
We know by Lemma \ref{ll1},  
that an object $F \in \mathrm{Ho}(SPr_{\tau}(T))$ is $n$-truncated if and only if $F\simeq t_{\leq n}(F)$.   
Therefore, if a morphism $f : F \longrightarrow G$ satisfies the condition $(3)$ of Definition \ref{d19},  
one has   
$$[t_{\leq n}(F),H] \simeq [F,H] \simeq  
[G,H] \simeq [t_{\leq n}(G),H]$$  
for any $n$-truncated object $H \in \mathrm{Ho}(SPr_{\tau}(T))$.  
This implies that for any $n$, the induced morphism $t_{\leq n}(F) \longrightarrow t_{\leq n}(G)$ is  
an isomorphism in $\mathrm{Ho}(SPr_{\tau}^{\leq n}(T))$, and hence  
in $\mathrm{Ho}(SPr_{\tau}(T))$. In other words, $f$ is an $\pi_{\leq n}$-equivalence for any   
$n$, and hence a $\pi_{*}$-equivalence. This shows that the model category  
$SPr_{\tau}(T)$ is a $t$-complete model category and allows us  
to define the map $\phi : \mathcal{T}(T) \longrightarrow \mathcal{B}_{t}(T)$   
by the formula $\phi(\tau)=SPr_{\tau}(T)$. \\  
  
\begin{center} \textit{The map $\psi:\mathcal{B}_{t}(T) \rightarrow \mathcal{T}(T)$} \end{center}  
  
Let $\mathrm{L}_{S}SPr(T) \in \mathcal{B}_{t}(T)$, and let us consider the derived Quillen adjunction given by the  
identity functor $\mathrm{Id}  : SPr(T)  \longrightarrow \mathrm{L}_{S}SPr(T)$  
$$a:=\mathbb{L}\mathrm{Id} :   
\mathrm{Ho}(SPr(T)) \longrightarrow \mathrm{Ho}(\mathrm{L}_{S}SPr(T))   
\qquad \mathrm{Ho}(SPr(T)) \longleftarrow \mathrm{Ho}(\mathrm{L}_{S}SPr(T)) : \mathbb{R}\mathrm{Id}=:i.$$  
The reader should note that the above functor $a$ is not equal a priori to the associated stack functor of   
Definition \ref{d16} $(5)$, as no $S$-topology on $T$ has been given yet.   
We know that $j$ is fully faithful and identifies $\mathrm{Ho}(\mathrm{L}_{S}SPr(T))$   
with the full subcategory of  
$\mathrm{Ho}(SPr(T))$ consisting of $S$-local objects (see \cite[Def. 3.2.4 1(a); Th. 4.1.1 (2)]{hi}).   
  
We consider the full subcategory $\mathrm{Ho}_{\leq 0}(\mathrm{L}_{S}SPr(T))$ (resp.   
$\mathrm{Ho}_{\leq 0}(SPr(T))$) of $\mathrm{Ho}(\mathrm{L}_{S}SPr(T))$ (resp.   
of $\mathrm{Ho}(SPr(T))$) consisting of $0$-truncated objects. Note that in   
general, an object $x$ in a model category  
is $0$-truncated if and only if for any $n\geq1$, the natural morphism $x^{\mathbb{R}\Delta^{n}}   
\longrightarrow x^{\mathbb{R}\partial\Delta^{n}}$  
is an equivalence. As both $a$ and $i$ preserve homotopy fiber products,   
they also preserve $0$-truncated objects. Therefore we have an   
induced adjunction  
$$a_{0} : \mathrm{Ho}_{\leq 0}(SPr(T)) \longrightarrow \mathrm{Ho}_{\leq 0}(\mathrm{L}_{S}SPr(T))   
\qquad \mathrm{Ho}_{\leq 0}(SPr(T)) \longleftarrow \mathrm{Ho}_{\leq 0}(\mathrm{L}_{S}SPr(T)) : i_{0}.$$  
Now, the functor $\pi^{pr}_{0} : \mathrm{Ho}(SPr(T)) \longrightarrow Set^{\mathrm{Ho}(T)^{op}}$   
induces an equivalence of categories  
$$\mathrm{Ho}_{\leq 0}(SPr(T)) \simeq Set^{\mathrm{Ho}(T)^{op}}=:Pr(\mathrm{Ho}(T)),$$  
and so the adjunction $(a_{0},i_{0})$ is in fact equivalent to an adjunction  
$$a_{0} : Pr(\mathrm{Ho}(T)) \longrightarrow \mathrm{Ho}_{\leq 0}(\mathrm{L}_{S}SPr(T))   
\qquad Pr(\mathrm{Ho}(T)) \longleftarrow \mathrm{Ho}_{\leq 0}(\mathrm{L}_{S}SPr(T)) : i_{0}$$  
where, of course, the functor $i_{0}$ is still fully faithful and the functor $a_{0}$ is exact.   
By \cite[Thm. 20.3.7]{schu}, there exists then a unique  
Grothendieck topology $\tau$ on $\mathrm{Ho}(T)$ such that the essential image of $i_{0}$ is exactly   
the full subcategory of sheaves on $\mathrm{Ho}(T)$ for the  
topology $\tau$. The functor $a_{0}$ is then equivalent to the asscociated sheaf functor. Thus, we define  
$\psi : \mathcal{B}_{t}(T) \longrightarrow \mathcal{T}(T)$ by the formula $\psi(\mathrm{L}_{S}SPr(T)):=  
\tau \in \mathcal{T}(T)$. \\  
    
\begin{center} \textit{Proof of } $\phi \circ \psi=\mathrm{Id}$ \end{center}  
   
Let $\mathrm{L}_{S}SPr(T) \in \mathcal{B}_{t}(T)$ be a left exact Bousfield localization of $SPr(T)$ and  
$\tau=\psi(\mathrm{L}_{S}SPr(T))$ the corresponding topology on $T$. We need to prove that the set  
of $S$-local equivalences equal the set of $\pi_{*}$-equivalences. Recall that we have denoted by  
$$a:=\mathbb{L}\mathrm{Id} :   
\mathrm{Ho}(SPr(T)) \longrightarrow \mathrm{Ho}(\mathrm{L}_{S}SPr(T))   
\qquad \mathrm{Ho}(SPr(T)) \longleftarrow \mathrm{Ho}(\mathrm{L}_{S}SPr(T)) : \mathbb{R}\mathrm{Id}=:i,$$  
the adjunction induced by the identity functor $\mathrm{Id} : L_{S}SPr(T) \longrightarrow SPr(T)$.  
  
Let us first prove that $S$-local equivalences are $\pi_{*}$-equivalences. Equivalently, we need to prove that for any   
morphism $f : F \longrightarrow G$ which is an equivalence in $\mathrm{L}_{S}SPr(T)$,   
$f$ is an hypercover in $SPr_{\tau}(T)$.   
For this we may assume that $F$ and $G$ are both objectwise fibrant objects.   
As the identity functor $\mathrm{Id} : SPr(T) \longrightarrow SPr_{\tau}(T)$ preserves   
homotopy fiber products, the induced morphism  
$$F^{\Delta^{n}} \longrightarrow F^{\partial\Delta^{n}}\times_{G^{\partial \Delta^{n}}}G^{\Delta^{n}}$$  
is still an $S$-local equivalence. Using this fact and Lemma \ref{l3}, one sees that it is enough to show that   
$f$ is a covering in $SPr_{\tau}(T)$.   
  
Recall that the topology $\tau$ is defined in such a way that   
the associated sheaf to a presheaf of sets $E$ on $\mathrm{Ho}(T)$  
is $i_{0}a_{0}(E)$ (where the adjunction $(a_{0},i_{0})$ is the one considered above in the definition of the map $\psi$).   
It is therefore enough to prove that   
the induced morphism $a_{0}(\pi_{0}^{pr}(F)) \longrightarrow a_{0}(\pi_{0}^{pr}(G))$  
is an isomorphism\footnote{Recall that $\pi_{0}^{pr}(F)$ is a presheaf of sets  
on $\mathrm{Ho}(T)$, that is considered via the projection $p : T \longrightarrow \mathrm{Ho}(T)$  
as a presheaf of discrete simplicial sets on $T$, and therefore as an object  
in $SPr(T)$.}.  
  
\begin{lem}\label{tech2}  
For any $F  \in \mathrm{Ho}(SPr(T))$, one has  
$$a_{0}(\pi_{0}^{pr}(F))\simeq a_{0}\pi_{0}^{pr}(ia(F)).$$  
\end{lem}  
  
\textit{Proof:} This immediately follows from the adjunctions $(a,i)$ and $(a_{0},i_{0})$, and the fact that  
$\pi_{0}^{pr}$ is isomorphic to the $0$-th truncation functor $t_{\leq 0}$  
on $\mathrm{Ho}(SPr(T))$. \hfill $\Box$ \\
  
As $f$ is an $S$-local equivalence, the morphism   
$ia(F) \longrightarrow ia(G)$ is an isomorphism in $\mathrm{Ho}(SPr(T))$,   
and therefore the same is true for  
$$a_{0}(\pi_{0}^{pr}(F))\simeq a_{0}\pi_{0}^{pr}(ia(F))\longrightarrow   
a_{0}\pi_{0}^{pr}(ia(G))\simeq a_{0}(\pi_{0}^{pr}(G)).$$  
  
We have thus shown that the $S$-local equivalences are $\pi_{*}$-equivalences.   
Conversely, to show that $\pi_{*}$-equivalences  
are $S$-local equivalences it is enough to show that for any $x \in Ob(T)$ and any   
hypercover $F_{*} \longrightarrow \underline{h}_{x}$  
in $SPr_{\tau}(T)$, the natural morphism   
$$ia(|F_{*}|) \longrightarrow ia(\underline{h}_{x})$$  
is an isomorphism in $\mathrm{Ho}(SPr(T))$ (see Corollary \ref{c4}).  As $a$ preserves   
homotopy fibered products, one has  
$(ia(G))^{\mathbb{R}K}\simeq ia(G^{\mathbb{R}K})$,   
for any $G \in \mathrm{Ho}(SPr(T))$ and any finite simplicial set $K$ (here  
$(-)^{\mathbb{R}K}$ is computed in the model category $SPr(T)$). Therefore, by $t$-completeness one has, for any $n$  
$$t_{\leq n-1}(ia(|F_{*}|))\simeq t_{\leq n-1}(ia(|\mathbb{R}Cosk_{n}F_{*}|)).$$  
This shows that one can assume that $F_{*}=\mathbb{R}Cosk_{n}(F_{*}/\underline{h}_{x})$, for some $n$  
(i.e. that $F_{*} \longrightarrow \underline{h}_{x}$ is relatively $n$-bounded). Furthermore, the same argument as   
in the proof of Theorem \ref{t3}, but relative to $\underline{h}_{x}$,   
shows that, by induction, one can assume $n=0$. In other words, one can assume  
that $F_{*}$ is the derived nerve of a covering $F_{0} \longrightarrow \underline{h}_{x}$ (which will be assumed to be  
an objectwise fibration).   
  
By the left exactness property of $a$ and $i$, the object $ia(|F_{*}|)$ is isomorphic in $\mathrm{Ho}(SPr(T))$ to   
the geometric realization of the derived nerve of $ia(F_{0}) \longrightarrow ia(\underline{h}_{x})$. This implies that   
for any $y \in Ob(T)$, the morphism $ia(|F_{*}|)(y) \longrightarrow ia(\underline{h}_{x})(y)$  
is isomorphic in $\mathrm{Ho}(SSet)$ to the geometric realization of the   
nerve of a fibration between simplicial sets. It is well  
known that such a morphism is isomorphic in $\mathrm{Ho}(SSet)$ to an   
inclusion of connected components. Therefore it is enough to show that the morphism    
$$\pi_{0}^{pr}(ia(|F_{*}|)) \longrightarrow \pi_{0}^{pr}(ia(\underline{h}_{x}))$$  
induces an isomorphism on the associated sheaves. By Lemma \ref{tech2}, this is equivalent to showing that the morphism  
$$i_{0}a_{0}\pi_{0}^{pr}(ia(|F_{*}|)) \longrightarrow i_{0}a_{0}\pi_{0}^{pr}(ia(\underline{h}_{x}))$$  
is an isomorphism of presheaves of sets on $\mathrm{Ho}(T)$.  
This morphism is also isomorphic to  
$$i_{0}a_{0}(\pi_{0}^{pr}(|F_{*}|)) \longrightarrow i_{0}a_{0}\pi_{0}^{pr}(\underline{h}_{x})$$  
whose left hand side is the sheaf associated to the co-equalizer of the two projections  
$$pr_{1},pr_{2} : \pi_{0}^{pr}(F_{0})\times_{\pi_{0}^{pr}(\underline{h}_{x})}\pi_{0}^{pr}(F_{0}) \longrightarrow  
\pi_{0}^{pr}(\underline{h}_{x}),$$  
whereas the right hand side is the sheaf associated to $\pi_{0}^{pr}(\underline{h}_{x})$. To conclude the proof, it is  
enough to notice that   
$\pi_{0}^{pr}(F_{0}) \longrightarrow \pi_{0}^{pr}(\underline{h}_{x})$ induces an epimorphism of sheaves (because  
$F_{*}$ is a hypercover) and that epimorphisms of sheaves are always effective (see \cite[Exp. II, Th\'eor\`eme 4.8]{sga4}).\\    
  
\begin{center} \textit{Proof of } $\psi \circ \phi=\mathrm{Id}$ \end{center}  
   
Let $\tau$ be a topology on $T$.   
By definition of the maps $\psi$ and $\phi$, to prove that $\psi \circ \phi=\mathrm{Id}$,   
it is equivalent to show that the functor $\pi^{pr}_{0} : \mathrm{Ho}(SPr_{\tau}(T)) \longrightarrow  
Pr(\mathrm{Ho}(T))$, when restricted to the full subcategory of $0$-truncated objects in   
$\mathrm{Ho}(SPr_{\tau}(T))$, induces an equivalence to the category of sheaves on the   
site $(\mathrm{Ho}(T),\tau)$. But this follows from Proposition \ref{pp2}. \hfill \textbf{$\Box$} 
  
\begin{cor}\label{c7}  
Let $M$ be a model category in $\mathbb{U}$.   
The following conditions are equivalent:  
\begin{enumerate}  
\item The model category $M$ is a $t$-complete $\mathbb{U}$-model topos.  
  
\item The model category $M$ is $t$-complete and there exists   
a $\mathbb{U}$-small category $C$ and a subcategory $S \subset C$, such that   
$M$ is Quillen equivalent to a left exact Bousfield localization of $M^{C,S}$ (see Def. \ref{d7}).  
  
\item There exists a $\mathbb{U}$-small $S$-site $(T,\tau)$ such that $M$ is Quillen equivalent to   
$SPr_{\tau}(T)$.  
  
\end{enumerate}  
\end{cor}  
  
\textit{Proof:} The equivalence of $(2)$ and $(3)$ follows immediately   
from Theorem \ref{t2} and the \textit{delocalization theorem} \cite[Thm. $2.5$]{dk2}, while $(1)$ and $(3)$ are equivalent   
by Theorem \ref{t4}. \hfill \textbf{$\Box$} \\  
  
The previous results imply in particular the following interesting rigidity property for $S$-groupoids.  
  
\begin{cor}  
Let $T$ be a $\mathbb{U}$-small $S$-category such that $\mathrm{Ho}(T)$ is a groupoid (i.e.  
every morphism in $T$ is invertible up to homotopy). Then, there is no non-trivial  
$t$-complete left exact Bousfield localization of $SPr(T)$.  
\end{cor}  
  
\textit{Proof:} In fact, there is no non-trivial topology on a groupoid, and therefore there  
is no non-trivial $S$-topology on $T$. \hfill $\Box$ 
  
\begin{rmk} \label{lurie}   
\begin{enumerate}  
\item \emph{There exist $t$-complete $\mathbb{U}$-model topoi which are not Quillen   
equivalent to some $SPr_{\tau}(T)$, for $T$ a $\mathbb{U}$-small category. Indeed,   
when $T$ is a category, the model category $SPr_{\tau}(T)$ is such that   
any object is a homotopy colimits of $0$-truncated objects (this is because representable  
objects are $0$-truncated). It is not difficult to see that this last property  
is not satisfied when $T$ is a general $S$-category.   
For example, let  
$T=BK(\mathbb{Z},1)$ be the $S$-category with a unique object and   
the simplicial monoid $K(\mathbb{Z},1)$ as simplicial set of endomorphisms.  
Then, $SPr(T)$ is the model category of simplicial sets together with an   
action of $K(\mathbb{Z},1)$, and $0$-truncated objects in $SPr(T)$ are  
all equivalent to discrete simplicial set with a trivial action of $K(\mathbb{Z},1)$.  
Therefore any homotopy colimit of such will be a simplicial set with a trivial  
action by $K(\mathbb{Z},1)$. However, the action of $K(\mathbb{Z},1)$ on itself  
by left translations is \textit{not} equivalent to a trivial one. } 
\item \emph{As observed by J. Lurie, there are examples of   
left exact Bousfield localization of $SPr(T)$ which are \textit{not}  
of the form $SPr_{\tau}(T)$. To see this, let $(T,\tau)$ be   
a Grothendieck site and consider the left Bousfield localization $\mathrm{L}_{cov}SPr(T)$ of $SPr(T)$ along   
only those hypercovers which are nerves of coverings (obviously, not all   
hypercovers are of this kind). Now, an example due to   
C. Simpson shows that there are Grothendieck sites $(T,\tau)$ such that $\mathrm{L}_{cov}SPr(T)$   
is not the same as $SPr_{\tau}(T)$ (see for example \cite[Ex. $A.10$]{dhi}). However,   
$\mathrm{L}_{cov}SPr(T)$ is a left exact Bousfield localization of $SPr(T)$, and the topology it   
induces on $T$ via the procedure used   
in the proof of Theorem \ref{t4}, coincides with $\tau$. Of course, the point   
here is that $\mathrm{L}_{cov}SPr(T)$ is \textit{not} a   
$t$-complete model category. This shows that one can not omit the hypothesis of   
$t$-completeness in Theorem \ref{t4}.  } 
\item \emph{Though the hypothesis of $t$-completeness in Theorem \ref{t4} is quite natural, and allows for a clean explanation in terms of $S$-topologies, it could be interesting to look for a similar   
comparison result without such an assumption. One way to proceed would be to introduce a notion of   
\textit{hyper-topology} on a  
category (or more generally on an $S$-category), a notion which was suggested to us by some independent remarks of V. Hinich, A. Joyal and C. Simpson.   
A hyper-topology on a category would be essentially the same thing as a topology with the difference   
that one specifies directly the hypercovers  
and not only the coverings; the conditions it should satisfy are analogous to  
the conditions imposed on the family of coverings in the usual definition of a Grothendieck (pre)topology.   
The main point here is that for a given Grothendieck site $(T,\tau)$, the two hyper-topologies  
defined using \textit{all} $\tau$-hypercovers on one side or only \textit{bounded} $\tau$\textit{-coverings}   
on the other side, will not be  equivalent in general.  
It seems reasonable to us that our Theorem \ref{t4} can be generalized to a correspondence between hyper-topologies   
on $T$ and arbitrary left exact Bousfield localizations of $SPr(T)$. This notion of hypertopology seems to be closely related to Cisinki's results in \cite{cis}.}  
\item \emph{Theorem \ref{t4} suggests also a way to think of \textit{higher topologies} on $n$-categories   
(and of \textit{higher topoi}) for $n \geq 1$ as appropriate \textit{left exact localizations}. In this case,   
the explicit notion of higher topology (that one has to reconstruct e.g. assuming the Theorem still holds   
for higher categories), will obviously depend on more then the associated homotopy category.   
For example, for the case of $2$-categories, as opposed to the case when all $i$-morphisms are invertible for $i>1$   
(see Remark \ref{infty}), a topology should give rise to some kind of topologies   
on the various categories of $1$-morphisms and these topologies should   
satisfy some compatibility condition with respect to   
the composition.}  
\end{enumerate}  
\end{rmk}  
  
\medskip  
  
We finish this paragraph with the following definition.  
  
\begin{df}\label{dstop}  
An $\mathbb{U}-S$\emph{-topos} is an $S$-category which is isomorphic in $\mathrm{Ho}(S-Cat)$ to   
some $LSPr_{\tau}(T)$, for $(T,\tau)$ a $\mathbb{U}$-small $S$-site.   
\end{df}  
  
\end{subsection}

\end{section}  

\begin{section}{Stacks over pseudo-model categories}\label{stacksonpseudo}  
  
In this Section we define the notion of a \textit{model pre-topology} on a model category and the   
notion of \textit{stacks} on such \textit{model sites}. A model pre-topology is a homotopy variation   
of the usual notion of a Grothendieck pre-topology and it reduces to the latter when the model   
structure is trivial (i.e. when equivalences are isomorphisms and any map is a fibration and a   
cofibration). We develop the theory in the slightly more general context of \textit{pseudo-model   
categories}, i.e of full subactegories of model categories that are closed under equivalences   
and homotopy pull-backs (see Definition \ref{pseudo}). We have chosen to work in this more general   
context because in some applications we will need to   
use subactegories of model categories defined by \textit{homotopy invariant conditions} but not   
necessarily closed under small limits and/or colimits (e.g., certain  
subcategories of objects of \textit{finite presentation}). The reader is however strongly encouraged to  
cancel everywhere the word \textit{pseudo-} in the following and to restore it only when interested in   
some application that requires such a degree of generality (as for example, the problem of defining   
\'etale $K$-theory on the pseudo-model category of connective commutative $\mathbb{S}$-algebras, see Prop.   
\ref{pseudoconn}). On the other hand, the theory itself presents no   
additional difficulty, except possibly for the linguistic one.   
  
\begin{subsection}{Model categories of pre-stacks on a pseudo-model category}\label{Mprestacks}  
  
In this subsection we will define the (model)category of pre-stacks on a \textit{pseudo-model}   
category which is essentially a category with weak equivalences   
that admits a nice embedding into a model category.   
  
\begin{df} \label{pseudo}  
A $\mathbb{U}$-small \emph{pseudo-model category} is a triple   
$(C,S,\iota)$ where $C$ is a $\mathbb{U}$-small category, $S\subset C$ is a  
subcategory of $C$ and   
$\iota : C \rightarrow M$ is a functor to a model $\mathbb{U}$-category $M$ satisfying the following  
four conditions.  
  
\begin{enumerate}  
\item The functor $\iota$ is fully faithful.  
\item One has $\iota(S)=W\cap \iota(C)$, where $W$ is the set of weak equivalences   
in the model category $M$.  
\item The category $C$ is closed under equivalences in $M$, i.e. if $x \rightarrow y$ is   
an equivalence in $M$ and $x$ (resp. $y$) is in the image of $\iota$,   
then so is $y$ (resp. $x$).  
\item The category $C$ is closed under homotopy pullbacks in $M$.  
\end{enumerate}  
  
\noindent The localization $S^{-1}C$ will be called the \emph{homotopy category} of $(C,S)$   
and often denoted by $\mathrm{Ho}(C,S)$ or simply $\mathrm{Ho}(C)$ when the choice of $S$ is unambiguous.  
  
\end{df}  
  
Condition $(4)$ of the previous definition can be precised as follows. Denoting by   
$\mathrm{Ho}(\iota) : S^{-1}C \rightarrow \mathrm{Ho}(M)$ the functor   
induced by $\iota$ (due to $(2)$.), which is fully faithful due to $(1)$ and $(3)$,   
the image of $\mathrm{Ho}(\iota)$, that coincides with its essential immage, is closed under homotopy pullbacks. 
  
Note also that because of  
condition $(3)$ of Definition \ref{pseudo}, the functor $\iota$ is an isomorphism  
from $C$ to its essential image in $M$.  
Hence we will most of the time identify $C$ with its image $\iota(C)$ in the model   
category $M$; therefore an object $x\in C$ will be called \textit{fibrant} (respectively,   
\textit{cofibrant}) in $C$ if $\iota(x)$ is fibrant (resp. cofibrant) in $M$. Moreover,   
we will sometimes call the maps in $S$ simply equivalences.   
  
Conditions $(3)$ and $(4)$ imply in particular that for any diagram  
$$\xymatrix{x \ar[r]^{p} & y \\  
 & z \ar[u]}$$  
of fibrant objects in $C$, such that $p$ is a fibration, the fibered product  
$x\times_{z}y$ exists. Indeed, this fibered product exists in   
the ambient model category $M$, and being equivalent to the   
homotopy fibered product, it also belongs to $C$ by condition $(3)$ and $(4)$.  
  
\begin{rmk}\label{trivia}  
\emph{  
\begin{enumerate}  
\item Being a pseudo-model category is not a self-dual property, in the sense   
that if $M$ is a pseuo-model category, then   
$M^{op}$ is not pseudo-model in general. Objects satisfying Definition \ref{pseudo} should be called more  
correctly \textit{right pseudo-model categories} and the dual definition (i.e. closure by   
homotopy push-outs) should deserve the name of \textit{left pseudo-model category}. However, to   
simplify the terminology, we fix once for all definition \ref{pseudo} as it is stated.  
\item Note that if $M$ is a model category with weak equivalences   
$W$, the triple $(M,W,\mathrm{Id}_{M})$ is a pseudo-model category. Moreover,   
a pseudo-model category is \textit{essentially} a model category. In fact, conditions   
$(1)-(3)$ imply that $C$ satisfies conditions $(1)$, $(2)$ and $(4)$ of   
the definition of a \textit{model structure} in the sense of \cite[Def. 1.1.3]{ho}.   
However, $C$ is not exactly a model category in   
general, since it is not required to be complete and co-complete (see \cite[Def. 1.1.4]{ho}),  
and the lifting property $(3)$ of \cite[Def. 1.1.3]{ho} is not necessarily satisfied.  
\item If $C$ is a complete and co-complete category and $S$ consists of all isomorphisms   
in $C$, then $(C,S, \mathrm{Id}_{C})$ is a \textit{trivial} pseudo-model category, where   
we consider on $C$ the trivial model structure with equivalences consisting of all   
isomorphisms and any map being a fibration (and a cofibration). If $C$ is not necessarily complete and co-complete but has finite limits,   
then we may view it as a \textit{trivial} pseudo-model category by replacing it with its   
essential image in $Pr(C)$ or $SPr(C)$, endowed with the trivial model stuctures,   
and taking $S$ to be all the isomorphisms.    
\end{enumerate}  
}  
\end{rmk}  
  
\begin{ex}\label{pseudoex}  
\emph{  
\begin{enumerate}  
\item Let $k$ be a commutative ring and $M:=Ch(k)^{op}$ the opposite model category of   
unbounded chain complexes of $k$-modules (see \cite[Def. 2.3.3]{ho}). The full subcategory   
$C \hookrightarrow M$ of \textit{homologically positive} objects (i.e. objects $P_{\bullet}$   
such that $H_{i}(P_{\bullet})=0$ for $i<0$) is a pseudo-model category.  
\item Let $k$ be a commutative ring (respectively, a field of characteristic zero) and let   
$M:=(E_{\infty}-\mathrm{Alg}_{k})^{op}$ (respectively, $M=\mathrm{CDGA}_{k}^{op}$) be the   
opposite model category of $E_{\infty}$-algebras over the   
category of unbounded cochain complexes of $k$-modules (resp., the  
opposite model category of commutative and unital differential  
graded $k$-algebras in non-positive degrees) which belong to  
$\mathbb{U}$ (see for example \cite{hin} for a description of  
these model structures). We say that an object $A$ of $M$ is \textit{finitely presented} if for any   
filtered direct diagram $C : J \rightarrow M^{op}$, with $J \in \mathbb{U}$, the natural map   
$$  
\mathrm{hocolim}_{j \in J}\mathrm{Map}_{M^{op}}(A,C_{j})   
\longrightarrow \mathrm{Map}_{M^{op}}(A, \mathrm{hocolim}_{j \in J}C_{j})  
$$  
is an equivalence of simplicial sets. Here $\mathrm{Map}_{M^{op}}(-,-)$ denotes the   
mapping spaces (or function complexes) in the model category $M^{op}$ (see \cite[\S $5.4$]{ho}).   
The reader will check that   
the full subcategory $C \hookrightarrow M$   
of finitely presented objects is a pseudo-model category.  
\item Let $A$ be a commutative $\mathbb{S}$-algebra as defined in \cite[Ch. 2, \S 3]{ekmm}.   
Let $M$ be the opposite category of the comma model category of commutative $\mathbb{S}$-algebras   
under $A$: an object in $M$ is then a map of commutative $\mathbb{S}$-algebras $A \rightarrow B$.   
Then, the full subcategory $C \hookrightarrow M$ consisting of finitely presented   
$A$-algebras (see the previous example or Definition \ref{fp}) is a pseudo-model category.  The full subcategory $C \hookrightarrow M$ consisting of \textit{\'etale maps} $A \rightarrow B$ (see Definition \ref{etdf}) is also a pseudo-model category. This pseudo-model category will be called the \textit{small \'etale site} over $A$.   
\item Let $X$ be a scheme and $C(X,\mathcal{O})$ be the category of unbounded cochain complexes of  
$\mathcal{O}$-modules. There exists a model structure on $C(X,\mathcal{O})$ where the equivalences are  
the local quasi-isomorphisms. Then, the full subcategory of $C(X,\mathcal{O})$ consisting of  
\textit{perfect complexes} is a pseudo-model category.  
\end{enumerate}  
}  
\end{ex}  
  
Recall from Subsection \ref{resdia} that for any category $C$ in $\mathbb{U}$ and any subcategory   
$S \subset C$, we have defined (Definition \ref{d7}) the model category $SSet_{\mathbb{U}}^{C,S}$ of   
\textit{restriced diagrams} on $(C,S)$ of simplicial sets.
Below, we will consider restricted diagrams on $(C^{op},S^{op})$, where $(C,S,\iota)$ is a pseudo-model category.  
  
\begin{df}\label{d20}  
\begin{enumerate}  
\item  
Let $(C,S)$ be a category with a distinguished subset of morphisms. The model category   
$SSet_{\mathbb{U}}^{C^{op},S^{op}}$, of restricted diagrams of simplicial sets   
on $(C^{op},S^{op})$ will be denoted by $(C,S)^{\wedge}$ and called the \emph{model   
category of pre-stacks on} $(C,S)$ (note that if $(C,S,\iota)$ is a pseudo-model category, $(C,S)^{\wedge}$   
does not depend on $\iota$).  
\item  
Let $(C,S,\iota)$ be a pseudo-model category and   
let $C^{c}$ (resp. $C^{f}$, resp. $C^{cf}$) be the full subcategory of $C$  
consisting of cofibrant (resp. fibrant, resp cofibrant and fibrant) objects, and $S^{c}:=C^{c}\cap S$  
(resp. $S^{f}:=C^{f}\cap S$, resp. $S^{cf}:=C^{cf}\cap S$).   
We will denote by $((C,S)^{c})^{\wedge}$ (resp.  $((C,S)^{f})^{\wedge}$, resp. $((C,S)^{cf})^{\wedge}$)  
the model category of restricted diagrams of $\mathbb{U}$-simplicial sets on $(C^{c},S^{c})^{op}$  
(resp. on $(C^{f},S^{f})^{op}$, resp. on $(C^{cf},S^{cf})^{op}$).  
  
\end{enumerate}  
\end{df}  
  
Objects of $(C,S)^{\wedge}$ are simply functors $F:C^{op} \longrightarrow SSet_{\mathbb{U}}$ and, as observed in Subsection \ref{resdia}, $F$ is fibrant in  $(C,S)^{\wedge}$ if and only if it is objectwise fibrant and preserves equivalences.

The category $(C,S)^{\wedge}$ is naturally tensored and co-tensored over $SSet_{\mathbb{U}}$,   
with external products and exponential objects defined objectwise. This makes $(C,S)^{\wedge}$ into a   
\textit{simplicial closed model category}. This model category is furthermore left proper,  $\mathbb{U}$-cellular and  
$\mathbb{U}$-combinatorial (see \cite[Ch. 14]{hi}, \cite{du2} and Appendix A).   
The derived simplicial $Hom$'s of the model category $(C,S)^{\wedge}$ will be denoted by  
$$\mathbb{R}_{w}\underline{Hom}(-,-) :   
\mathrm{Ho}((C,S)^{\wedge})^{op}\times \mathrm{Ho}((C,S)^{\wedge}) \longrightarrow  
\mathrm{Ho}((C,S)^{\wedge}).$$  
The derived simplicial $Hom$'s of the model categories $((C,S)^{c})^{\wedge}$, $((C,S)^{f})^{\wedge}$ and  
$((C,S)^{cf})^{\wedge}$, will be denoted similarly by   
$$\mathbb{R}_{w,c}\underline{Hom}(-,-) \qquad \mathbb{R}_{w,f}\underline{Hom}(-,-)  \qquad  
\mathbb{R}_{w,cf}\underline{Hom}(-,-).$$  
  
For an object $x \in C$, the evaluation functor $j_{x}^{*} : (C,S)^{\wedge} \longrightarrow SSet_{\mathbb{U}}$  
is a right Quillen functor. Its left adjoint is denoted by $(j_{x})_{!} :   
SSet_{\mathbb{U}} \longrightarrow (C,S)^{\wedge}$.   
We note that there is a canonical isomorphism $h_{x}\simeq(j_{x})_{!}(*)$ in $(C,S)^{\wedge}$,   
where $h_{x} : C^{op} \longrightarrow  
SSet_{\mathbb{U}}$ sends an object $y \in C$ to the constant simplicial set $Hom(y,x)$. More generally, for any  
$A \in SSet_{\mathbb{U}}$, one has $(j_{x})_{!}(A)\simeq A\otimes h_{x}$.   
  
As $(C,S)^{\wedge}$ is a left Bousfield localization of $SPr(C)$, the identity functor  
$\mathrm{Id} : SPr(C) \longrightarrow (C,S)^{\wedge}$ is left Quillen. In particular,   
homotopy colimits of diagrams in $(C,S)^{\wedge}$ can be computed in the objectwise model   
category $SPr(C)$. On the contrary, homotopy limits in $(C,S)^{\wedge}$ are not  
computed in the objectwise model structure; moreover, the identity functor $\mathrm{Id} : (C,S)^{\wedge} \longrightarrow  
SPr(C)$ does not preserve homotopy fibered products in general.   
  
As explained in Subsection \ref{resdia} (before Corollary \ref{c0}), if $(C,S)$ and $(C',S')$ are categories with distinguished subsets of morphisms   
(e.g., pseudo-model categories) and $f : C \rightarrow C'$ is a functor sending $S$ into $S'$, then we have a direct and inverse image Quillen adjunction 
$$f_{!} : (C,S)^{\wedge} \longrightarrow (C,S')^{\wedge} \qquad (C,S)^{\wedge} \longleftarrow (C',S')^{\wedge} : f^{*}.$$

In particular, if $(C,S,\iota)$ is a pseudo-model category, we may consider the inclusions   
$$(C^{c},S^{c}) \subset (C,S) \qquad (C^{f},S^{f}) \subset (C,S)  
\qquad (C^{cf},S^{cf}) \subset (C,S).$$  
As a consequence of Theorem \ref{t2} (or by a direct check), we get   
  
\begin{prop}\label{junction}  
Let $(C,S,\iota)$ be a pseudo-model category. The natural inclusions  
$$i_{c} : (C,S)^{c} \hookrightarrow (C,S) \qquad i_{f} : (C,S)^{f} \hookrightarrow (C,S) \qquad  
i_{cf} : (C,S)^{cf} \hookrightarrow (C,S),$$  
induce right Quillen equivalences   
$$i_{c}^{*} : (C,S)^{\wedge} \simeq ((C,S)^{c})^{\wedge} \qquad  
i_{f}^{*} : (C,S)^{\wedge} \simeq ((C,S)^{f})^{\wedge} \qquad i_{cf}^{*} : (C,S)^{\wedge} \simeq ((C,S)^{cf})^{\wedge}.$$  
These equivalences are furthermore compatible with derived simplicial $Hom$,   
in the sense that there exist natural isomorphisms  
$$\mathbb{R}_{w,c}\underline{Hom}(\mathbb{R}(i_{c})^{*}(-),\mathbb{R}(i_{c})^{*}(-))\simeq  
\mathbb{R}_{w}\underline{Hom}(-,-) $$  
$$\mathbb{R}_{w,f}\underline{Hom}(\mathbb{R}(i_{f})^{*}(-),\mathbb{R}(i_{f})^{*}(-))\simeq  
\mathbb{R}_{w}\underline{Hom}(-,-)$$  
$$\mathbb{R}_{w,cf}\underline{Hom}(\mathbb{R}(i_{cf})^{*}(-),\mathbb{R}(i_{cf})^{*}(-))\simeq  
\mathbb{R}_{w}\underline{Hom}(-,-).$$  
\end{prop}  
  
\end{subsection}  
  
\begin{subsection}{The Yoneda embedding of a pseudo-model category}  
  
Let us fix a pseudo-model category $(C,S,\iota:C\rightarrow M)$.  
Throughout this subsection we will also fix a \textit{cofibrant resolution  
functor} $(\Gamma : M \longrightarrow M^{\Delta},i)$ in the model category $M$   
(see \cite[$17.1.3, (1)$]{hi}). This means that for any object $x \in M$, $\Gamma(x)$ is a co-simplicial object  
in $M$, which is cofibrant for the Reedy model structure on $M^{\Delta}$, together  
with a natural equivalence $i(x) : \Gamma(x) \longrightarrow c^{*}(x)$,   
$c^{*}(x)$ being the constant co-simplicial object  
in $M$ at $x$. Let us remark that when the model category $M$ is simplicial, one can   
use the standard cofibrant resolution functor   
$\Gamma(x):=\Delta^{*}\otimes Q(x)$, where $Q$ is a cofibrant replacement  
functor in $M$.

We define the functor $\underline{h} : C \longrightarrow SPr(C)$, by sending each $x \in C$  
to the simplicial presheaf  
$$\begin{array}{cccc}  
\underline{h}_{x} : & M^{op} & \longrightarrow & SSet_{\mathbb{U}} \\  
& y & \mapsto & Hom_{M}(\Gamma(y),x),  
\end{array}$$  
where, 
to be more explicit, the presheaf  
of $n$-simplices of $\underline{h}_{x}$ is given by the formula  
$$(\underline{h}_{x})_{n}(-):=Hom_{M}(\Gamma(-)_{n},x).$$  
Note that for any $y \in M$, $\Gamma(y)_{n} \rightarrow y$ is   
an equivalence in $M$, therefore $y \in C$ implies that $\Gamma(y)   
\in C^{\Delta}$ (since $C$ is a pseudo-model category).  
  
We warn the reader that the two functors $h$ and $\underline{h}$ from $C$ to $(C,S)^{\wedge}$   
are different and should not be confused.  
For any $x \in C$, $h_{x}$ is a presheaf of discrete   
simplicial sets (i.e. a presheaf of sets) whereas $\underline{h}_{x}$  
is an actual simplicial presheaf. The natural equivalence $i(-) : \Gamma(-) \longrightarrow c^{*}(-)$ induces  
a morphism in $(C,S)^{\wedge}$  
$$h_{x}=Hom(c^{*}(-),x) \longrightarrow Hom(\Gamma(-),x)=\underline{h}_{x},$$   
which is functorial in $x \in M$.  
  
If, for a moment we denote by $\underline{h}^{C} : C \longrightarrow (C,S)^{\wedge}$ and by   
$\underline{h}^{M}:M \longrightarrow (M,W)^{\wedge}$ the functor defined for the pseudo-model categories   
$(C,S,\iota)$ and $(M,W,\mathrm{Id})$, respectively, we have a commutative diagram  
  
\begin{equation} \label{red}  
\xymatrix {  
C \ar[d]_{\iota} \ar[r]^{\underline{h}^{C}} & (C,S)^{\wedge}\\  
M \ar[r]^{\underline{h}^{M}} & M^{\wedge} \ar[u]_{\iota^{*}}  
}  
\end{equation}  
  
\noindent where $\iota^{*}$ is the restriction, right Quillen functor.\\  
  
\begin{lem}\label{l10}  
Both functors $\underline{h} : C \longrightarrow SPr(C)$   
and $\underline{h} : C \longrightarrow (C,S)^{\wedge}$ preserves fibrant objects and equivalences between them.  
\end{lem}  
  
\textit{Proof:} The statement for $\underline{h} : M \longrightarrow SPr(M)$ follows   
from the standard properties of mapping spaces, see \cite[\S 5.4]{ho} or \cite[Prop. 18.1.3,   
Thm. 18.8.7]{hi}. The statement for $\underline{h} : M \longrightarrow M^{\wedge}$ follows from the previous one and from \cite[Thm. 18.8.7 (2)]{hi}. Finally, the statements for   
$\underline{h}^{C} : C \longrightarrow SPr(C)$ and $\underline{h} : C \longrightarrow   
(C,S)^{\wedge}$ follow from the previous ones for $M$ and from the commutativity   
of diagram (\ref{red}), since $\iota^{*}$ is right Quillen. \hfill \textbf{$\Box$} \\  
  
Lemma \ref{l10} enables us to define a right derived functor of $\underline{h}$ as  
  
$$\begin{array}{cccc}  
\mathbb{R}\underline{h} : & S^{-1}C & \longrightarrow & \mathrm{Ho}((C,S)^{\wedge}) \\  
& x & \mapsto (\underline{h}\circ R \circ \iota)(x).  
\end{array}$$  
where $R$ denotes a fibrant replacement functor in   
$M$ and we implicitly used the fact that $R\iota(x)$ is still in $C$ for $x \in C$. Also note that, by definition of $(C,S)^{\wedge}$, the functor $h:C\longrightarrow (C,S)^{\wedge}$ preserves equivalences, hence induces a functor $\mathrm{Ho}(h):S^{-1}C \longrightarrow \mathrm{Ho}((C,S)^{\wedge})$.   
  
The reader should notice that if  
$(\Gamma',i')$ is another   
cofibrant resolution functor in $M$, then  
the two derived functor $\mathbb{R}\underline{h}$ and $\mathbb{R}\underline{h'}$ obtained using respectively $\Gamma$  
and $\Gamma'$, are naturally isomorphic. Therefore, our construction does not   
depend on the choice of $\Gamma$   
once one moves to the homotopy category.  
  
\begin{lem} \label{strange} The functors $\mathrm{Ho}(h)$ and $\mathbb{R}\underline{h}$ from $S^{-1}C$ to $\mathrm{Ho}((C,S)^{\wedge})$ are canonically isomorphic. More precisely, if $R$ be a fibrant replacement functor in $M$, then the natural equivalence $i(-): \Gamma(-) \longrightarrow c^{*}(-)$   
induces, for any $x\in C$, an equivalence in $(C,S)^{\wedge}$ (hence a fibrant replacement, by Lemma \ref{l10})  
$$h_{x}=Hom(-,x) \longrightarrow Hom(\Gamma(-),R(x))=\underline{h}_{R(x)}.$$  
\end{lem}  
\textit{Proof:}    
First we show that if $x$ is a fibrant and cofibrant object in $C$,   
then the natural morphism $h_{x} \longrightarrow \underline{h}_{x}$ is an equivalence in  
$((C,S)^{c})^{\wedge}$.  
To see this, let   
$x \longrightarrow x_{*}$ be a simplicial resolution of $x$ in $M$, hence in $C$ (see \cite[$17.1.2$]{hi}).  
We consider the following two simplicial presheaves  
$$\begin{array}{cccc}  
h_{x_{*}} : & (C^{c})^{op} & \longrightarrow & SSet_{\mathbb{U}} \\  
& y & \mapsto & Hom(y,x_{*}),  
\end{array}$$  
$$\begin{array}{cccc}  
\underline{h}_{x_{*}} : & (C^{c})^{op} & \longrightarrow & SSet_{\mathbb{U}} \\  
& y & \mapsto & \mathrm{diag}(Hom(\Gamma(y),x_{*})).  
\end{array}$$  
The augmentation $\Gamma(-) \longrightarrow c(-)$ and co-augmentation  
$x \longrightarrow x_{*}$ induce a commutative diagram in $((C,S)^{cf})^{\wedge}$  
$$\xymatrix{  
h_{x} \ar[r]^-{a} \ar[d]_-{b} & \underline{h}_{x} \ar[d]^-{d} \\  
h_{x_{*}}  \ar[r]^-{c} & \underline{h}_{x_{*}}. }$$  
By the properties of mapping spaces (see \cite[\S $5.4$]{ho}), both  morphisms $c$ and $d$ are equivalences  
in $SPr(C^{c})$. Furthermore, the morphism $h_{x} \longrightarrow h_{x_{*}}$ is  
isomorphic in $\mathrm{Ho}(SPr(C^{c}))$ to the induced morphism  
$h_{x} \longrightarrow \mathrm{hocolim}_{[n] \in \Delta}h_{x_{n}}$.   
As each morphism $h_{x} \longrightarrow h_{x_{n}}$ is  
an equivalence in $((C,S)^{c})^{\wedge}$, this implies that $d$ is an  
equivalence in $((C,S)^{c})^{\wedge}$. We deduce from this that  
also the natural morphism $h_{x} \longrightarrow \underline{h}_{x}$ is an equivalence in $((C,S)^{c})^{\wedge}$.  
Let us show how this implies that for any $x \in C$,  
the natural morphism $h_{x} \longrightarrow \underline{h}_{Rx}$ is an equivalence   
in $(C,S)^{\wedge}$.  
  
Since for any equivalence $z \rightarrow z'$ in $C$, the induced map   
$h_{z} \rightarrow h_{z'}$ is an equivalence in $(C,S)^{\wedge}$ (see Remark \ref{earth}), it is enough to show that, for any $x\in C$, the canonical map   
$h_{Rx} \longrightarrow \underline{h}_{Rx}$ is an equivalence.   
By the Yoneda lemma for $\mathrm{Ho}((C,S)^{\wedge})$, it is enough to show that the induced map   
$Hom_{\mathrm{Ho}((C,S)^{\wedge})}(\underline{h}_{Rx},F)\rightarrow   
Hom_{\mathrm{Ho}((C,S)^{\wedge})}(h_{Rx},F)$ is a bijection for any $F\in \mathrm{Ho}((C,S)^{\wedge})$. Now,  
$$Hom_{\mathrm{Ho}((C,S)^{\wedge})}(G,F)\simeq \pi_{0}(\mathbb{R}_{w}\underline{Hom}(G,F))$$  
for any $G$ and $F$ in $(C,S)^{\wedge}$, hence it is enough   
to show that we have an induced equivalence of simplicial sets   
$$\mathbb{R}_{w}\underline{Hom}(\underline{h}_{Rx},F))\simeq   
\mathbb{R}_{w}\underline{Hom}(h_{x},F).$$  
By the properties of mapping spaces (see \cite[\S $5.4$]{ho}), if $Q$ denotes a   
cofibrant replacement functor in $M$, the map $\underline{h}_{Rx} \longrightarrow \underline{h}_{QRx}$   
is an equivalence in $(C,S)^{\wedge}$; therefore, if we denote by   
$(-)_{c}$ th restriction to $C^{c}$, we have an equivalence of simplicial sets  
$$  
\mathbb{R}_{w}\underline{Hom}(\underline{h}_{Rx},F))\simeq   
\mathbb{R}_{w,c}\underline{Hom}((\underline{h}_{QRx})_{c},F_{c}).  
$$  
Since $QR(x)$ is fibrant and cofibrant, we have already proved that   
$$  
\mathbb{R}_{w,c}\underline{Hom}((\underline{h}_{QRx})_{c},F_{c})\longrightarrow   
\mathbb{R}_{w,c}\underline{Hom}((h_{x})_{c},F_{c})  
$$  
is an equivalence of simplicial sets and we conclude since 
$  
\mathbb{R}_{w,c}\underline{Hom}((h_{x})_{c},F_{c}))\simeq   
\mathbb{R}_{w}\underline{Hom}(h_{x},F)  
$ by Proposition \ref{junction}. 
\hfill \textbf{$\Box$} \\  
  
The main result of this subsection is the following theorem.  
  
\begin{thm}\label{t10}  
If $(C,S,\iota:C \rightarrow M)$ is a pseudo-model category, the functor   
$\mathbb{R}\underline{h} : S^{-1}C \longrightarrow \mathrm{Ho}((C,S)^{\wedge})$  
is fully faithful.  
\end{thm}  
  
\textit{Proof:}  
We will identify $C$ as a full subcategory of $M$ and $S^{-1}C$ as a   
full subcategory of $\mathrm{Ho}(M)$ using $\iota$.    
For any $x$ and $y$ in $S^{-1}C$, letting $R$ be a fibrant replacement functor in $M$, one has   
$$Hom_{S^{-1}C}(x,y)\simeq \pi_{0}(Hom_{M}(\Gamma(x),R(y))  
$$  
since $\mathrm{Ho}(\iota)$ is fully faithful and $Hom(\Gamma(-),R(-))$ is a homotopy mapping   
complex in $M$ (see \cite[$5.4$]{ho}). As $(C,S,\iota)$ is a pseudo-model  category, we have   
$Hom_{M}(\Gamma(x)),R(y))=Hom_{C}(\Gamma(x),R(y))$. But, by definition of $\underline{h}$   
and the enriched Yoneda lemma in $(C,S)^{\wedge}$, we have isomorphisms of simplicial sets  
$$Hom_{C}(\Gamma(x),R(y)) \simeq \underline{h}_{R(y)}(x) \simeq  
\underline{Hom}_{(C,S)^{\wedge}}(h_{x},\underline{h}_{R(y)}).$$  
Now, $h_{x}$ is cofibrant in $(C,S)^{\wedge}$ and, by Lemma \ref{l10},   
$\underline{h}_{R(y)}$ is fibrant in $(C,S)^{\wedge}$, so that    
$$  
\pi_{0}(\underline{Hom}_{(C,S)^{\wedge}}(h_{x},\underline{h}_{R(y)})) \simeq   
Hom_{\mathrm{Ho}((C,S)^{\wedge})}(h_{x},\underline{h}_{R(y)})  
$$  
since $(C,S)^{\wedge}$ is a simplicial model category. Finally, by Lemma \ref{strange} we have  
$$Hom_{\mathrm{Ho}((C,S)^{\wedge})}(h_{x},\underline{h}_{R(y)})\simeq   
Hom_{\mathrm{Ho}((C,S)^{\wedge})}(\mathbb{R}\underline{h}_{x},\mathbb{R}\underline{h}_{y})$$  
showing that $\mathbb{R}\underline{h}$ is fully faithful.   
\hfill \textbf{$\Box$}  

\begin{cor}
For any $x\in C$ and any $F\in SPr(C)$, there is an isomorphism in $\mathrm{Ho}(SSet)$ $$\mathbb{R}_{w}\underline{Hom}_{(C,S)^{\wedge}}(\underline{h}_{x},F)\simeq F(x).$$
\end{cor}  

\begin{df}\label{d21'}  
For any pseudo-model category $(C,S,\iota)$ which is $\mathbb{U}$-small, the fully faithful emdedding  
$$\mathbb{R}\underline{h} : \mathrm{Ho}(C,S) \longrightarrow \mathrm{Ho}((C,S)^{\wedge})$$  
is called the \emph{Yoneda embedding} of $(C,S,\iota)$.  
\end{df}

\begin{rmk}\emph{\begin{enumerate}
	\item According to Definition \ref{d21'}, the Yoneda embedding of a pseudo-model category a priori depends on the embedding $\iota : C \hookrightarrow M$. However, it will be shown in \ref{c16} that it \textit{only} depends on the pair $(C,S)$.
	\item The Yoneda embedding for (pseudo-)model categories is one of the key ingredients used in \cite{kthe} to prove that, for a large class of Waldhausen categories, the $K$-theory only depends on the Dwyer-Kan simplicial localization (though it is known to depend on strictly more than the usual localization).
\end{enumerate}
}  
\end{rmk}  
  
\end{subsection}  
  
\begin{subsection}{Model pre-topologies and pseudo-model sites}  
  
\begin{df}\label{d21}  
A \emph{model pre-topology} $\tau$ on a $\mathbb{U}$-small pseudo-model  
category $(C,S,\iota)$, is the datum for any object $x \in C$, of  
a set $Cov_{\tau}(x)$ of subsets  
of objects in $\mathrm{Ho}(C,S)/x$, called $\tau$\emph{-covering families} of  
$x$, satisfying the following three conditions.  
  
\begin{enumerate}  
  
\item \emph{(Stability)} For all $x \in C$ and any isomorphism $y \rightarrow x$ in $\mathrm{Ho}(C,S)$,  
the one-element set $\{y \rightarrow x\}$ is in $Cov_{\tau}(x)$.  
  
\item \emph{(Composition)} If $\{u_{i} \rightarrow x\}_{i \in I} \in Cov_{\tau}(x)$, and for any $i \in I$,  
$\{v_{ij} \rightarrow u_{i}\}_{j \in J_{i}} \in Cov-{\tau}(u_{i})$,   
the family $\{v_{ij} \rightarrow x\}_{i \in I, j\in J_{i}}$  
is in $Cov_{\tau}(x)$.  
  
\item \emph{(Homotopy base change)} Assume the two previous conditions hold. For any   
$\{u_{i} \rightarrow x\}_{i \in I} \in Cov_{\tau}(x)$, and any morphism in $\mathrm{Ho}(C,S)$,  
$y \rightarrow x$, the family $\{u_{i}\times^{h}_{x}y \rightarrow y\}_{i \in I}$ is in $Cov_{\tau}(y)$.  
  
\end{enumerate}  
A $\mathbb{U}$-small pseudo-model category $(C,S,\iota)$  
together with a model pre-topology $\tau$ will be called a $\mathbb{U}$-small \emph{pseudo-model site}.  
\end{df}  
  
\begin{rmk}  
\emph{  
\begin{enumerate}  
\item Note that in the third condition (\textit{Homotopy base-change}) we used the   
homotopy fibered product of diagrams $\xymatrix{x\ar[r] & z & \ar[l] y}$  
in $\mathrm{Ho}(M)$. By this we mean the homotopy fibered product of a lift (up to equivalence)  
of this diagram to $M$. This is a well defined object in $\mathrm{Ho}(M)$ but only up to a   
\textit{non-canonical} isomorphism in $\mathrm{Ho}(M)$  
(in particular it is not functorially defined). However, condition $(3)$ of the previous   
definition still makes sense because we assumed the two previous conditions $(1)$ and $(2)$  
hold.  
\item When the pseudo-model structure on $(C,S)$ is trivial as in Remark \ref{trivia} 2, a   
model pre-topology on $(C,S)$ is the same thing as a Grothendieck   
pre-topology on the category $C$ as defined in \cite[Exp. II]{sga4}. Indeed,   
in this case  we have a canonical identification $\mathrm{Ho}(C,S)=C$ under which homotopy   
fibered products correspond to fibered products.  
\end{enumerate}  
}  
\end{rmk}  
  
\medskip  
  
Let $(C,S,\iota;\tau)$ be a $\mathbb{U}$-small pseudo-model site and $\mathrm{Ho}(C,S)=S^{-1}C$ the homotopy category   
of $(C,S)$. A sieve $R$ in $\mathrm{Ho}(C,S)$ over an object $x \in   
\mathrm{Ho}(C,S)$ will be called a $\tau$\textit{-covering sieve}   
if it contains a $\tau$-covering family.   
  
\begin{lem}\label{l11}  
For any $\mathbb{U}$-small pseudo-model site $(C,S,\iota;\tau)$, the $\tau$-covering   
sieves form a Grothendieck topology on $\mathrm{Ho}(C,S)$.  
\end{lem}  
  
\textit{Proof:} The stability and composition axioms of  
Definition \ref{d21} clearly imply conditions (i') and (iii') of \cite[Ch. III, \S 2,  Def. 2]{mm}.   
It is also clear that if $u : y \rightarrow x$ is any  
morphism in $\mathrm{Ho}(C,S)$, and if $R$ is a sieve on $x$ which contains a $\tau$-covering family  
$\{u_{i} \rightarrow x\}_{i \in I}$, then the pull-back sieve $u^{*}(R)$ contains the  
family $\{u_{i}\times^{h}_{x}y \rightarrow y\}_{i \in I}$. Therefore, the homotopy base change axiom of  
Definition \ref{d21} implies condition (ii') of \cite[Ch. III, \S 2, Def. 2]{mm}. \hfill \textbf{$\Box$} \\  
  
The previous lemma shows that any ($\mathbb{U}$-small) pseudo-model site  
$(C,S,\iota\;\tau)$ gives rise to a ($\mathbb{U}$-small) $S$-site $(L(C,S),\tau)$, where $\mathrm{L}(C,S)$  
is the Dwyer-Kan localization of $C$ with respect to $S$ and   
$\tau$ is the Grothendieck topology on $\mathrm{Ho}(\mathrm{L}(C,S))=\mathrm{Ho}(C,S)$   
defined by $\tau$-covering sieves. We will say that   
the $S$-topology $\tau$ on $\mathrm{L}(C,S)$ is \textit{generated} by the pre-topology $\tau$ on $(C,S)$.   
  
\textit{Conversely}, a topology on $\mathrm{Ho}(C,S)$ induces a model pre-topology   
on the pseudo-model category $(C,S,\iota)$ in  
the following way. A subset of objects $\{u_{i}\rightarrow x\}_{i \in I}$   
in $\mathrm{Ho}(C,S)/x$ is defined   
to be a $\tau$-covering family if the sieve it generates is a covering   
sieve (for the given topology on $\mathrm{Ho}(C,S)$).   
  
\begin{lem}\label{l12}  
Let $(C,S,\iota)$ be a $\mathbb{U}$-small pseudo-model category and   
let $\tau$ be a Grothendieck topology on $\mathrm{Ho}(C,S)$.   
Then, the $\tau$-covering   
families in $\mathrm{Ho}(C,S)$ defined above form a model pre-topology on $(C,S,\iota)$,   
called the \emph{induced} model pre-topology.   
\end{lem}   
  
\textit{Proof:} Conditions $(1)$ and $(2)$ of Definition \ref{d21} are   
clearly satisfied and it only remains to check   
condition $(3)$. For this, let us recall that the homotopy fibered products   
have the following semi-universal property in $\mathrm{Ho}(C,S)$.  
For any commutative diagram in $\mathrm{Ho}(C,S)$  
$$\xymatrix{  
x \ar[r] \ar[d] & y \ar[d] \\  
z \ar[r] & t,}$$  
there exists a morphism $x \rightarrow z\times_{t}^{h}y$ compatible with   
the two projections to $z$ and $y$. Using this property  
one sees that for any subset of objects $\{u_{i} \rightarrow x\}_{i \in I}$   
in $\mathrm{Ho}(C,S)/x$, and any morphism $u : y \rightarrow x$,   
the sieve over $y$ generated by the family $\{u_{i}\times^{h}_{x}y \rightarrow y\}_{i \in I}$   
coincides with the pull-back by  
$u$ of the sieve generated by $\{u_{i}\rightarrow x\}_{i \in I}$. Therefore,   
the base change axiom (ii') of   
\cite[Ch. III, \S 2, Def. 2]{mm} implies the homotopy base change property   
$(3)$ of Definition \ref{d21}. \hfill \textbf{$\Box$} \\  
  
Lemmas \ref{l11} and \ref{l12} show that model pre-topologies on a   
pseudo-model category $(C,S)$ are essentially the same  
as Grothendieck topologies on $\mathrm{Ho}(C,S)$, and therefore the   
same thing as $S$-topologies on the   
$S$-category $\mathrm{L}(C,S)$. As in the usual case (i.e. for the  
trivial model structure on $(C,S)$) the two above constructions are   
not exactly mutually inverse but we have the following   
  
\begin{prop}\label{upanddown}  
Let $(C,S,\iota)$ be a pseudo-model category. The rule assigning to a model pre-topology $\tau$ on $(C,S,\iota)$ the   
$S$-topology on $\mathrm{L}(C,S)$ generated by $\tau$ and the rule assigning to an $S$-topology on   
$\mathrm{L}(C,S)$ the induced model pre-topology on $(C,S,\iota)$, induce a bijection  
  
$$  
\xymatrix{  
\left\{ \txt{\emph{Saturated model}\\ \emph{pre-topologies on (C,S,$\iota$)} } \right\}   
\ar@{<->}[rr] &&  \left\{ \txt{$S$\emph{-topologies}\\\emph{on L(C,S)}}\right\}}  
$$  
  
\noindent where we call a model pretopology $\tau$ \emph{saturated} if any family of   
morphisms in $\mathrm{Ho}(C,S)/x$ that contains a $\tau$-covering family   
for $x$ is again a $\tau$-covering family for $x$.  
\end{prop}  
  
\textit{Proof:} Straightforward from Lemma \ref{l11} and \ref{l12}. \hfill \textbf{$\Box$} 
  
\begin{ex} \label{extop}   
\emph{  
\begin{enumerate}  
\item \textit{Topological spaces.} Let us take as $C=M$ the model category of  
$\mathbb{U}$-topological spaces, $Top$, with $S=W$ consisting of the usual weak equivalences. We define  
a model pre-topology $\tau$ in the following way.  
A family of morphism in $\mathrm{Ho}(Top)$, $\{X_{i} \rightarrow X\}_{i \in I}$,  
$I \in \mathbb{U}$, is defined to be in $Cov_{\tau}(X)$ if  
the induced map $\coprod_{i\in I}\pi_{0}(X_{i}) \longrightarrow \pi_{0}(X)$ is surjective.   
The reader will check easily that this defines a topology on $Top$ in the sense of Definition \ref{d21}.  
\item \textit{Strong model pre-topologies for $E_{\infty}$-algebras over $k$.} Let $k$ be a   
commutative ring (respectively, a field of characteristic zero) and let   
$C=M:=(E_{\infty}-\mathrm{Alg}_{k})^{op}$ (resp. $C=M:=(\mathrm{CDGA}_{\leq0 ;\;  k})^{op}$)   
be the opposite model category of $E_{\infty}$-algebras over the category of (unbounded)   
complexes of $k$-modules (resp., the opposite model category of commutative and unital   
differential graded $k$-algebras in negative degrees) which belong to  
$\mathbb{U}$; see for example \cite{hin} and \cite{bg} for a description of  
these model structures. Let $\tau$ be one of the usual topologies  
defined on $k$-schemes (e.g. Zariski, Nisnevich, \'{e}tale, ffpf  
or ffqc). Let us define the \textit{strong topology} $\tau_{\textrm{str}}$ on $M$ in the sense of   
Definition \ref{d21}, as follows. A family of morphisms in $\mathrm{Ho}(M^{op})$, $\{B  
\rightarrow A_{i}\}_{i \in I}$, $I \in \mathbb{U}$, is defined to  
be in $Cov_{\tau_{\textrm{str}}}(B)$ if it satisfies the two following  
conditions.  
\begin{itemize}  
\item The induced family of morphisms of affine $k$-schemes   
$\{Spec\, H^{0}(A_{i}) \rightarrow Spec H^{0}(B)\}_{i \in I}$  
is a $\tau$-covering.  
\item For any $i \in I$, one has $H^{*}(A_{i})\simeq H^{*}(B)\otimes_{H^{0}(B)}H^{0}(A_{i}).$  
\end{itemize}  
\noindent In the case of negatively graded commutative differential graded algebras over a field of   
characteristic zero, the strong \'etale topology $(\textrm{\'et})_{\textrm{str}}$ has been considered in   
\cite{be}. We will use these kind of model pre-topologies in \cite{partII} to give another approach to the theory of DG-schemes of \cite{ck1} and \cite{ck2} (or, more generally, to  the theory of $E_{\infty}$-schemes, when the base ring is not a   
field of characteristic zero) by viewing them as \textit{geometric stacks over the category of   
complexes of $k$-modules}.  
\item \textit{Semi-strong model pre-topologies for   
$E_{\infty}$-algebras over $k$.} With the same notations   
as in the previous example, we define the \textit{semi-strong topology}  
$\tau_{\textrm{sstr}}$ on $M$ by stipulating that a family of morphisms in $\mathrm{Ho}(M^{op})$,   
$\{B\rightarrow A_{i}\}_{i \in I}$, $I \in \mathbb{U}$, is in $Cov_{\tau_{\textrm{sstr}}}(B)$ if   
the induced family of morphisms of affine $k$-schemes $$\{Spec\, H^{*}(A_{i}) \rightarrow Spec H^{*}(B)\}_{i \in I}$$  
is a $\tau$-covering.  
\item \textit{The $Tor_{\geq0}$ model pre-topology for   
$E_{\infty}$-algebras over $k$.} Let $k$ be a commutative ring and $C=M:=(E_{\infty}-\mathrm{Alg}_{k})^{op}$ be the   
opposite model category of $E_{\infty}$-algebras over the category of (unbounded)   
complexes of $k$-modules which belong to $\mathbb{U}$. For any $E_{\infty}$-algebra   
$A$, we denote by $\mathrm{Mod}_{A}$ the model category of $A$-modules (see \cite{hin}   
or \cite{sp}). We define the \textit{positive $Tor$-dimension} pre-topology, $Tor_{\geq0}$,   
on $M$, as follows. A family of morphisms in   
$\mathrm{Ho}(M^{op})$, $\{f_{i}:B\rightarrow A_{i}\}_{i \in I}$, $I \in \mathbb{U}$, is defined to  
be in $Cov_{Tor_{\geq 0}}(B)$ if it satisfies the two following  
conditions.  
\begin{itemize}  
\item For any $i \in I$, the derived base-change functor $\mathbb{L}f_{i}^{*}=-\otimes_{B}^{\mathbb{L}}A_{i} :   
\mathrm{Ho}(\mathrm{Mod}_{B}) \longrightarrow \mathrm{Ho}(\mathrm{Mod}_{A_{i}})$ preserves the subcategories of   
\textit{positive modules} (i.e. of modules $P$ such that $H^{i}(P)=0$ for any $i\leq 0$).  
\item The family of derived base-change functors  
$$\{\mathbb{L}f_{i}^{*} : \mathrm{Ho}(\mathrm{Mod}_{B}) \longrightarrow \mathrm{Ho}(\mathrm{Mod}_{A_{i}})\}_{i \in I}$$  
is conservative (i.e. a morphism in $\mathrm{Ho}(\mathrm{Mod}_{B})$ is an isomorphism if and only if,  
for any $i \in I$, its image in $\mathrm{Ho}(\mathrm{Mod}_{A_{i}})$ is an isomorphism).  
\end{itemize}  
\noindent This positive $Tor$-dimension pre-topology is particularly relevant in interpreting \textit{higher   
tannakian duality} (\cite{to1}) as a part of algebraic geometry over the category of unbounded   
complexes of $k$-modules. We will come back on this in \cite{partII}.  
\end{enumerate}  
}  
\end{ex}

We fix a model pre-topology $\tau$ on a pseudo-model category   
$(C,S,\iota)$ and consider the pseudo-model site $(C,S,\iota;\tau)$.   
The induced Grothendieck topology on $\mathrm{Ho}(C,S)$   
described in the previous  paragraphs will still be denoted by $\tau$.   
  
Let $F \in (C,S)^{\wedge}$ be a pre-stack on the pseudo-model site $(C,S,\iota;\tau)$,   
and let $F \rightarrow RF$ be a fibrant replacement of $F$ in $(C,S)^{\wedge}$. We may    
consider the presheaf of connected components of $RF$, defined as  
$$\begin{array}{cccc}  
\pi_{0}^{pr}(RF) : C^{op} & \longrightarrow & Set \\  
 x & \mapsto & \pi_{0}(RF(x)).  
\end{array}$$  
Since any other fibrant model of $F$ in $(C,S)^{\wedge}$ is actually objectwise equivalent to $RF$, the presheaf   
$\pi_{0}^{pr}(RF)$ is  
well defined up to a unique isomorphism. This defines a functor  
$$\begin{array}{cccc}  
\pi_{0}^{eq} : (C,S)^{\wedge} & \longrightarrow & Pr(C) \\  
 F & \mapsto & \pi_{0}^{pr}(RF).  
\end{array}$$  
As $RF$ is fibrant, it sends equivalences in $C$ to equivalences of simplicial sets, hence the presheaf $\pi_{0}^{eq}(F)$ always sends equivalences in $C$ to  
isomorphisms, so it factors through $\mathrm{Ho}(C,S)^{op}$. Again, this defines a functor $$\begin{array}{cccc}  
\pi_{0}^{eq} : (C,S)^{\wedge} & \longrightarrow & Pr(\mathrm{Ho}(C,S)) \\  
 F & \mapsto & \pi_{0}^{eq}(F).  
\end{array}$$  
Finally, if $F \longrightarrow G$ is an equivalence in $(C,S)^{\wedge}$, the   
induced morphism $RF \longrightarrow RG$ is an objectwise equivalence, and therefore the induced morphism   
$\pi_{0}^{eq}(F) \longrightarrow \pi_{0}^{eq}(G)$ is an isomorphism of presheaves of sets.   
In other words, the functor $\pi_{0}^{eq}$ factors through the homotopy category $\mathrm{Ho}((C,S)^{\wedge})$ as   
$$\pi_{0}^{eq} : \mathrm{Ho}((C,S)^{\wedge}) \longrightarrow Pr(\mathrm{Ho}(C,S)).$$  
  
\begin{df}\label{d22}  
Let $(C,S,\iota;\tau)$ be a pseudo-model site in $\mathbb{U}$.  
\begin{enumerate}  
\item For any object $F \in (C,S)^{\wedge}$, the sheaf associated to the presheaf $\pi_{0}^{eq}(F)$  
is denoted by $\pi_{0}^{\tau}(F)$ (or $\pi_{0}(F)$ if the topology $\tau$ is   
clear from the context). It is a usual sheaf on the site   
$(\mathrm{Ho}(C,S),\tau)$, and is called the \emph{sheaf of connected components} of $F$;  
\item   
A morphism $f : F \longrightarrow G$ in $\mathrm{Ho}((C,S)^{\wedge})$ is   
called a $\tau$\emph{-covering} (or just a \emph{covering}   
if the topology $\tau$ is clear from the context) if the induced morphism of presheaves  
$\pi_{0}^{eq}(F) \longrightarrow \pi_{0}^{eq}(G)$  
induces an epimorphism of sheaves on $\mathrm{Ho}(C,S)$ for the topology $\tau$;  
\item A morphism $F \longrightarrow G$ in $(C,S)^{\wedge}$   
is called a $\tau$\emph{-covering} (or just a \emph{covering}  
if the topology $\tau$ is clear) if the induced morphism in   
$\mathrm{Ho}((C,S)^{\wedge})$ is a $\tau$\emph{-covering} according to the previous definition.   
\end{enumerate}  
\end{df}  
  
Coverings in the model category $(C,S)^{\wedge}$ behave exactly as   
coverings in the model category of pre-stacks over an $S$-site (see Subsection \ref{stop}).  
It is easy to check (compare to Proposition \ref{hls}) that   
a morphism $F \longrightarrow G$ between fibrant objects in $(C,S)^{\wedge}$ is a $\tau$-covering  
iff for any object $x \in C$ and any morphism $h_{x} \longrightarrow G$ in $(C,S)^{\wedge}$, there exists  
a covering family $\{u_{i} \rightarrow x\}_{i \in I}$ in $C$ (meaning that its   
image in $\mathrm{Ho}(C,S)$ is a $\tau$-covering family),   
and for each $i \in I$, a commutative diagram in $\mathrm{Ho}((C,S)^{\wedge})$  
$$\xymatrix{  
F \ar[r] & G \\  
h_{u_{i}} \ar[r] \ar[u] & h_{x}. \ar[u]}$$  
Moreover, we have the following analog of proposition \ref{p1}.   
  
\begin{prop}\label{p5}  
Let $(C,S,\iota;\tau)$ be a pseudo-model site.  
\begin{enumerate}  
  
\item A morphism in $SPr(T)$ which is a composition of coverings is a covering.  
  
\item Let   
$$\xymatrix{  
F' \ar[r]^{f'} \ar[d] & G' \ar[d] \\  
F \ar[r]_{f} & G }$$  
be a homotopy cartesian diagram in $(C,S)^{\wedge}$. If $f$ is a covering then so is $f'$.   
\item Let $\xymatrix{F \ar[r]^-{u} & G \ar[r]^-{v} & H}$ be two morphisms in $(C,S)^{\wedge}$. If   
the morphism $v\circ u$ is a covering then so is $v$.  
  
\item Let   
$$\xymatrix{  
F' \ar[r]^{f'} \ar[d] & G' \ar[d]^{p} \\ F \ar[r]_{f} & G }$$  
be a homotopy cartesian diagram in $(C,S)^{\wedge}$. If $p$ and $f'$ are coverings then so is $f$.   
  
\end{enumerate}  
\end{prop}  
  
\textit{Proof:} Easy exercise left to the reader. \hfill \textbf{$\Box$} 
  
\end{subsection}  
  
\begin{subsection}{Simplicial objects and hypercovers}\label{Mhhc}  
  
In this subsection we fix a pseudo-model site $(C,S,\iota; \tau)$ in $\mathbb{U}$ and keep the notations of Subsection \ref{sobj}, with $SPr(T)$ replaced by $(C,S)^{\wedge}$; more precisely we take $T=L(C^{op},S^{op})$ (with the induced $S$-topology, see Prop. \ref{upanddown}) and use Theorem \ref{t2} with $M=SSet$ to have definitions and results of Subsection \ref{sobj} available for $(C,S)^{\wedge}=SSet_{\mathbb{U}}^{C^{op},S^{op}}$. \\

We introduce a nice class of hypercovers that will be used in the proof  
of the existence of the local model structure; this class will replace our distinguised set of hypercovers $H$ used in the proof of Theorem \ref{t3}.   
  
\begin{df}\label{dpseudorep}  
\begin{enumerate}  
  
\item An object $F\in (C,S)^{\wedge}$ is called \emph{pseudo-representable}  
if it is a $\mathbb{U}$-small disjoint union of representable presheaves  
$$F\simeq \coprod_{u\in I}h_{u}.$$  
\item A morphism between pseudo-representable objects  
$$f : \coprod_{u\in I}h_{u} \longrightarrow \coprod_{v\in J}h_{v}$$  
is called a \emph{pseudo-fibration}  
if for all $u \in I$, the corresponding projection  
$$f \in \prod_{u\in I}\coprod_{v \in J}Hom(h_{u},h_{v}) \longrightarrow \coprod_{v \in J}Hom(h_{u},h_{v})  
\simeq \coprod_{v \in J}Hom_{C}(u,v)$$  
is represented by a fibration in $C$.  
  
\item Let   
$$f : \coprod_{u\in I}h_{u} \longrightarrow \coprod_{v\in J}h_{v}$$  
be a morphism between pseudo-representable objects, and for any $v \in J$ let   
$I_{v}$ be the sub-set of $I$ of components $h_{u}$ which are sent to $h_{v}$.   
The morphism is called a \emph{pseudo-covering} if for any $v \in J$,   
the family of morphisms   
$$\{h_{u} \rightarrow h_{v}\}_{u \in I_{v}}$$  
corresponds to a covering family in the pseudo-model site $(C,S)$.   
  
\item Let $x$ be a fibrant object in $C$. A \emph{pseudo-representable hypercover} of $x$ is an object  
$F_{*} \longrightarrow h_{x}$ in $s(C,S)^{\wedge}/h_{x}$ such that for any integer $n\geq 0$  
the induced morphism  
$$F_{n} \longrightarrow F_{*}^{\partial \Delta^{n}}\times_{h_{x}^{\partial \Delta^{n}}}h_{x}^{\Delta^{n}}$$  
is a pseudo-fibration and a pseudo-covering between pseudo-representable objects.  
\end{enumerate}  
\end{df}  
  
The first thing to check is that pseudo-representable hypercovers are   
hypercovers.  
  
\begin{lem}  
A pseudo-representable hypercover $F_{*} \longrightarrow h_{x}$ is a  
$\tau$-hypercover (see Definition \ref{d12}).  
\end{lem}  
  
\textit{Proof:} It is enough to check that the natural morphism  
$$F_{*}^{\partial \Delta^{n}}\times_{h_{x}^{\partial \Delta^{n}}}h_{x}^{\Delta^{n}} \longrightarrow  
F_{*}^{\mathbb{R}\partial \Delta^{n}}\times^{h}_{h_{x}^{\mathbb{R}\partial \Delta^{n}}}h_{x}^{\mathbb{R}\Delta^{n}}$$  
is an isomorphism in $\mathrm{Ho}((C,S)^{\wedge})$. But this   
follows from the fact that $h$ preserves finite limits (when they  
exists) and the fact that $(C,S)$ is a pseudo-model category. \hfill $\Box$ 
  
\end{subsection}  
  
\begin{subsection}{Local equivalences}  
  
This subsection is completely analogous (actually a bit   
easier, because the notion of comma site is completely harmless here) to Subsection \ref{loceq}.  
  
Let $(C,S,\iota; \tau)$ be a $\mathbb{U}$-small pseudo-model site, and $x$ be a fibrant object in $C$.  
The comma category $(C/x,S,\iota)$ is then endowed with its natural structure of a pseudo-model category.  
The underlying category is $C/x$, the category of objects over $x$. The equivalences $S$ in $C/x$ are simply the  
morphisms whose images in $C$ are equivalences. Finally, the embedding $\iota : C \longrightarrow M$ induces  
an embedding $\iota : C/x \longrightarrow M/\iota(x)$. The comma category  
$M/\iota(x)$ is endowed with its natural model category structure (see \cite[\S $1$]{ho}). It is easy to check that  
$(C/x,S,\iota)$ is a pseudo-model category in the sense of Definition \ref{pseudo}.  
  
We define a model pre-topology, still denoted by $\tau$, on the   
comma pseudo-model category $(C/x,S,\iota)$ by declaring that a family  
$\{u_{i} \rightarrow y\}_{i \in I}$ of objects in $\mathrm{Ho}((C/x,S))/y$ is a $\tau$-covering family   
if its image family under the natural functor $\mathrm{Ho}((C/x,S))/y \longrightarrow \mathrm{Ho}((C,S))/y$  
is a $\tau$-covering family for $y$. As the object $x$ is fibrant in $(C,S)$ the forgetful functor  
$(C/x,S) \longrightarrow (C,S)$ preserves homotopy fibered products,   
and therefore one checks immediately that this defines   
a model pre-topology $\tau$ on $(C/x,S,\iota)$.  
  
\begin{df}\label{d24}  
The pseudo-model site $(C/x,S,\iota;\tau)$ will be called the \emph{comma pseudo-model site}  
of $(C,S,\iota;\tau)$ over the (fibrant) object $x$.  
\end{df}  
  
\begin{rmk} \emph{Note that in the case where $(C,S,\iota)$ is a   
right proper pseudo-model category, the hypothesis that $x$ is fibrant is unnecessary.}   
\end{rmk}  
  
For any object $x \in C$, the evaluation functor  
$$\begin{array}{cccc}  
j_{x}^{*} : & (C,S)^{\wedge} & \longrightarrow & SSet_{\mathbb{U}} \\  
 & F & \mapsto & F(x)   
\end{array}$$  
has a left adjoint $(j_{x})_{!}$. The adjunction  
$$(j_{x})_{!} : SSet_{\mathbb{U}} \longrightarrow (C,S)^{\wedge} \qquad   
SSet_{\mathbb{U}} \longleftarrow (C,S)^{\wedge} : j_{x}^{*}$$  
is clearly a Quillen adjunction.  
  
Let $F \in (C,S)^{\wedge}$, $x$ a fibrant object  in $(C,S)$ and   
$s \in \pi^{eq}_{0}(F(x))$ be represented by a morphism    
$s : h_{x} \longrightarrow F$   
in $\mathrm{Ho}((C,S)^{\wedge})$.  
By pulling-back this morphism through the functor   
$$\mathbb{R}j_{x}^{*} : \mathrm{Ho}((C,S)^{\wedge}) \longrightarrow \mathrm{Ho}((C/x,S)^{\wedge})$$   
one gets a morphism in $\mathrm{Ho}((C/x,S)^{\wedge})$  
$$s : \mathbb{R}j_{x}^{*}(h_{x}) \longrightarrow \mathbb{R}j_{x}^{*}(F).$$  
By definition of the comma pseudo-model category $(C/x,S)$, it is immediate that   
$\mathbb{R}j_{x}^{*}(h_{x})$ has a natural global point  
$* \longrightarrow \mathbb{R}j_{x}^{*}(h_{x})$ in $\mathrm{Ho}((C/x,S)^{\wedge})$.   
Observe that the morphism $* \longrightarrow \mathbb{R}j_{x}^{*}(h_{x})$  
can also be seen as induced by adjunction from the identity of $h_{x}\simeq \mathbb{L}(j_{x})_{!}(*)$.  
We therefore obtain a global point  
$$s : * \longrightarrow \mathbb{R}j_{x}^{*}(h_{x})\longrightarrow \mathbb{R}j_{x}^{*}(F).$$   
  
\begin{df}\label{d25}  
\begin{enumerate}  
\item   
For an integer $n>0$, the sheaf $\pi_{n}(F,s)$ is defined to be   
$$\pi_{n}(F,s):=\pi_{0}(\mathbb{R}j_{x}^{*}(F)^{\mathbb{R}\Delta^{n}}  
\times_{\mathbb{R}j_{x}^{*}(F)^{\mathbb{R}\partial \Delta^{n}}}*).$$  
It is a usual sheaf on the site $(\mathrm{Ho}(C/x,S),\tau)$ called the $n$\emph{-th  
homotopy sheaf} of $F$ \emph{pointed at} $s$.   
  
\item A morphism $f : F \longrightarrow G$ in $(C,S)^{\wedge}$ is called a  
$\pi_{*}$\emph{-equivalence} (or equivalently a  
\emph{local equivalence}) if the following two conditions are satisfied:  
\begin{enumerate}  
\item The induced morphism   
$\pi_{0}(F) \longrightarrow \pi_{0}(G)$  
is an isomorphism of sheaves on $\mathrm{Ho}(C,S)$;  
\item For any fibrant object $x \in (C,S)$, any section $s \in \pi_{0}^{eq}(F(x))$  
and any integer $n>0$, the induced morphism  
$\pi_{n}(F,s) \longrightarrow \pi_{n}(G,f(s))$  
is an isomorphism of sheaves on $\mathrm{Ho}(C/x,S)$.   
\end{enumerate}  
\end{enumerate}  
\end{df}  
  
\medskip  
  
As observed in Subsection \ref{loceq}, an equivalence in the model category $(C,S)^{\wedge}$ is always a $\pi_{*}$-equivalence, for any model pre-topology $\tau$ on $(C,S)$.   
   
The $\pi_{*}$-equivalences in $(C,S)^{\wedge}$ behave the same way as the  
$\pi_{*}$-equivalences in $SPr(T)$ (see Subsection \ref{loceq}). We   
will therefore state the following basic facts without repeating their proofs.  
  
\begin{lem}\label{l13}  
A morphism $f : F \longrightarrow G$ in $(C,S)^{\wedge}$ is a $\pi_{*}$-equivalence if and only if  
for any $n\geq 0$, the induced morphism   
$$F^{\mathbb{R}\Delta^{n}} \longrightarrow F^{\mathbb{R}\partial   
\Delta^{n}}\times^{h}_{G^{\mathbb{R}\partial \Delta^{n}}}G^{\mathbb{R}\Delta^{n}}$$  
is a covering.  
In other words $f$ is a $\pi_{*}$-equivalence if and only if it is a $\tau$-hypercover   
when considered as a morphism between constant simplicial objects in $(C,S)^{\wedge}$.  
\end{lem}  
  
\begin{cor}\label{c12}  
Let $f : F \longrightarrow G$ be a morphism in $(C,S)^{\wedge}$ and   
$G' \longrightarrow G$ be a covering. Then, if the induced morphism  
$$f' : F\times_{G}^{h}G' \longrightarrow G'$$  
is a $\pi_{*}$-equivalence, so is $f$.  
\end{cor}  
  
\medskip  
  
Let $f : F \longrightarrow G$ be a morphism in $(C,S)^{\wedge}$.  
For any fibrant object $x \in (C,S)$ and any morphism $s : h_{x}   
\longrightarrow G$ in $\mathrm{Ho}((C,S)^{\wedge})$, let   
us define $F_{s} \in \mathrm{Ho}(((C,S)/x)^{\wedge})$ via the following homotopy cartesian square  
$$\xymatrix{  
\mathbb{R}j_{x}^{*}(F) \ar[r]^-{\mathbb{R}j_{x}^{*}(f)} & \mathbb{R}j_{x}^{*}(G) \\  
F_{s} \ar[u] \ar[r] & \bullet \ar[u]}$$  
where the morphism $* \longrightarrow \mathbb{R}j_{x}^{*}(G)$ is adjoint to the   
morphism $s : \mathbb{L}(j_{x})_{!}(*)\simeq h_{x} \longrightarrow G$.   
  
\begin{cor}\label{c13}  
Let $f : F \longrightarrow G$ be a morphism in $(C,S)^{\wedge}$. With the same  
notations as above, the morphism $f$ is a  
$\pi_{*}$-equivalence if and only if for any $s : h_{x} \longrightarrow G$ in  
$\mathrm{Ho}((C,S)^{\wedge})$,   
the induced morphism $F_{s} \longrightarrow *$ is a $\pi_{*}$-equivalence in $\mathrm{Ho}(((C,S)/x)^{\wedge})$.  
\end{cor}  
  
\begin{prop}\label{p6}  
Let $f : F \longrightarrow G$ be a $\pi_{*}$-equivalence in $(C,S)^{\wedge}$ and   
$F \longrightarrow F'$ be an objectwise cofibration (i.e. a monomorphism). Then, the induced morphism  
$$f' : F' \longrightarrow F'\coprod_{F}G'$$  
is a $\pi_{*}$-equivalence.  
\end{prop}  
  
\textit{Proof:} As $F \longrightarrow F'$ is an objectwise monomorphism,   
$F'\coprod_{F}G'$ is a homotopy coproduct in $SPr(C)$, and therefore in $(C,S)^{\wedge}$. One can therefore  
replace $F$, $G$ and $F'$ by their fibrant models in $(C,S)^{\wedge}$ and   
suppose therefore that they preserve equivalences. The proof is then the same as in \cite[Prop. $2.2$]{ja}.  
\hfill \textbf{$\Box$} 
  
\end{subsection}  
  
\begin{subsection}{The local model structure}  
  
The following result is completely similar to Theorem \ref{t3}, also as far as  
the proof is concerned. Therefore we will omit to repeat the complete   
proof below, only mentioning how to replace the set $H$ used in the proof of Theorem \ref{t3}.\\  
  
\begin{thm}\label{t5}  
Let $(C,S,\iota;\tau)$ be a pseudo-model site. There exists a closed model   
structure on $SPr(C)$, called the \emph{local projective model structure},   
for which the equivalences are the  
$\pi_{*}$-equivalences and the cofibrations are the cofibrations for the projective model  
structure on $(C,S)^{\wedge}$. Furthermore the local projective model structure is  
$\mathbb{U}$-combinatorial and left proper.   
  
The category $SPr(C)$ together with its local projective model structure will be denoted by $(C,S)^{\sim,\tau}$.  
\end{thm}  
  
\textit{Proof:} It is essentially the same as the proof of \ref{t3}. We will however give  
the set of morphism $H$ that one needs to use. We choose $\alpha$ to be a $\mathbb{U}$-small cardinal which is bigger than the   
cardinality of the set of morphisms in $C$ and than $\aleph_{0}$. Let $\beta$  
be a $\mathbb{U}$-small cardinal such that $\beta > 2^{\alpha}$.   
  
For a fibrant object $x \in C$,   
we consider a set $H_{\beta}(x)$, of representatives of the set of isomorphism classes of objects  
$F_{*} \longrightarrow h_{x}$ in $s(C,S)^{\wedge}/h_{x}$   
satisfying the following two conditions  
  
\begin{enumerate}  
  
\item The morphism  $F_{*} \longrightarrow h_{x}$ is a pseudo-representable hypercover  
in the sense of Definition \ref{dpseudorep}.  
  
\item For all $n\geq 0$, one has $\mathrm{Card}(F_{n})<\beta$.   
  
\end{enumerate}  
  
\noindent We set $H=\coprod_{x \in C^{f}}H_{\beta}(x)$, which is clearly a $\mathbb{U}$-small set.  
   
The main point of the proof is then to check that   
equivalences in the left Bousfield localization $L_{H}(C,S)^{\wedge}$ are exactly  
local equivalences.  
The argument follows exactly the main line of the proof of Theorem   
\ref{t3} and we leave details to the interested reader.  
\hfill $\Box$ \\  
  
The following corollaries and definitions are the same as the ones following Theorem \ref{t3}.  
  
\begin{cor}\label{c42}  
The model category $(C,S)^{\sim,\tau}$ is the left Bousfield localization of  
$(C,S)^{\wedge}$ with respect to the set of morphisms   
$$\left\{|F_{*}| \longrightarrow h_{x} | x \in Ob(C^{f}), \;  F_{*}\in \mathcal{H}_{\beta}(x) \right\}.$$  
\end{cor}  
  
\textit{Proof:} This is exactly the way we proved Theorem \ref{t5}. \hfill \textbf{$\Box$} \\  
  
\begin{cor}\label{c62}  
An object $F \in (C,S)^{\sim,\tau}$ is fibrant if and only if it is objectwise fibrant, preserves equivalences and satisfies the following \emph{hyperdescent} condition \begin{itemize}  
  \item  For any fibrant object $x\in C$ and any $H_{*}\in H_{\beta}(x)$,    
the natural morphism   
$$F(x)\simeq \mathbb{R}_{w}\underline{Hom}(h_{x},F) \longrightarrow \mathbb{R}_{w}\underline{Hom}(|H_{*}|,F)$$  
is an isomorphism in $\mathrm{Ho}(SSet)$.  
\end{itemize}  
\end{cor}  
  
\textit{Proof:} This follows from Corollary \ref{c42} and from the explicit description of fibrant objects in a   
left Bousfield localization  
(see \cite[Thm. 4.1.1]{hi}). \hfill \textbf{$\Box$} 
  
\begin{rmk} \emph{As we did in Remark \ref{twoinone}, we would like to stress here that the proof of Theorem   
\ref{t5} (i.e. of Theorem \ref{t3}) proves actaully both Theorem \ref{t5} and Corollary \ref{c42}, in that   
it gives \textit{two descriptions} of the same model category $(C,S)^{\sim,\tau}$:   
one as the left Bousfield localization of $(C,S)^{\wedge}$ with respect to \textit{local equivalences}   
and the other as the left Bousfield localization of the same $(C,S)^{\wedge}$ but this time   
with respect to \textit{hypercovers} (precisely, with respect to the set of   
morphisms defined in the statement of Corollary \ref{c42}).}  
\end{rmk}  
  
\begin{df}\label{d152}  
An object $F \in (C,S)^{\wedge}$ is said to \emph{have hyperdescent} (or $\tau$-hyperdescent if the topology  
$\tau$ has to be reminded)   
if for any fibrant object $x \in C$ and any pseudo-representable hypercover $H_{*} \longrightarrow h_{x}$, the induced  
morphism  
$$F(x)\simeq \mathbb{R}_{w}\underline{Hom}(h_{x},F)\longrightarrow    
\mathbb{R}_{w}\underline{Hom}(|H_{*}|,F)$$  
is an isomorphism in $\mathrm{Ho}(SSet_{\mathbb{U}})$.  
\end{df}  
  
From now on, we will adopt the following terminology and notations.  
  
\begin{df}\label{d162}  
  
Let $(C,S,\iota;\tau)$ be a pseudo-model site in $\mathbb{U}$.  
\begin{itemize}  
  
\item A \emph{stack} on $(C,S,\iota;\tau)$ is a pre-stack $F \in (C,S)^{\wedge}$ that has  
$\tau$-hyperdescent (Definition \ref{d152}).   
  
\item The model category $(C,S)^{\sim,\tau}$ is called the \emph{model category of stacks} on the   
pseudo-model site $(C,S,\iota;\tau)$. The category $\mathrm{Ho}((C,S)^{\wedge})$ (resp. $\mathrm{Ho}((C,S)^{\sim,\tau})$) is called the \emph{homotopy category of pre-stacks}, (resp. the \emph{homotopy category of stacks}). Objects of $\mathrm{Ho}((C,S)^{\wedge})$ (resp. $\mathrm{Ho}((C,S)^{\sim,\tau})$) will simply be called \emph{pre-stacks} on   
$(C,S,\iota)$ (resp., \emph{stacks} on $(C,S,\iota;\tau)$). The functor $a : \mathrm{Ho}((C,S)^{\wedge}) \longrightarrow   
\mathrm{Ho}((C,S)^{\sim,\tau})$ will be called the \emph{associated stack functor}.  
   
\item The topology $\tau$ is said to be \emph{sub-canonical} if for any $x \in C$ the  
pre-stack $\mathbb{R}\underline{h}_{x} \in \mathrm{Ho}((C,S)^{\wedge})$ is   
a stack (in other words if the Yoneda   
embedding $\mathbb{R}\underline{h}_{x} : \mathrm{Ho}(C,S) \longrightarrow   
\mathrm{Ho}((C,S)^{\wedge})$ factors through the subcategory  
of stacks).  
    
\item For pre-stacks $F$ and $G$ on $(C,S,\iota;\tau)$, we will denote by  
$\mathbb{R}_{w}\underline{Hom}(F,G) \in \mathrm{Ho}(SSet_{\mathbb{U}})$ (resp. by $\mathbb{R}_{w,\tau}\underline{Hom}(F,G) \in \mathrm{Ho}(SSet_{\mathbb{U}})$) 
the simplicial derived $Hom$-simplicial set computed in the simplicial model category $(C,S)^{\wedge}$ (resp. $(C,S)^{\sim,\tau}$).  
  
\end{itemize}  
\end{df}  
  
\medskip  
  
As $(C,S)^{\sim,\tau}$ is a left Bousfield localization of $(C,S)^{\wedge}$, the identity functor  
$(C,S)^{\wedge} \longrightarrow (C,S)^{\sim,\tau}$ is left Quillen and its right adjoint (which is still the  
identity functor) induces by right derivation a fully faithful functor  
$$j : \mathrm{Ho}((C,S)^{\sim,\tau}) \longrightarrow \mathrm{Ho}((C,S)^{\wedge}).$$  
Furthermore, the essential image of this inclusion functor is exactly the full subcategory consisting  
of objects having the hyperdescent property. The left adjoint    
$$a : \mathrm{Ho}((C,S)^{\wedge}) \longrightarrow \mathrm{Ho}((C,S)^{\sim,\tau})$$  
to the inclusion $j$, is a left inverse to $j$.  
  
We will finish this paragraph by the following proposition.  
  
\begin{prop}\label{p32}  
\begin{enumerate}  
\item Let $F$ and $G$ be two pre-stacks on $(C,S,\iota;\tau)$. If $G\!$ is a stack then the natural morphism  
$$\mathbb{R}_{w}\underline{Hom}(F,G) \longrightarrow \mathbb{R}_{w,\tau}\underline{Hom}(F,G)$$  
is an isomorphism in $\mathrm{Ho}(SSet)$.  
  
\item The functor $\mathrm{Id} : (C,S)^{\wedge} \longrightarrow (C,S)^{\sim,\tau}$ preserves   
homotopy fibered products.  
  
\end{enumerate}  
\end{prop}  
  
\textit{Proof:} $(1)$ follows formally from Corollary \ref{c42} while $(2)$ follows from Corollary \ref{c12}.   
\hfill \textbf{$\Box$} 
  
\end{subsection}  
  
\begin{subsection}{Comparison between the $S$-theory and the pseudo-model theory}  
  
In this Subsection, we fix a pseudo-model category $(C,S,\iota)$ in $\mathbb{U}$,   
together with a pre-topology $\tau$ on it. The natural induced topology on $\mathrm{Ho}(C,S)$ will be denoted  
again by $\tau$.   
We let $T$ be $L(C,S)$, the simplicial localization of $(C,S)$ along the set $S$ of its equivalences.   
As $\mathrm{Ho}(T)=\mathrm{Ho}(C,S)$ (though the two $\mathrm{Ho}(-)$'s here have  
different meanings), the topology $\tau$ may also be considered as an $S$-topology on $T$.   
Therefore, we have  
on one side a pseudo-model site $(C,S,\iota; \tau)$, and on the other side  
an $S$-site $(T,\tau)$, and we wish to compare the two corresponding model categories of stacks.   
  
\begin{thm}\label{t6}  
The two model categories $(C,S)^{\sim,\tau}$ and $SPr_{\tau}(T)$ are Quillen equivalent.  
\end{thm}  
  
\textit{Proof:} By Theorem \ref{t2}, the model categories of pre-stacks   
$SPr(T)$ and $(C,S)^{\wedge}$ are Quillen equivalent.  
Furthermore, it is quite clear that through this equivalence the notions of local equivalences in   
$SPr(T)$ and $(C,S)^{\wedge}$ coincide.  
As the local model structures are both left Bousfield localizations with respect  
to local equivalences, this shows that this Quillen equivalence between $(C,S)^{\wedge}$ and  
$SPr(T)$ induces a Quillen equivalence on the model categories of stacks. \hfill \textbf{$\Box$} \\  
\noindent Then, corollaries \ref{cc2} and \ref{c7} imply the following\\  
   
\begin{cor}\label{c15}  
\begin{enumerate}  
\item The model category $(C,S)^{\sim,\tau}$ is a $t$-complete $\mathbb{U}$-model topos.   
  
\item The homotopy category $\mathrm{Ho}((C,S)^{\sim,\tau})$ is internal.  
  
\item  
There exists an isomorphism of $S$-categories in $\mathrm{Ho}(S-Cat_{\mathbb{U}})$  
$$LSPr_{\tau}(T)\simeq L(C,S)^{\sim,\tau}.$$  
  
\end{enumerate}  
\end{cor}  
  
\medskip  
  
Now we want to compare the two Yoneda embeddings (the simplicial one and the   
pseudo-model one). To do this, let us suppose now that the topology $\tau$ is  
\textit{sub-canonical} so that the two Yoneda embeddings factor through   
the embeddings of the homotopy categories of stacks:  
$$\mathbb{R}\underline{h} : \mathrm{Ho}(C,S) \longrightarrow \mathrm{Ho}((C,S)^{\sim,\tau}) \qquad  
L\underline{h} : \mathrm{Ho}(T) \longrightarrow \mathrm{Ho}(Int(SPr_{\tau}(T)))\simeq \mathrm{Ho}(SPr_{\tau}(T)).$$  
One has $\mathrm{Ho}(C,S)=\mathrm{Ho}(T)$, and Corollary \ref{c15}  
gives an equivalence of categories between   
$\mathrm{Ho}(SPr_{\tau}(T))$ and $\mathrm{Ho}((C,S)^{\sim,\tau})$.   
\begin{cor}\label{c16}  
The following diagram commutes up to an isomorphism  
$$\xymatrix{  
\mathrm{Ho}(C,S) \ar[d]_-{\sim} \ar[r]^-{\mathbb{R}\underline{h}} & \mathrm{Ho}((C,S)^{\sim,\tau}) \ar[d]^-{\sim} \\  
\mathrm{Ho}(T) \ar[r]_-{L\underline{h}} & \mathrm{Ho}(SPr_{\tau}(T)).}$$  
\end{cor}  
  
\textit{Proof:} This follows from the fact that for any $x\in M$, one has natural isomorphisms  
$$[\mathbb{R}\underline{h}_{x},F]_{\mathrm{Ho}((C,S)^{\sim,\tau})}  
\simeq F(x)\simeq [L\underline{h}_{x},F]_{\mathrm{Ho}(SPr_{\tau}(T))}.$$  
This implies that $\mathbb{R}\underline{h}_{x}$ and $L\underline{h}_{x}$ are naturally isomorphic  
as objects in $\mathrm{Ho}((C,S)^{\wedge})$.  \hfill \textbf{$\Box$} 
  
\end{subsection}  
  
\begin{subsection}{Functoriality}\label{Mfunctoriality}  
  
In this subsection, we state and prove in detail the functoriality results and some useful criteria for  
continuous morphisms and continuous equivalences between pseudo-model sites, in such a way that the   
reader only interested in working with stacks over pseudo-model sites will find here a more or less   
self-contained treatment. However, at the end of the subsection and in occasionally scattered   
remarks, we will also mention the comparison between functoriality on pseudo-model sites and   
the corresponding functoriality on the associated Dwyer-Kan localization S-sites.  
  
Recall from Subsection \ref{Mprestacks} (or Subsection \ref{resdia} before Corollary \ref{c0}) that if $(C,S)$ and $(C',S')$ are categories with a distinguished subset of morphisms (e.g.,   
pseudo-model categories) and $f : C \rightarrow C'$ is a functor sending $S$ into $S'$, we have a Quillen adjunction   
$$f_{!} : (C,S)^{\wedge} \longrightarrow (C,S')^{\wedge} \qquad (C,S)^{\wedge}   
\longleftarrow (C',S')^{\wedge} : f^{*}$$  
  
If $(C,S,\iota)$ is a pseudo-model category, by Proposition \ref{junction}, we have in particular the following Quillen equivalences  
$$i_{c}^{*} : (C,S)^{\wedge} \simeq ((C,S)^{c})^{\wedge} \qquad  
i_{f}^{*} : (C,S)^{\wedge} \simeq ((C,S)^{f})^{\wedge} $$  
$$i_{cf}^{*} : (C,S)^{\wedge} \simeq ((C,S)^{cf})^{\wedge},$$  
which will be useful to establish functorial properties of the homotopy category  
$\mathrm{Ho}((C,S)^{\wedge})$. Indeed, if $f : (C,S) \longrightarrow (C',S')$ is a functor  
\textit{such that} $f(S^{cf})\subset S'$ (e.g. a left or right Quillen functor), then $f$ induces well defined functors  
$$\mathbb{R}f^{*} : \mathrm{Ho}((C',S')^{\wedge}) \longrightarrow  
\mathrm{Ho}(((C,S)^{cf})^{\wedge})\simeq \mathrm{Ho}((C,S)^{\wedge}),$$  
$$\mathbb{L}f_{!} : \mathrm{Ho}((C,S)^{\wedge})\simeq \mathrm{Ho}(((C,S)^{cf})^{\wedge}) \longrightarrow  
\mathrm{Ho}((C',S')^{\wedge}).$$  
The (derived) \textit{inverse image} functor $\mathbb{R}f^{*}$ is clearly right adjoint to the (derived) \textit{direct image} functor $\mathbb{L}f_{!}$.

The reader should be warned that the direct and inverse image functors are not, in general, functorial in $f$.
However, the following proposition ensures in many cases the functoriality of these  constructions.  
  
\begin{prop}\label{p1p1}  
Let $(C,S)$, $(C',S')$ and $(C'',S'')$ be pseudo-model categories and  
$$\xymatrix{(C,S) \ar[r]^-{f} & (C',S') \ar[r]^-{g} & (C'',S'')}$$  
be two functors preserving fibrant or cofibrant objects and equivalences between them.   
Then, there exist natural isomorphisms  
$$\mathbb{R}(g\circ f)^{*}\simeq \mathbb{R}f^{*}\circ \mathbb{R}g^{*} :  
\mathrm{Ho}((C'',S'')^{\wedge}) \longrightarrow \mathrm{Ho}((C,S)^{\wedge}),$$  
$$\mathbb{L}(g\circ f)_{!}\simeq \mathbb{L}g_{!}\circ \mathbb{L}f_{!} :  
\mathrm{Ho}((C,S)^{\wedge}) \longrightarrow \mathrm{Ho}((C'',S'')^{\wedge}).$$  
These isomorphisms are furthermore associative and unital in the arguments $f$ and $g$.  
\end{prop}  
  
\textit{Proof:} The proof is the same as that of the usual property of  
composition for derived Quillen functors (see \cite[Thm. $1.3.7$]{ho}), and is left to the reader. \hfill $\Box$ \\  
  
Examples of pairs of functors to which the previous proposition applies are given by pairs of right or left Quillen functors. 
  
\begin{prop}\label{p2p2}  
If $f : (C,S) \longrightarrow (C,S)$ is a (right or left) Quillen equivalence between  
pseudo-model categories, then the induced functors  
$$\mathbb{L}f_{!} : \mathrm{Ho}((C,S)^{\wedge}) \longrightarrow \mathrm{Ho}((C',S')^{\wedge}) \qquad  
\mathrm{Ho}((C,S)^{\wedge}) \longleftarrow \mathrm{Ho}((C',S')^{\wedge}) : \mathbb{R}f^{*},$$  
are equivalences, quasi-inverse of each others.  
\end{prop}  
  
\textit{Proof:} This is a straightforward application of Corollary \ref{c0}.   
\hfill $\Box$ \\  
  
Let $(C,S)$ and $(C',S')$ be pseudo-model categories and let us consider a functor   
$f : C \longrightarrow C'$ such that $f(S^{cf})\subset S'$. We will denote by   
$f_{cf} : (C,S) \longrightarrow (C',S')$ the composition  
$$f_{cf} : \xymatrix{  
(C,S)\ar[r]^-{RQ} & (C,S)^{cf} \ar[r]^-{f} & (C',S'),}$$  
\noindent where $R$ (respectively, $Q$) denotes the fibrant (resp., cofibrant) replacement functor in $(C,S)$. We deduce an adjunction on the model categories of pre-stacks  
$$(f_{cf})_{!} : (C,S)^{\wedge} \longrightarrow (C',S')^{\wedge} \qquad  
(C,S)^{\wedge} \longleftarrow (C',S')^{\wedge} : f_{cf}^{*}.$$  
Note that the right derived functor $\mathbb{R}f_{cf}^{*}$   
is isomorphic to the functor $\mathbb{R}f^{*}$ defined above.   
  
\begin{prop} \label{cont}  
Let $(C,S;\tau)$ and $(C',S';\tau')$ be pseudo-model sites and $f : C \longrightarrow C'$ a functor such   
that $f(S^{cf})\subset S'$. Then the following properties are equivalent.  
\begin{enumerate}  
\item The right derived functor $\mathbb{R}f_{cf}^{*}\simeq\mathbb{R}f^{*}:\mathrm{Ho}((C',S')^{\wedge})  
\rightarrow \mathrm{Ho}((C,S)^{\wedge})$ sends the subcategory $\mathrm{Ho}((C',S')^{\sim, \tau'})$   
into the subcategory $\mathrm{Ho}((C,S)^{\sim, \tau})$.  
\item If $F \in (C',S')^{\wedge}$ has $\tau'$-hyperdescent, then $f^{*}F \in SPr(C)$ has $\tau$-hyperdescent.  
  
\item For any pseudo-representable hypercover $H_{*} \longrightarrow h_{x}$ in   
$(C,S)^{\wedge}$ (see Def. \ref{dpseudorep}), the morphism  
$$\mathbb{L}(f_{cf})_{!}(H_{*}) \longrightarrow \mathbb{L}(f_{cf})_{!}(h_{x})\simeq h_{f_{cf}(x)}$$  
is a local equivalence in $(C',S')^{\wedge}$.  
  
\item The functor $f_{cf}^{*}:(C',S')^{\sim, \tau} \longrightarrow (C,S)^{\sim, \tau}$ is right Quillen.  
\end{enumerate}  
\end{prop}  
  
\textit{Proof:} The equivalence between $(1)$, $(2)$ and $(3)$  
follows immediately from the fact that fibrant objects in   
$(C,S)^{\sim, \tau}$ (resp. in $(C',S')^{\sim, \tau}$) are exactly those fibrant objects in   
$(C,S)^{\wedge}$ (resp. in $(C',S')^{\wedge}$) which satisfy $\tau$-hyperdescent (resp.   
$\tau'$-hyperdescent) (see Corollary \ref{c62}). Finally, $(4)$ and $(2)$ are   
equivalent by adjunction.  \hfill \textbf{$\Box$}   
   
\begin{df}\label{dfcont}  
Let $(C,S;\tau)$ and $(C',S';\tau')$ be pseudo-model sites. A functor $f:C \rightarrow C'$   
such that $f(S^{cf})\subseteq S'$, is  said to be \emph{continuous} or \emph{a morphism of   
pseudo-model sites}, if it satisfies one of the equivalent conditions of Proposition \ref{cont}.  
\end{df}  
  
\begin{rmk}  
\emph{By the comparison Theorem \ref{t6}, a functor $f:(C,S;\tau) \rightarrow   
(C',S';\tau')$ such that $f(S^{cf})\subseteq S'$, is continuous if and only if the induced   
functor $(L(C,S),\tau)\simeq (L(C^{cf},S^{cf}),\tau)\rightarrow (L(C',S'),\tau')$ between   
the simplicially localized associated $S$-sites is continuous according to Definition \ref{d17}.}  
\end{rmk}  
  
\medskip  
  
It is immediate to check that if $f$ is a continuous functor, then  
the functor  
$$\mathbb{R}f^{*} : \mathrm{Ho}((C',S')^{\sim,\tau'}) \longrightarrow  
\mathrm{Ho}((C,S)^{\sim,\tau})$$  
has as left adjoint  
$$\mathbb{L}(f_{!})^{\sim}\simeq\mathbb{L}(f_{cf\,!}) : \mathrm{Ho}((C,S)^{\sim,\tau}) \longrightarrow  
\mathrm{Ho}((C',S')^{\sim,\tau'}),$$  
the functor defined by the formula  
$$\mathbb{L}(f_{!})^{\sim}(F):=a(\mathbb{L}f_{!}(F)),$$  
for $F \in \mathrm{Ho}((C,S)^{\sim,\tau}) \subset  
\mathrm{Ho}((C,S)^{\wedge})$, where $a:\mathrm{Ho}((C,S)^{\wedge})\rightarrow   
\mathrm{Ho}((C,S)^{\sim,\tau})$ is the associated stack functor.  
  
The basic properties of the associated stack functor $a$ imply that the functoriality result of Proposition \ref{p1p1} still holds by replacing the model categories of pre-stacks with the model categories of stacks, if $f$ and $g$ are continuous. 
  
Now we define the obvious notion of continuous equivalence between pseudo-model sites.  
  
\begin{df}\label{conteq}  
A continuous functor $f:(C,S;\tau) \rightarrow (C',S';\tau')$ is said to be a \emph{continuous   
equivalence} or an \emph{equivalence of pseudo-model sites} if the induced right Quillen functor   
$f_{cf}^{*}:(C',S')^{\sim, \tau'}\rightarrow (C,S)^{\sim, \tau}$ is a Quillen equivalence.  
\end{df}  
  
The following criterion will be useful in the next section. 
  
\begin{prop}\label{critconteq}  
Let $(C,S;\tau)$ and $(C',S';\tau')$ be pseudo-model sites, $f : C \longrightarrow C'$ a functor    
such that $f(S^{cf})\subseteq S'$ and   
$f_{cf} : (C,S)\longrightarrow (C',S')$   
the induced functor. Let us denote by $\tau$ (resp. by $\tau'$) the induced   
Grothendieck topology on the homotopy category $\mathrm{Ho}(C,S)$ (resp. $\mathrm{Ho}(C',S')$). Suppose that  
\begin{enumerate}  
\item The induced morphism $Lf_{cf} : L(C,S) \longrightarrow L(C',S')$   
between the Dwyer-Kan localizations is an equivalence of $S$-categories.  
\item The functor   
$$\mathrm{Ho}(f_{cf}) : \mathrm{Ho}(C,S) \longrightarrow \mathrm{Ho}(C',S')$$  
reflects covering sieves (i.e., a sieve $R$ over $x \in   
\mathrm{Ho}(C,S)$ is $\tau$-covering iff the sieve generated by $\mathrm{Ho}(f_{cf})(R)$   
is a $\tau'$-covering sieve over $f_{cf}(x)$.   
\end{enumerate}  
Then $f$ is a continuous equivalence.  
\end{prop}  
  
\textit{Proof:} This follows easily from the comparison statement Theorem \ref{t6} and from Theorem \ref{t1}.  
\hfill \textbf{$\Box$}  
  
\end{subsection}  
  
\begin{subsection}{A Giraud's theorem for model topoi}

In this section we prove a Giraud's type theorem characterizing
model topoi internally. Applied to $t$-complete model topoi, this
will give an internal description of model categories that 
are Quillen equivalent to some model category 
of stacks over an $S$-site. We like to consider
this result as an extension of D. Dugger 
characterization of combinatorial model categories
(\cite{du2}), and
as a model category analog of J. Lurie's theorem
characterizing $\infty$-topoi (see \cite[Thm. 2.4.1]{lu}). 
Using the strictification theorem 
of A. Hirschowitz and C. Simpson (stated in \S 4.2 of \cite{msri}) it also 
gives a proof of the Giraud's theorem for Segal topoi
conjectured in \cite[Conj. 5.1.1]{msri}. The statement presented here
is very close in spirit to the statement presented
in \cite{rez2}, with some minor differences in that our conditions
are weaker than \cite{rez2}, and closer to the original ones stated by Giraud (see \cite{sga4} Exp. IV, Th\'eor\`eme 1.2). \\

We start with some general definitions.

\begin{df}\label{dgiraud}
Let $M$ be any $\mathbb{U}$-model category. 
\begin{enumerate}
\item The model category has
\emph{disjoint homotopy coproducts} if for any
$\mathbb{U}$-small family of objects
$\{x_{i}\}_{i\in I}$, and any $i\neq j$ in $I$, the following 
square is homotopy cartesian
$$\xymatrix{
\emptyset \ar[r] \ar[d] & x_{i} \ar[d] \\
x_{j} \ar[r] & \coprod^{\mathbb{L}}_{i\in I}x_{i}.}$$

\item \emph{The homotopy colimits are stable under pullbacks in $M$}
if for any morphism $y\longrightarrow z$ in $M$, such that 
$z$ is fibrant, and any
$\mathbb{U}$-small diagram $x_{*} : I \longrightarrow M/z$
of objects over $z$, the natural morphism
$$\mathrm{hocolim}_{i\in I}(x_{i}\times^{h}_{z}y)\longrightarrow
(\mathrm{hocolim}_{i\in I}x_{i})\times^{h}_{z}y$$
is an isomorphism in $\mathrm{Ho}(M)$. 

\item A \emph{Segal groupoid object in $M$} is a 
simplicial object 
$$X_{*} : \Delta^{op} \longrightarrow M,$$
such that 
\begin{itemize}
\item for any $n>0$, the natural morphism
$$X_{n} \longrightarrow \underbrace{X_{1}\times^{h}_{X_{0}}X_{1}\times^{h}_{X_{0}}
\dots \times^{h}_{X_{0}}X_{1}}_{n\; times}$$
induced by the $n$ morphisms $s_{i} : [1] \longrightarrow [n]$,
defined as $s_{i}(0)=i$, $s_{i}(1)=i+1$, is an 
isomorphism in $\mathrm{Ho}(M)$.
\item The morphism
$$d_{0}\times d_{1} : X_{2} \longrightarrow X_{1}\times^{h}_{d_{0},X_{0},d_{0}}X_{1}$$
is an equivalence in $\mathrm{Ho}(M)$.

\end{itemize}

\item We say that \emph{Segal equivalences relation are homotopy effective in $M$} if
for any Segal groupoid object $X_{*}$ in $M$ with homotopy colimit
$$|X_{*}|:=\mathrm{hocolim}_{n\in \Delta}X_{n},$$
and any $n>0$, the
natural morphism
$$X_{n} \longrightarrow 
\underbrace{X_{0}\times^{h}_{|X_{*}|}X_{0}\times^{h}_{|X_{*}|}
\dots \times^{h}_{|X_{*}|}X_{0}}_{n\; times}$$
induced by the $n$ distinct morphisms $[0] \rightarrow [n]$,
is an isomorphism in $\mathrm{Ho}(M)$. 

\end{enumerate}
\end{df}

We are now ready to state our version of Giraud's theorem
for model topoi.

\begin{thm}\label{tgiraud}
Let $M$ be a $\mathbb{U}$-combinatorial model
category (see Definition \ref{db3}). Then, $M$ is 
a $\mathbb{U}$-model topos if and only
if it satisfies the following conditions.
\begin{enumerate}
\item $M$ has disjoint homotopy coproducts.

\item Homotopy colimits in $M$ are stable under homotopy 
pullbacks.

\item Segal equivalence relations are homotopy effective in $M$.

\end{enumerate}
\end{thm}

\textit{Proof:} The fact that the conditions are 
satisfied in any model topos follows easily from
the well known fact that they are satisfied 
in the model category $SSet$.
The hard point is to 
prove they are sufficient conditions. \\
Let $M$ be a $\mathbb{U}$-model category satisfying 
the conditions of the theorem. 

We chose 
a regular cardinal $\lambda$ as in the proof of 
\cite[Prop. 3.2]{du2}, and let $C:=M_{\lambda}$ be a 
$\mathbb{U}$-small full sub-category of 
$M$ consisting of a set of representatives
of $\lambda$-small objects in $M$.
By increasing $\lambda$ if necessary, one
can assume that the full sub-category
$C$ of $M$ is $\mathbb{U}$-small, and
is stable under fibered products in $M$ and under the
fibrant and cofibrant replacement functors (let us
suppose these are fixed once for all). By this last
condition we mean that for any morphism $x \rightarrow y$
in $C$, the functorial 
factorizations $x \rightarrow x' \rightarrow y$
are again in $C$.
Let $\Gamma_{*}$ and $\Gamma^{*}$ be
fibrant and cofibrant resolution functors on $M$ (\cite{hi}, Ch. 16).
We can also assume that $C$ is stable by $\Gamma_{*}$
and $\Gamma^{*}$ (i.e. that 
for any $x\in C$ and any $[n]\in \Delta$, 
$\Gamma_{n}(x)$ and $\Gamma^{n}(x)$ belong to $C$).
We note that $C$ is not stricly speaking a 
pseudo-model category but will 
behave pretty much the same way. 

We consider
the functor
$$\underline{h}^{C} : M \longrightarrow SPr(C),$$
sending an object $x\in M$ to the simplicial 
presheaf
$$\begin{array}{cccc}
\underline{h}^{C}_{x} : & C^{op} & \longrightarrow & SSet_{\mathbb{U}} \\
 & y & \longmapsto & Hom(\Gamma^{*}(y),x). 
\end{array}$$
The functor $\underline{h}$ has a left adjoint
$$L : SPr(C) \longrightarrow M,$$
sending a $\mathbb{U}$-simplicial presheaf $F$ to 
its geometric realization with respect to 
$\Gamma$. By the standard properties of mapping
spaces, one sees that for any fibrant object $x\in M$
the simplicial presheaf $\underline{h}^{C}_{x}$ 
is fibrant in the model category of restricted diagrams
$(C,W)^{\wedge}$. This, and the
general properties of left Bousfield
localizations imply that 
the pair $(\underline{h}^{C},L)$ 
defines a Quillen adjunction
$$L : (C,W)^{\wedge} \longrightarrow M \qquad
(C,W)^{\wedge} \longleftarrow M : \underline{h}^{C}.$$

\begin{lem}\label{lgiraud1}
The right derived functor
$$\mathbb{R}\underline{h}^{C} : \mathrm{Ho}(M) \longrightarrow \mathrm{Ho}((C,W)^{\wedge})$$
is fully faithful. 
\end{lem}

\textit{Proof:} By the choice of $C$, any object $x\in M$ is
a $\lambda$-filtered colimit $x\simeq \mathrm{colim}_{i\in I}x_{i}$
of objects $x_{i}\in C$. As all objects in $C$ are 
$\lambda$-small, this implies that 
$$\mathbb{R}\underline{h}^{C}_{x}\simeq 
\mathrm{hocolim}_{i\in I}\mathbb{R}\underline{h}^{C}_{x_{i}}.$$
From this, one sees that to prove that $\mathbb{R}\underline{h}^{C}$
is fully faithful, it is enough to prove it is
fully faithful when restricted to 
objects of $C$. This last case can be treated exactely as
in the proof of our Yoneda Lemma \ref{t10}. \hfill $\Box$ \\

By the previous lemma and 
by Proposition 3.2 of \cite{du2}, we can conclude that 
there is a $\mathbb{U}$-small set of morphisms
$S$ in $(C,W)^{\wedge}$ such that the above adjunction induces
a Quillen equivalence
$$L : L_{S}(C,W)^{\wedge} \longrightarrow M \qquad
L_{S}(C,W)^{\wedge} \longleftarrow M : \underline{h}^{C}.$$
By Corollary \ref{c7} $(2)$, it only remains to show that 
the left Bousfield localization
of $(C,W)^{\wedge}$ along $S$ is exact, 
or equivalently  that the functor 
$\mathbb{L}L$ commutes with homotopy pull backs. \\

We start by the following particular 
case. Let $c \in C$ and $h_{c}$ be the
presheaf represented by $c$. One can see
$h_{c}$ as an object in $(C,W)^{\wedge}$
by considering it as a presheaf of
discrete simplicial sets. Let $F \longrightarrow h_{c}$
and $G\longrightarrow h_{c}$ be two 
morphisms in $(C,W)^{\wedge}$. 

\begin{lem}\label{lgiraud2}
The natural morphism
$$\mathbb{L}L(F\times^{h}_{h_{c}}G) \longrightarrow
\mathbb{L}L(F)\times^{h}_{\mathbb{L}L(h_{c})}\mathbb{L}L(G)$$
is an isomorphism in $\mathrm{Ho}(M)$.
\end{lem}

\textit{Proof:} Up to an equivalence, we can write $F$
as a homotopy colimit $\mathrm{hocolim}_{i\in I}h_{x_{i}}$ for
some $x_{i} \in C$. As homotopy pull-backs commutes
with homotopy colimits this shows that one can
suppose $F$ and $G$ of the form $h_{a}$ and $h_{b}$, 
for $a$ and $b$ two objects in $C$. 

Now, as in lemma \ref{strange}, one checks that
$h_{x}$ and $\mathbb{R}\underline{h}^{C}_{x}$
are naturally isomorphic in $Ho((C,W)^{\wedge})$. For this, 
we easily deduce that the natural morphism
$$h_{a}\times^{h}_{h_{c}}h_{b}\longrightarrow h_{a\times^{h}_{c}b},$$
is an equivalence in $(C,W)^{\wedge}$ (here 
$h_{a\times^{h}_{c}b}$ can be seen as an object of $C$ 
because of our stabilty assumptions). Therefore, to prove the lemma
it is enough to check that for any $x \in C$ the natural morphism
$h_{x} \longrightarrow \underline{h}^{C}(x)$ induces
by adjunction a morphism $L(h_{x}) \longrightarrow x$ which is
an equivalence in $M$. But, as $h_{x}$ is always a cofibrant
object in $(C,W)^{\wedge}$, one has
$$L(h_{x})\simeq \mathbb{L}L(h_{x})
\simeq \mathbb{L}L(\underline{h}^{C}_{x})\simeq x$$
by lemma \ref{lgiraud1}. \hfill $\Box$ \\

Let $\coprod_{i\in I}h_{c_{i}}$ be a
coproduct with $c_{i} \in C$, and
$$F \longrightarrow \coprod_{i\in I}h_{c_{i}} \longleftarrow
G$$
be two morphisms in $(C,W)^{\wedge}$.

\begin{lem}\label{lgiraud3}
The natural morphism
$$\mathbb{L}L(F\times^{h}_{\coprod_{i\in I}h_{c_{i}}}G) \longrightarrow
\mathbb{L}L(F)\times^{h}_{\mathbb{L}L(\coprod_{i\in I}h_{c_{i}})}\mathbb{L}L(G)$$
is an isomorphism in $\mathrm{Ho}(M)$.
\end{lem}

\textit{Proof:} As for lemma \ref{lgiraud2}, one can reduce to the
case where $F$ and $G$ are of the form $h_{a}$ and $h_{b}$. 
Lemma \ref{lgiraud3} will then follows easily from
our assumption $(1)$ on $M$. \hfill $\Box$ \\

We are now ready to treat the general case. 

\begin{lem}\label{lgiraud4}
The functor $\mathbb{L}L$ preserves homotopy pull-backs.
\end{lem}

\textit{Proof:} Let $\xymatrix{F \ar[r] & H & \ar[l] G}$
be two morphisms in $(C,W)^{\wedge}$. One can, as for
lemma \ref{lgiraud2} suppose that $F$ and $G$ are of the form
$h_{a}$ and $h_{b}$. We can also suppose that $H$ is fibrant
in $(C,W)^{\wedge}$. 

We let $\coprod_{i}h_{x_{i}} \longrightarrow H$
be an epimorphism of simplicial presheaves with
$x_{i}\in C$, and we replace it by an equivalent fibration  
$p : X_{0} \longrightarrow H$.
We set $X_{*}$ the nerve of $p$, which is the
simplicial object of $(C,W)^{\wedge}$ given by
$$X_{n}:=
\underbrace{X_{0}\times_{H}X_{0}\times_{H}\dots \times_{H}X_{0}}_{n\; times},$$
and for which faces and degeneracies are given by the various
projections and generalized diagonals. As $p$ is a fibration between
fibrant objects one sees that 
$X_{*}$ is a Segal groupoid object in $(C,W)^{\wedge}$. 
Furthermore, as $p$ is homotopycally surjective
(as a morphism of simplicial presheaves), the natural morphism
$$|X_{*}| \longrightarrow H$$
is an equivalence in $(C,W)^{\wedge}$. Finally, 
as $X_{0}$ is equivalent to $\coprod_{i}h_{x_{i}}$, lemma
\ref{lgiraud3} implies that $\mathbb{L}L(X_{*})$ is 
a Segal groupoid object in $M$, and one has
$|\mathbb{L}L(X_{*})|\simeq \mathbb{L}L(H)$ as
$L$ is left Quillen.
The assumption $(3)$ on $M$
implies that 
$$\mathbb{L}L(X_{0}\times^{h}_{H}X_{0})
\simeq \mathbb{L}L(X_{1})\simeq \mathbb{L}L(X_{0})\times^{h}_{\mathbb{L}L(H)}
\mathbb{L}L(X_{0}).$$
To finish the proof of lemma \ref{lgiraud4} it is then enough to notice 
that since $X_{0} \longrightarrow H$ is surjective up to homotopy, 
the morphisms $h_{a},h_{b} \longrightarrow H$ can be lifted
up to homotopy to morphisms to $X_{0}$ (because they
correspond to elements in $H(a)$ and $H(b)$), and therefore
$$h_{a}\times^{h}_{H}h_{b}\simeq
h_{a}\times^{h}_{X_{0}}(X_{0}\times^{h}_{H}X_{0})\times_{H}^{h}h_{b}.$$
One can then apply lemma \ref{lgiraud3}. \hfill $\Box$ \\

Theorem \ref{tgiraud} is proven. \hfill $\Box$ \\

The following corollary is an internal classification
of $t$-complete model topoi.

\begin{cor}\label{cgiraud}
Let $M$ be a $\mathbb{U}$-combinatorial model category. Then 
the following are equivalent. 
\begin{enumerate}
\item The model category $M$ satisfies the conditions
of theorem \ref{tgiraud} and is furthermore
$t$-complete.
\item There exists a $\mathbb{U}$-small 
$S$-site $T$ such that $M$ is Quillen 
equivalent to $SPr_{\tau}(T)$. 
\end{enumerate}
\end{cor}

\textit{Proof:} $(1)$ and $(2)$ follow from
Theorem \ref{tgiraud} combined with 
our Theorem \ref{t4}.  \hfill $\Box$ \\

From the proof of theorem \ref{tgiraud} one 
also extracts the following consequence.

\begin{cor}\label{cgiraud1'}
Let $M$ be a $\mathbb{U}$-combinatorial model category. Then 
the following are equivalent. 
\begin{enumerate}
\item The model category $M$ satisfies the conditions
of theorem \ref{tgiraud} and is furthermore
$t$-complete.
\item There exists a $\mathbb{U}$-model category $N$, and a 
$\mathbb{U}$-small full subcategory of
cofibrant object $C \subset N^{c}$, and 
a topology $\tau$ on $\mathrm{Ho}(C):=(W\cap C)^{-1}C$, 
such that $M$ is Quillen 
equivalent to $(C,W)^{\sim,\tau}$. Furthermore, 
the natural functor $\mathrm{Ho}(C) \longrightarrow \mathrm{Ho}(N)$
is fully faithful and its image is stable
under homotopy pull backs.
\end{enumerate}
\end{cor}

This last corollary states that $M$ is Quillen 
equivalent to the model category of stacks over
something which is ``almost'' a pseudo-model site. However, the sub-category $C$ 
produced during the proof of Theorem \ref{tgiraud} is not
a pseudo-model site as it is not stable by equivalences in $N$.
On the other hand, one can show that
the closure $\overline{C}$ of $C$ by equivalences in $N$ is a pseudo-model
site, and that the natural morphism $LC \longrightarrow L\overline{C}$
is an equivalence of $S$-categories. 

\begin{cor}\label{cgiraud2}
If $M$ is a $\mathbb{U}$-model topos 
(resp. a $t$-complete $\mathbb{U}$-model topos) then so is
$M/x$ for any fibrant object $x \in M$.
\end{cor}

\textit{Proof:} Indeed, if $M$ is a $\mathbb{U}$-combinatorial 
model category satisfying the conditions of Theorem \ref{tgiraud}
then so does $M/x$ for any fibrant object $x$. Furthermore,
one can check that 
for any $S$-site $(T,\tau)$, and any 
object $F$ the model category $SPr_{\tau}(T)/F$
is $t$-complete. This implies that if $M$ 
is furthermore $t$-complete then so is $M/x$.  \hfill $\Box$ \\

\begin{cor}\label{cgiraud3}
\begin{enumerate}

\item Any $\mathbb{U}$-model topos $M$ is Quillen equivalent to 
a left proper model category for which avery object
is cofibrant and which is furthermore internal (i.e. is 
a symmetric monoidal model category for the direct product
moniodal structure). 
\item For any $\mathbb{U}$-model topos $M$ and any fibrant object
$x\in M$, the category $Ho(M/x)$ is
cartesian closed. 

\end{enumerate}
\end{cor}

\textit{Proof:} It is enough to check this for $M=L_{S}SPr(T)$, for 
some $\mathbb{U}$-small $S$-category $T$ and
some 
$\mathbb{U}$-small set of morphisms $S$ in $SPr(T)$ such that 
$Id : SPr(T) longrightarrow L_{S}SPr(T)$ preserves
homotopy fiber products. We can also replace 
the projective model structure $SPr(T)$ by the injective one
$SPr_{inj}(T)$ (see Prop. \ref{pp1}), and therefore can suppose 
$M$ of the form $L_{S}SPr_{inj}(T)$, again with
$Id : SPr_{inj}(T) longrightarrow L_{S}SPr_{inj}(T)$ preserving
homotopy fiber products. We know that 
$SPr_{inj}(T)$ is an internal model category in which 
every object is cofibrant, and from this one
easily deduces that the same is true for the exact localization
$L_{S}SPr_{inj}(T)$. \\

$(2)$ follows from $(1)$ and Cor. \ref{cgiraud2}. \hfill $\Box$ 

\end{subsection}

\end{section}  
  
\begin{section}{\'{E}tale $K$-theory of commutative $\mathbb{S}$-algebras}  
  
In this section we apply the theory of stacks over pseudo-model sites   
developed in the previous section to the problem of defining a notion of \'{e}tale   
$K$-theory of a commutative $\mathbb{S}$-algebra i.e. of a commutative monoid in   
Elmendorf-Kriz-Mandell-May's category of $\mathbb{S}$-modules (see \cite{ekmm}).   
The idea is very simple. We only need two ingredients: the first is a notion of an   
\textit{\'etale topology} on the model category $(\mathrm{Alg}_{\mathbb{S}})$ of   
commutative $\mathbb{S}$-algebras and the second is the corresponding model category of   
\textit{\'etale stacks} on $(\mathrm{Alg}_{\mathbb{S}})$. Then, in analogy with the   
classical situation (see \cite[\S 3]{ja}), \textit{\'etale $K$-theory} will be just defined as a   
fibrant replacement of algebraic $K$-theory in the category of \'etale stacks over   
$(\mathrm{Alg}_{\mathbb{S}})$. The first ingredient is introduced in Subsection   
\ref{ettop} as a natural generalization of the conditions defining \'etale coverings   
in Algebraic Geometry; the second ingredient is contained in the general theory   
developed in Section \ref{stacksonpseudo}. We also study some basic properties of this   
\'etale $K$-theory and suggest some further lines of investigation.\\  
  
A remark on the choice of our setting for \textit{commutative ring spectra} is in order.   
Although we choosed to build everything in this Section starting from  \cite{ekmm}'s   
category $\mathcal{M}_{\mathbb{S}}$ of $\mathbb{S}$-modules, completely analogous   
constructions and results continue to hold if one replaces  
from the very beginning $\mathcal{M}_{\mathbb{S}}$ with any other model for    
spectra having a well behaved smash product. Therefore, the reader could replace   
$\mathcal{M}_{\mathbb{S}}$ with Hovey-Shipley-Smith's category $\mathbf{Sp}^{\Sigma}$   
of symmetric spectra (see \cite{hss}) or with Lydakis' category $\mathbf{SF}$ of   
simplicial functors (see \cite{ly}), with no essential changes.  
  
Moreover, one could also apply the constructions we give below for commutative   
$\mathbb{S}$-algebras, to the category of $E_{\infty}$-algebras over any symmetric   
monoidal model category of the type considered by Markus Spitzweck in \cite[\S 8, 9]{sp}.   
In particular, one can repeat with almost no changes what is in this Section starting from Spitzweck's generalization of $\mathbb{S}$-modules as presented in \cite[\S 9]{sp}.    
  
The problem of defining an \'{e}tale $K$-theory of ring spectra was suggested to us   
by Paul-Arne Ostv\ae r and what we give below is a possible answer to his question. We were very delighted by the question since it looks as a particularly good test of applicability of our theory. For other applications
of the theory developed in this paper to moduli spaces in algebraic topology we refer the reader to \cite{newton}.   
  
\begin{subsection}{$\mathbb{S}$-modules, $\mathbb{S}$-algebras and their algebraic $K$-theory}  
  
The basic reference for what follows is \cite{ekmm}. We fix two universes   
$\mathbb{U}$ and $\mathbb{V}$ with $\mathbb{U}\in \mathbb{V}$.   
These universes are, as everywhere else in this paper, to be understood   
in the sense of \cite[Exp. I, Appendice]{sga4} and \textit{not} in the sense of \cite[1.1]{ekmm}.  
  
\begin{df}   
\begin{itemize}  
  
\item We will denote by $\mathcal{M}_{\mathbb{S}}$ the category of  $\mathbb{S}$-modules in   
the sense of \emph{\cite[II, Def. 1.1]{ekmm}} which belong to $\mathbb{U}$.  
  
\item $\mathrm{Alg}_{\mathbb{S}}$ will denote the category of commutative  $\mathbb{S}$-algebras   
in $\mathbb{U}$, i.e. the category of commutative monoids in $\mathcal{M}_{\mathbb{S}}$.   
Its opposite category will be denoted by $\mathrm{Aff}_{\mathbb{S}}$. Following the standard   
usage in algebraic geometry, an object $A$ in $\mathrm{Alg}_{\mathbb{S}}$, will be formally   
denoted by $\mathrm{Spec}A$ when considered as an object in $\mathrm{Aff}_{\mathbb{S}}$.  
  
\item If $A$ is a commutative $\mathbb{S}$-algebra, $\mathcal{M}_{A}$ will denote the   
category of $A$-modules belonging to $\mathbb{U}$ and $\mathrm{Alg}_{A}$ the category   
of commutative $A$-algebras belonging to $\mathbb{U}$ (i.e. the comma category   
$A/\mathrm{Alg}_{\mathbb{S}}$ of objects in $\mathrm{Alg}_{\mathbb{S}}$ under $A$ or   
equivalently the category of commutative monoids in $\mathcal{M}_{A}$).  
\item We denote by $\mathrm{Alg}_{\mathrm{conn},\;\mathbb{S}}$ the full subcategory of   
$\mathrm{Alg}_{\mathbb{S}}$ consisting of \emph{connective algebras}; its opposite   
category will be denoted by $\mathrm{Aff}_{\mathrm{conn},\;\mathbb{S}}$. If $A$ is a   
(connective) algebra, we denote by $\mathrm{Alg}_{\mathrm{conn},\;A}$ the full subcategory of   
$\mathrm{Alg}_{A}$ consisting of connective $A$-algebras; its   
opposite category will be denoted by $\mathrm{Aff}_{\mathrm{conn},\;A}$.  
\end{itemize}    
\end{df}  
  
Recall that $\mathcal{M}_{A}$ is a topologically enriched, tensored and cotensored   
over the category ($\mathrm{Top}$) of topological spaces in $\mathbb{U}$, left proper  
$\mathbb{U}$-cofibrantly generated $\mathbb{V}$-small model category where   
equivalences are morphisms inducing equivalences on the underlying spectra (i.e.   
equivalences are created by the forgetful functor $\mathcal{M}_{A} \rightarrow   
\mathcal{S}$, where $\mathcal{S}$ denotes the category of spectra \cite[I and VII,   
Th. 4.6]{ekmm} belonging to $\mathbb{U}$) and cofibrations are retracts of relative   
cell $A$-modules (\cite[III, Def. 2.1 (i), (ii); VII, Th. 4.15 ]{ekmm}). Note that   
since the realization functor $\left| - \right|: SSet \rightarrow \mathrm{Top}$ is   
monoidal, we can also view $\mathcal{M}_{\mathbb{S}}$ and $\mathcal{M}_{A}$ as   
tensored and cotensored over $SSet$.  
  
Moreover, a crucial property of $\mathcal{M}_{\mathbb{S}}$ and $\mathcal{M}_{A}$,   
for any commutative $\mathbb{S}$-algebra $A$, is that they admit a refinement of the usual   
``up to homotopy'' smash product of spectra giving them the structure of (topologically   
enriched, tensored and cotensored over the category ($\mathrm{Top}$) of topological spaces   
or over $SSet$) symmetric monoidal model categories (\cite[III, Th. 7.1]{ekmm}).   
  
Finally, both $\mathrm{Alg}_{\mathbb{S}}$ and $\mathrm{Alg}_{A}$ for any commutative   
$\mathbb{S}$-algebra $A$ are topologically or simplicially   
tensored and cotensored model categories (\cite[VII, Cor. 4.10]{ekmm}).  
   
\begin{prop}\label{pseudoconn}  
Let $\iota: \mathrm{Alg}_{\mathrm{conn},\;\mathbb{S}}\hookrightarrow   
\mathrm{Alg}_{\mathbb{S}}$ be the full subcategory of connective algebras and   
$W_{|}$ the set of equivalences in $\mathrm{Alg}_{\mathrm{conn},\;\mathbb{S}}$.   
Then $(\mathrm{Aff}_{\mathrm{conn},\;\mathbb{S}}=(\mathrm{Alg}_{\mathrm{conn},\;\mathbb{S}})^{op},W_{|}^{op},\iota^{op})$   
is a $\mathbb{V}$-small pseudo-model category (see Definition \ref{pseudo}).  
\end{prop}   
\textit{Proof:}   
The only nontrivial property to check is stability of $(\mathrm{Alg}_{\mathrm{conn},\;\mathbb{S}})^{op}$   
under homotopy pullbacks, i.e. stability of $\mathrm{Alg}_{\mathrm{conn},\;\mathbb{S}}$ under homotopy   
push-outs in $\mathrm{Alg}_{\mathbb{S}}$. Let $B\leftarrow A \rightarrow C$ be a diagram in   
$\mathrm{Alg}_{\mathrm{conn},\;\mathbb{S}}$; by \cite[p. 41, after Lemma 9.14]{sp}, there is   
an isomorphism $B\wedge^{\mathbb{L}}_{A}C \simeq B\coprod^{h}_{A}C$ in $\mathrm{\mathrm{Ho}(\mathcal{M}_{A})}$,   
where the left hand side is the derived smash product over $A$ while the right hand side is the   
homotopy pushout in $\mathrm{Alg}_{A}$. Therefore it is enough to know that for any connective   
$A$-modules $M$ and $N$, one has $\pi_{i}(M\wedge^{\mathbb{L}}_{A}N)\equiv \mathrm{Tor}^{A}_{i}(M,N)=0$   
if $i<0$; but this is exactly \cite[Ch. IV, Prop. 1.2 (i)]{ekmm}. \hfill \textbf{$\Box$}    
  
For any commutative $\mathbb{S}$-algebra $A$, the smash   
product $-\wedge_{A} -$ on $\mathcal{M}_{A}$ induces (by derivation) on the homotopy category  
$\mathrm{Ho}(\mathcal{M}_{A})$ the structure of a closed symmetric monoidal category (\cite[III, Th. 7.1]{ekmm}).  
One can therefore define the notion of \textit{strongly dualizable objects}  
in $\mathrm{Ho}(\mathcal{M}_{A})$ (as in \textrm{\cite[\S III.7, ($7.8$)]{ekmm}}).  
The full subcategory of the category $\mathcal{M}_{A}^{c}$ of cofibrant objects in $\mathcal{M}_{A}$, consisting of  
strongly dualizable objects will be denoted by  
$\mathcal{M}_{A}^{\mathrm{sd}}$, and will be endowed with the induced  
classes of cofibrations and equivalences coming from  
$\mathcal{M}_{A}$. It is not difficult to check that with this  
structure, $\mathcal{M}_{A}^{sd}$ is then a Waldhausen  
category (see \textrm{\cite[\S VI]{ekmm}}). Furthermore, if $A  
\longrightarrow B$ is a morphism of commutative $\mathbb{S}$-algebras, then the base change functor  
$$f^{*}:=B \wedge_{A}(-) : \mathcal{M}_{A}^{sd} \longrightarrow \mathcal{M}_{B}^{sd},$$  
being the restriction of a left Quillen functor, preserves equivalences and cofibrations.  
This makes the lax functor  
$$\begin{array}{cccc}  
\mathcal{M}_{-}^{sd} : & \mathrm{Aff}_{\mathbb{S}} & \longrightarrow & \mathrm{Cat}_{\mathbb{V}} \\  
& \mbox{Spec}\, A & \mapsto & \mathcal{M}_{A}^{sd} \\  
& (\mbox{Spec}f : \mbox{Spec}B \rightarrow \mbox{Spec}A) & \mapsto & f^{*}  
\end{array}$$  
into a lax presheaf of Waldhausen $\mathbb{V}$-small categories.  
Applying standard strictification techniques  
(e.g. \cite[Th. 3.4]{may}) and then taking the simplicial set   
(denoted by $\left|wS_{\bullet}\mathcal{M}_{A}^{sd}\right|$ in \cite{wa})   
whose $\Omega$-spectrum is the Waldhausen $K$-theory space,   
we deduce a presheaf of $\mathbb{V}$-simplicial sets  
of $K$-theory  
$$\begin{array}{cccc}  
K(-) : & \mathrm{Aff}_{\mathbb{S}} & \longrightarrow & SSet_{\mathbb{V}} \\  
& \mbox{Spec}\, A & \mapsto & K(\mathcal{M}_{A}^{sd}).  
\end{array}$$  
\noindent The restriction of the simplicial presheaf $K$ to the full subcategory   
$\mathrm{Aff}^{\mathrm{conn}}_{\mathbb{S}}$ of \textit{connective} affine objects will be denoted by   
$$K_{|}(-):\mathrm{Aff}^{\mathrm{conn}}_{\mathbb{S}} \longrightarrow SSet_{\mathbb{V}}.$$  
    
Following Subsection 4.1, we denote by $\mathrm{Aff}_{\mathbb{S}}^{\wedge}$   
(resp. by $\mathrm{Aff}^{\mathrm{conn}}_{\mathbb{S}})^{\wedge}$)  
the model category of pre-stacks over the $\mathbb{V}$-small   
pseudo-model categories $\mathrm{Aff}_{\mathbb{S}}$ (resp. $\mathrm{Aff}^{\mathrm{conn}}_{\mathbb{S}}$).  
  
\begin{df}\label{d34}  
The presheaf $K$ (respectively, the presheaf $K_{|}$) will be considered as an object in   
\mbox{$\mathrm{Aff}_{\mathbb{S}}^{\wedge}$} (resp. in \mbox{($\mathrm{Aff}^{\mathrm{conn}}_{\mathbb{S}})^{\wedge}$})   
and will be  
called the \emph{presheaf of algebraic} $K$-\emph{theory over the symmetric monoidal model category}   
$\mathcal{M}_{\mathbb{S}}$ (resp. the \emph{restricted presheaf of algebraic} $K$-\emph{theory over   
the category} \mbox{$\mathcal{M}^{\mathrm{conn}}_{\mathbb{S}}$} of connective $\mathbb{S}$-modules).  
For any $\mbox{Spec}\, A \in \mathrm{Aff}_{\mathbb{S}}$, we will write  
$$\mathbb{K}(A):=K(\mbox{Spec}\, A).$$  
\end{df}  
  
\begin{rmk} \label{differentK}  
\emph{  
\begin{enumerate}  
\item Note that we adopted here a slightly different definition of the algebraic   
$K$-theory space $\mathbb{K}(A)$ as compared to \cite[VI, Def. 3.2]{ekmm}.   
In fact our Waldhausen category $\mathcal{M}_{A}^{sd}$ (of strongly dualizable objects)   
contains \cite{ekmm} category $f\mathcal{C}_{A}$ of finite cell $A$-modules (\cite[III, Def. 2.1]{ekmm})   
as a full subcategory; this follows from \cite[III, Th. 7.9]{ekmm}. The Waldhausen structure   
on $f\mathcal{C}_{A}$ (\cite[VI, \S 3]{ekmm}) is however different from the one induced   
(via the just mentioned fully faithful embedding) by the Waldhausen structure we use on  
$\mathcal{M}_{A}^{sd}$: the cofibrations in $f\mathcal{C}_{A}$ are fewer. However, the   
same arguments used in \cite[p. 113]{ekmm} after Proposition 3.5, shows that the two   
definitions give isomorphic $K_{i}$ groups for $i>0$ while not, in general, for $i=0$.   
One should think of objects in $f\mathcal{C}_{A}$ as \textit{free modules} while objects in   
$\mathcal{M}_{A}^{sd}$ should be considered as \textit{projective modules}.  
\item Given any commutative $\mathbb{S}$-algebra $A$, instead of considering the simplicial   
set $\mathbb{K}(A)=\left|wS_{\bullet}\mathcal{M}_{A}^{sd}\right|$ whose $\Omega$-spectrum is   
the Waldhausen $K$-theory spectrum of the Waldhausen category $\mathcal{M}_{A}^{sd}$, we could   
as well have taken this spectrum itself and have defined a \textit{spectra}, or better an   
$\mathbb{S}$\textit{-modules valued presheaf} on $\mathrm{Aff}_{\mathbb{S}}$. Since   
$\mathbb{S}$-modules forms a nice simplicial model category, a careful inspection shows that   
all the constructions we made in the previous section still make sense if we replace from the   
very beginning the model category of simplicial presheaves (i.e. of contravariant functors from   
the source pseudo-model category  to simplicial sets in $\mathbb{V}$) with the model category   
of $\mathcal{M}_{\mathbb{S}}$-valued presheaves (i.e. of contravariant functors from the source  
pseudo-model category  to the simplicial model category of $\mathbb{S}$-modules). This leads   
naturally to a theory of \textit{prestacks} or, given a topology on the source pseudo-model or   
simplicial category, to a theory of \textit{stacks in $\mathbb{S}$-modules} (or in any other   
equivalent good category of spectra).  
\item The objects $K$ and $K_{|}$ are in fact underlying simplicial presheaves of   
presheaves of ring spectra, which encodes the ring structure on the $K$-theory spaces.  
We leave to the reader the details of this construction.  
\item A similar construction as the one given above, also yields a $K$-theory presheaf on the  
category of $E_{\infty}$-algebras in a general symmetric monoidal model category $\mathcal{M}$.   
It could be interesting to investigate further the output of this construction when $\mathcal{M}$   
is one of the \textit{motivic} categories considered in \cite[14.8]{sp}.  
\end{enumerate}  
}  
\end{rmk}  
  
\begin{df}\label{d35}  
Let $\tau$ (resp. $\tau'$) be a model pretopology on the model category \mbox{$\mathrm{Aff}_{\mathbb{S}}$}   
(resp. on the pseudo-model category $\mathrm{Aff}^{\mathrm{conn}}_{\mathbb{S}}$), as in Def. \ref{d21}, and let   
\mbox{$\mathrm{Aff}_{\mathbb{S}}^{\sim,\tau}$} (resp. \mbox{$(\mathrm{Aff}^{\mathrm{conn}}_{\mathbb{S}})^{\sim,\tau'}$})   
the associated model category of stacks (Thm. \ref{t5}).  
Let $K \longrightarrow K_{\tau}$ (resp. $K_{|} \longrightarrow K_{|\tau'}$)   
be a fibrant replacement of $K$ (resp. of $K_{|}$) in  
\mbox{$\mathrm{Aff}_{\mathbb{S}}^{\sim,\tau}$} (resp. in \mbox{  
$(\mathrm{Aff}^{\mathrm{conn}}_{\mathbb{S}})^{\sim,\tau'}$}).  
The $K_{\tau}$-\emph{theory space} of a commutative $\mathbb{S}$-algebra   
$A$ (resp. the restricted $K_{\tau'}$-\emph{theory space}   
of a commutative connective  $\mathbb{S}$-algebra $A$) is defined as  
$\mathbb{K}_{\tau}(A):=K_{\tau}(\mbox{Spec}\, A)$   
(resp. as $\mathbb{K}_{|\tau'}(A):=K_{|\tau'}(\mbox{Spec}\, A)$).  
The natural morphism $K \longrightarrow K_{\tau}$ (resp. $K_{|}\longrightarrow   
K_{|\tau'}$) induces a natural augmentation (localization morphism)  
$\mathbb{K}(A) \longrightarrow \mathbb{K}_{\tau}(A)$   
(resp. $\mathbb{K}_{|}(A) \longrightarrow \mathbb{K}_{|\tau'}(A)$).  
   
\end{df}  
  
\begin{rmk}  
\emph{Though we will not give all the details here, one can define also an algebraic   
$K$-theory and $K_{\tau}$-theory space of \textit{any stack} $X \in   
\mbox{$\mathrm{Aff}_{\mathbb{S}}^{\sim,\tau}$}$. The only new ingredient with respect   
to the above definitions is the notion of \textit{1-Segal stack} $\mathrm{Perf_{X}}$   
\textit{of perfect modules over} $X$, that replaces $\mathcal{M}_{A}^{sd}$ in the   
definition above. This notion is defined and studied in the forthcoming paper \cite{partII}.   
Of course, a similar construction is also available for the restricted $K$-theory.}  
  
\end{rmk}  
  
\end{subsection}  
  
\begin{subsection}{The \'{e}tale topology on commutative $\mathbb{S}$-algebras} \label{ettop}  
  
In this section we define an analog of the \'{e}tale topology in the   
category of commutative $\mathbb{S}$-algebras, by extending homotopically   
to these objects the notions of \textit{formally \'{e}tale} morphism and of   
morphism \textit{of finite presentation}.  
  
The notion of formally \'etale morphisms we will use has   
been previously considered by John Rognes \cite{ro} and by Randy McCarthy and Vahagn Minasian \cite{min}, \cite{min2}.  
  
We start with the following straightforward homotopical variation of   
the algebraic notion of finitely presented morphism between commutative   
rings (compare to \cite[Ch. 0, Prop. 6.3.11]{ega0}).  
  
\begin{df} \label{fp}  
A morphism $f:A \rightarrow B$ in $\mathrm{Ho}(\mathrm{Alg}_{\mathbb{S}})$   
will be said to be \emph{of finite presentation} if   
for any filtered direct diagram $C:J \rightarrow \mathrm{Alg}_{A}$, the natural map   
  
$$  
\mathrm{hocolim}_{j \in J}\mathrm{Map}_{\mathrm{Alg}_{A}}(B,C_{j})   
\longrightarrow \mathrm{Map}_{\mathrm{Alg}_{A}}(B, \mathrm{hocolim}_{j \in J}C_{j})  
$$  
  
\noindent is an equivalence of simplicial sets. Here $\mathrm{Map}_{\mathrm{Alg}_{A}}(-,-)$   
denotes the mapping space in the model category $\mathrm{Alg}_{A}$.    
\end{df}  
  
\begin{rmk}  
\emph{  
\begin{enumerate}  
\item It is immediate to check that the condition for   
$\mathrm{Map}_{\mathrm{Alg}_{A}}(-,-)$ of commuting (up to equivalences) with   
$hocolim$ is invariant under equivalences. Hence the definition of finitely   
presented is well posed for a map in the homotopy category   
$\mathrm{Ho}(\mathrm{Alg}_{\mathbb{S}})$.  
\item Since any commutative $A$-algebra can be written as a colimit of   
finite CW $A$-algebras, it is not difficult to show that $A \rightarrow B$   
is of finite presentation if and only if $B$ is a retract of a finite CW $A$-algebra.   
However, we will not use this characterization in the rest of this section.  
\end{enumerate}  
}  
\end{rmk}    
  
We refer to \cite{ba} for the definition and basic properties of topological   
Andr\'e-Quillen cohomology of commutative $\mathbb{S}$-algebras. Recall   
(\cite[Def. 4.1]{ba}) that if $A \rightarrow B$ is a map of commutative   
$\mathbb{S}$-algebras, and $M$ a $B$-module,   
the \textit{topological Andr\'e-Quillen cohomology}   
of $B$ relative to $A$ with coefficient in $M$ is defined as   
$$\mathrm{TAQ}^{*}(B|A,M):=\pi_{-*}F_{B}(\Omega_{B|A},M)=Ext_{B}^{*}(\Omega_{B|A},M),$$  
where $\Omega_{B|A}:=\mathbb{L}Q\mathbb{R}I(B\wedge^{\mathbb{L}}_{A}B)$, $Q$ being the \textit{module of indecomposables}   
functor (\cite[\S 3]{ba}) and $I$ the \textit{augmentation ideal} functor (\cite[\S 2]{ba}). We call $\Omega_{B|A}$ the \textit{topological cotangent complex} of $B$ over $A$. 
In complete analogy to the (discrete) algebraic situation where a morphism of commutative rings is formally   
\'etale if the cotangent complex is homologically trivial (or equivalently   
has vanishing Andr\'e-Quillen cohomology), we give the following (compare, on the algebro-geometric side, with \cite[Ch. III, Prop. 3.1.1]{ill})  
  
\begin{df}\label{etdf}  
\begin{itemize}  
\item A morphism $f:A \rightarrow B$ in $\mathrm{Ho}(\mathrm{Alg}_{\mathbb{S}})$ will be said to be   
\emph{formally \'{e}tale} if the associated topological cotangent complex $\Omega_{B|A}$ is weakly contractible.  
\item A morphism $f:A\rightarrow B$ $\mathrm{Ho}(\mathrm{Alg}_{\mathbb{S}})$ is \emph{\'etale} if it is of finite   
presentation and formally  \'etale. A morphism $\mathrm{Spec}\,B \rightarrow \mathrm{Spec}\,A$ in    
$\mathrm{Ho}(\mathrm{Aff}_{\mathbb{S}})$ is \emph{\'etale} if the map $A\rightarrow B$ in    
$\mathrm{Ho}(\mathrm{Alg}_{\mathbb{S}})$ inducing it, is \'etale.  
\end{itemize}  
\end{df}

\begin{rmk}\label{mandell}  
\emph{  
\begin{enumerate}  
\item Note that if $A' \rightarrow B'$ and $A'' \rightarrow B''$ are morphisms in   
$\mathrm{Alg}_{\mathbb{S}}$, projecting to isomorphic maps in $\mathrm{Ho}(\mathrm{Alg}_{\mathbb{S}})$,   
then $\Omega_{B'|A'}$ and $\Omega_{B''|A''}$ are isomorphic in the homotopy category of $\mathbb{S}$-modules. Therefore, the condition given above of   
being formally \'etale is well defined for a map in $\mathrm{Ho}(\mathrm{Alg}_{\mathbb{S}})$.  
\item \textit{$\mathrm{THH}$-\'etale morphisms.} If $A$ is a commutative $\mathbb{S}$-algebra,   
$B$ a commutative $A$-algebra, we recall that $\mathrm{Alg}_{A}$ is tensored and cotensored   
over $\mathrm{Top}$ or equivalently over $SSet$; therefore it makes sense to consider the object   
$S^{1}\otimes^{\mathbb{L}}B$ in $\mathrm{Ho}(\mathrm{Alg}_{A})$, where the derived tensor product   
is performed in $\mathrm{Alg}_{A}$. By a result of McClure, Schw\"anzl and Vogt (see \cite[IX, Th. 3.3]{ekmm}),   
$S^{1}\otimes^{\mathbb{L}}B$ is isomorphic to $\mathrm{THH}^{A}(B;B)\equiv \mathrm{THH}(B|A)$ in $\mathrm{Ho}(\mathrm{Alg}_{A})$   
and is therefore a model for \emph{topological Hochschild homology} as defined e.g in \cite[IX.1]{ekmm}. Moreover,   
note that any choice of a point $* \rightarrow S^{1}$ gives to $S^{1}\otimes^{\mathbb{L}}B$ a canonical   
structure of $A$-algebra.\\ A map $A \rightarrow B$ of commutative $\mathbb{S}$-algebras, will be   
called \textit{formally $\mathrm{THH}$-\'etale} if the  
canonical map $B \rightarrow S^{1}\otimes^{\mathbb{L}}B$ is an isomorphism in $\mathrm{Ho}(\mathrm{Alg}_{A})$;   
consequently, a map $A \rightarrow B$ of commutative $\mathbb{S}$-algebras, will be called \textit{$\mathrm{THH}$-\'etale}   
if it is formally $\mathrm{THH}$-\'etale and of finite presentation.   
As shown by Vahagn Minasian (\cite{min}) $\mathrm{THH}$-\'etale morphisms are in particular \'etale.  
\item It is easy to see that a morphism of commutative $\mathbb{S}$-algebras $A \rightarrow B$ is formally $\mathrm{THH}$-\'etale   
if and only if $B$ is a \textit{co-discrete} object in the model category $\mathrm{Alg}_{A}$ i.e., if for any   
$C \in \mathrm{Alg}_{A}$ the mapping space $\mathrm{Map}_{\mathrm{Alg}_{A}}(B,C)$ is a discrete (i.e. $0$-truncated)   
simplicial set. From this description, one can produce examples of \'etale morphisms of $\mathbb{S}$-algebras  
which are not $\mathrm{THH}$-\'etale. The following example was communicated to us by Michael Mandell.  
Let $A = H\mathbb{F}_{p}=K(\mathbb{F}_{p},0)$   
($H$ denotes the Eilenberg-Mac Lane $\mathbb{S}$-module functor, see \cite[IV, \S 2]{ekmm}),   
and perform the following construction. Start with $F_{1}(A)$, the free commutative $A$-algebra   
on a cell in degree $-1$. In $\pi_{-1}(F_{1}(A))$ there is a fundamental class   
but also lots of other linearly independent elements as for example the Frobenius $F$.  
We let $B$ to be the $A$-algebra defined by the following homotopy co-cartesian square  
$$\xymatrix{  
F_{1}(A) \ar[r]^-{1-F} \ar[d] & F_{1}A \ar[d] \\  
A \ar[r] & B.}$$  
The morphism $1-F$ being \'etale, we have that $B$ is an \'etale $A$-algebra. However,   
one has $\pi_{1}(Map_{Alg_{A}}(B,A))\simeq \mathbb{Z}/p \neq 0$, and therefore  
$A \longrightarrow B$ is not $\mathrm{THH}$-\'etale (because $Map_{Alg_{A}}(B,A)$ is not $0$-truncated).  
\end{enumerate}  
}   
\end{rmk}  
  
\begin{prop}\label{etalebasechange}  
If $C\leftarrow A\rightarrow B$ is a diagram in   
$\mathrm{Ho}(\mathrm{Alg}_{\mathbb{S}})$ and $A\rightarrow B$ is \'etale,   
then the homotopy co-base change map $C\rightarrow B\coprod_{A}^{h}C$ is again \'etale.  
\end{prop}  
  
\textit{Proof:} The co-base change invariance of the finite presentation property is easy and left to the reader.   
The co-base change invariance of the formally \'etale property follows at   
once from \cite[p. 41, after Lemma 9.14]{sp} and the ``flat base change'' formula for the cotangent complex (\cite[Prop. 4.6]{ba}) 
$$\Omega_{B\wedge^{\mathbb{L}}_{A}C|C}\simeq \Omega_{B|A}\wedge_{A}C.$$  
\hfill \textbf{$\Box$}  
  
As an immediate consequence we get the following corollary.  
  
\begin{cor}\label{smallsitesareok}  
Let $A$ be a commutative $\mathbb{S}$-algebra. The subcategory   
$\mathrm{Aff}^{\textrm{\'{e}t}}_{A}$ of $\mathrm{Aff}_{A}$ consisting of   
\'etale maps $\mathrm{Spec}\,B \rightarrow \mathrm{Spec}\,A$, is a pseudo-model category.  
\end{cor}  
  
For any (discrete) commutative ring $R$, we denote by $HR=\mathrm{K}(R,0)$ the   
\textit{Eilenberg-Mac Lane} commutative $\mathbb{S}$-algebra associated to $R$   
(\cite[IV, \S 2]{ekmm}). 

\begin{prop}\label{stefan}
A morphism of discrete commutative rings $R\rightarrow R'$ is \'etale iff $HR\rightarrow HR'$ is \'etale
\end{prop}
\textit{Proof:} By \cite{pr} and \cite{bm}, we can apply to topological Andr\'e-Quillen homology and Andr\'e-Quillen homology the two spectral sequences at the end of \cite[\S 7.9]{sch} to conclude that the algebraic cotangent complex $\mathbb{L}_{R'/R}$ is acyclic iff the topological cotangent complex $\Omega_{HR'|HR}$ is weakly contractible; therefore the two formal etaleness do imply each other. Also the two finite presentation condition easily imply each other, since the functor $\pi_{0}$ is left adjoint and therefore preserves finitely presented objects. So we only have to observe that a finitely presented morphism of discrete commutative rings $R\longrightarrow R'$ is \'etale iff it has an acyclic algebraic cotangent complex (\cite[Ch. III, Prop. 3.1.1]{ill}).  
\hfill \textbf{$\Box$}\\

The following proposition   
compare the notions of \'etale morphisms of commutative rings and commutative   
$\mathbb{S}$-algebras in the connective case.

\begin{prop}\label{pierre}  
Let $k$ be a commutative ring (in $\mathbb{U}$), and $Hk \longrightarrow B$ be an \'etale morphism of  
connective commutative $\mathbb{S}$-algebras. Then, the natural map $B \longrightarrow H(\pi_{0}(B))$ (\cite[Prop. IV.3.1]{ekmm})  
is an equivalence of commutative $\mathbb{S}$-algebras. Therefore, up to equivalences, $Hk \longrightarrow B$ is of the form $Hk \longrightarrow Hk'$ where $k\rightarrow k'$ is an \'etale extension of discrete commutative rings.  
\end{prop}  

\textit{Proof:} Consider the sequence of maps of commutative $\mathbb{S}$-algebras $Hk\longrightarrow B \longrightarrow H\pi_{0}(B)$; this gives a fundamental cofibration sequence (\cite[Prop. 4.3]{ba}) $$\Omega_{B|Hk}\wedge_{B}H\pi_{0}(B)\longrightarrow \Omega_{H\pi_{0}(B)|Hk}\longrightarrow \Omega_{H\pi_{0}(B)|B}.$$
Since $Hk\longrightarrow B$ is \'etale, by \cite[Prop. 3.8 (2)]{min2} also $Hk\longrightarrow H\pi_{0}(B)$ is \'etale; therefore the first two terms are contractible, hence  $\Omega_{H\pi_{0}(B)|B}\simeq *$, too. 
Now, the map $B\longrightarrow H\pi_{0}(B)$ is a $1$-equivalence (see also \cite[Proof of Thm. 8.1]{ba}) and therefore, $\Omega_{H\pi_{0}(B)|B}\simeq *$ and  \cite[Lemma 8.2]{ba}, tell us that $\pi_{1}B\simeq 0$. Then, $B\longrightarrow H\pi_{0}(B)$ is also a $2$-equivalence and the same argument shows then that $\pi_{2}B\simeq 0$, etc. Therefore $\pi_{i}B\simeq 0$, for any $i\geq 1$ and we get the first statement. The second one follows from this and Proposition \ref{stefan}. \hfill \textbf{$\Box$}

\begin{rmk}  
\emph{Note that Proposition \ref{pierre} is false if we remove the connectivity hypothesis.   
In fact, the $H\mathbb{F}_{p}$-algebra $B$ described in Remark \ref{mandell} (3) is \'etale   
but has, by construction, non-vanishing homotopy groups in  
infinitely many negative degrees. Actually, even restricting to
$thh$-etale characteristic 
zero will not be enough in order to avoid this 
kind of phenomenon (see e.g. \cite[Rem. 2.19]{newton}).}  
\end{rmk}    
  
\begin{df}\label{etalecoverings}  
For each $\mbox{Spec}A \in \mathrm{Ho}(\mathrm{Aff}_{\mathbb{S}})$,   
let us define $\mathrm{Cov}_{\textrm{\'{e}t}}(\mbox{Spec}A)$ as the set of \emph{finite} families $\{f_{i}  
: \mbox{Spec}\, B_{i} \longrightarrow \mbox{Spec}\, A\}_{i \in I}$ of morphisms    
in $\mathrm{Ho}(\mathrm{Aff}_{\mathbb{S}})$, satisfying the following two conditions:  
  
\begin{enumerate}  
\item  for any $i \in I$, the morphism $A \longrightarrow B_{i}$ is \emph{\'{e}tale};  
\item the family of base change functors  
  
$$\{\mathbb{L}f_{i}^{*} : \mathrm{Ho}(\mathcal{M}_{A}) \longrightarrow \mathrm{Ho}(\mathcal{M}_{B_{i}})\}_{i \in I}$$  
  
conservative, i.e. a morphism in $\mathrm{Ho}(\mathcal{M}_{A})$ is an isomorphism if and only if,  
for any $i \in I$, its image in $\mathrm{Ho}(\mathcal{M}_{B_{i}})$ is an isomorphism.   
  
\end{enumerate}  
  
\end{df}  
  
We leave to the reader the easy task of checking that this actually defines a model   
pre-topology (\textit{\'{e}t}) (see Def. \ref{d21}), called the \textit{\'{e}tale topology}   
on \mbox{$\mathrm{Aff}_{\mathbb{S}}$}. By restriction to the sub-pseudo-model category   
(see Prop. \ref{pseudoconn}) \mbox{$\mathrm{Aff}_{\mathrm{conn},\;\mathbb{S}}$} of   
connective objects, we also get a pseudo-model site $(\mbox{$\mathrm{Aff}_{\mathrm{conn},\;\mathbb{S}}$},\textrm{\'et})$,   
called the \textit{restricted \'etale site}.  
   
If $A$ is a commutative (resp. commutative and connective) $\mathbb{S}$-algebra,   
the pseudo-model category (see Cor. \ref{smallsitesareok}) \mbox{$\mathrm{Aff}_{\textrm{\'{e}t}/A}$}   
(resp. \mbox{$\mathrm{Aff}_{\textrm{conn,\;\'{e}t}/A}$}), together with the ``restriction'' of   
the \'etale topology, will be called the \textit{small \'etale site}   
(resp. the \textit{restricted small \'etale site}) over $A$.   
More precisely, let us consider the obvious forgetful functors  
$$F:\mbox{$\mathrm{Aff}_{\textrm{\'{e}t}/A}$} \longrightarrow \mbox{$\mathrm{Aff}_{\mathbb{S}}$}$$  
$$F':\mbox{$\mathrm{Aff}_{\textrm{conn,\;\'{e}t}/A}$}\longrightarrow   
\mbox{$\mathrm{Aff}_{\mathbb{S}}$}.$$  
By definition of the pseudo-model structures on \mbox{$\mathrm{Aff}_{\textrm{\'{e}t}/A}$}   
(resp. on \mbox{$\mathrm{Aff}_{\textrm{conn,\;\'{e}t}/A}$}), $F$ (resp. $F'$) preserves   
(actually, creates) equivalences. Therefore, we say that family of morphisms  
$\left\{\mathrm{Spec}(C_{i})\rightarrow \mathrm{Spec}(B) \right\}$ in   
$\mathrm{Ho}(\mbox{$\mathrm{Aff}_{\textrm{\'{e}t}/A}$})$ (resp. in   
$\mathrm{Ho}(\mbox{$\mathrm{Aff}_{\textrm{conn,\;\'{e}t}/A}$})$) is an   
\'etale covering family of $(\mathrm{Spec}B\rightarrow \mathrm{Spec}A)$ in   
\mbox{$\mathrm{Aff}_{\textrm{\'{e}t}/A}$} (resp. \mbox{$\mathrm{Aff}_{\textrm{conn,\;\'{e}t}/A}$})   
iff its image via $\mathrm{Ho}(F)$ (resp. via $\mathrm{Ho}(F')$)   
is an \'etale covering family of $\mathrm{Spec}(B)$ in  
\mbox{$\mathrm{Aff}_{\mathbb{S}}$} i.e. belongs to  $\mathrm{Cov}_{\textrm{\'{e}t}}(\mbox{Spec}A)$   
(Def. \ref{etalecoverings}). \\  
  
We finish this paragraph by the following corollary that compare the  
small \'etale sites of a ring $k$ and of its associated Eilenberg-Mac Lane $\mathbb{S}$-algebra $Hk$.  
  
\begin{cor}\label{sesae}  
Let $k$ be a discrete commutative ring, $(\mathrm{aff_{\textrm{\'et}/k}},\textrm{\'et})$   
be the small \'etale affine site over $\mathrm{Spec}(k)$ consisting of affine \'etale   
schemes $\mathrm{Spec}(k')\rightarrow \mathrm{Spec}(k)$, and   
$H : \mathrm{aff_{\textrm{\'et}/k}}\longrightarrow \mbox{$\mathrm{Aff}_{\textrm{conn,\;\'{e}t}/Hk}$}$   
be the Eilenberg-Mac Lane space functor.  
Then $H$ induces a continuous equivalence of \'etale pseudo-model sites $$H:(\mathrm{aff_{\textrm{\'et}/k}},  
\textrm{\'et})\rightarrow (\mbox{$\mathrm{Aff}_{\textrm{conn,\;\'{e}t}/A}$},\textrm{\'et}).$$  
\end{cor}  
  
\textit{Proof:} Propositions \ref{pierre} and \ref{stefan} imply that the conditions of Proposition \ref{critconteq} are satisfied. \hfill \textbf{$\Box$}  

\end{subsection}  
  
\begin{subsection}{\'{E}tale $K$-theory of commutative $\mathbb{S}$-algebras}  
  
The following one is the main definition of this section.  
  
\begin{df}\label{dket}  
\begin{itemize}  
\item For any $A \in \mbox{$\mathrm{Alg}_{\mathbb{S}}$}$, we define   
its \textit{\'{e}tale} $K$-\textit{theory space} $\mathbb{K}_{\textrm{\'{e}t}}(A)$   
by applying Definition \ref{d35} to $\tau =\emph{(\'{e}t})$.  
\item For any $A \in \mbox{$\mathrm{Alg}^{\mathrm{conn}}_{\mathbb{S}}$}$,   
we define its \textit{restricted \'{e}tale} $K$-\textit{theory space}   
$\mathbb{K}_{|\textrm{\'{e}t}}(A)$ by applying Definition \ref{d35} to $\tau' =\emph{(\'{e}t})$.  
\end{itemize}  
\end{df}  
    
The following Proposition shows that, as in the algebraic case (compare   
\cite[Thm. 3.10]{ja}), also in our context, \'etale $K$-theory can be computed on the small \'etale sites.  
  
\begin{prop}\label{smallsitescompute}  
Let $A$ be a commutative (resp. commutative and connective) \   
$\mathbb{S}$-algebra and  $(\mbox{$\mathrm{Aff}_{\textrm{\'{e}t}/A}$})^{\sim,\;{\textrm{\'{e}t}}}$   
(resp. $(\mbox{$\mathrm{Aff}_{\textrm{conn,\;\'{e}t}/A}$})^{\sim,\;{\textrm{\'{e}t}}}$) the   
model category of stacks on the small \'etale site (resp. on the restricted small \'etale site)   
over $A$. For any presheaf $F$ on $\mathrm{Aff}_{\mathbb{S}}$, we denote by $F^{sm}$ (resp.   
$F_{|}^{sm}$) its restriction to $\mbox{$\mathrm{Aff}_{\textrm{\'{e}t}/A}$}$ (resp. to   
$\mbox{$\mathrm{Aff}_{\textrm{conn,\;\'{e}t}/A}$}$). Then the map $K^{sm} \rightarrow   
K_{\textrm{\'{e}t}}^{sm}$ (resp. $K_{|}^{sm} \rightarrow K_{\textrm{\'{e}t}\,|}^{sm}$)   
induced via restriction by a fibrant replacement $K \rightarrow K_{\textrm{\'{e}t}}$   
(resp. $K_{|} \rightarrow K_{|\textrm{\'{e}t}}$) in \mbox{$(\mathrm{Aff}_{\mathbb{S}})^{\sim,\,\textrm{\'{e}t}}$}   
(resp. in \mbox{$(\mathrm{Aff}_{\mathrm{conn},\!\mathbb{S}})^{\sim,\,\textrm{\'{e}t}}$}), is a   
fibrant replacement in $(\mbox{$\mathrm{Aff}_{\textrm{\'{e}t}/A}$})^{\sim,\;{\textrm{\'{e}t}}}$   
(resp. in $(\mbox{$\mathrm{Aff}_{\textrm{conn,\;\'{e}t}/A}$})^{\sim,\;{\textrm{\'{e}t}}}$).  
\end{prop}  
  
\textit{Proof:} We prove the proposition in the non-connective case, the connective case is the similar.  
  
\noindent Let us consider the natural functor  
$$f : Aff_{\textrm{\'et}/A} \longrightarrow Aff_{\mathbb{S}},$$  
from the small \'etale site of $Spec\, A$ to the big \'etale site.   
It is clear that the associated restriction functor  
$$f^{*} : Aff_{\mathbb{S}}^{\sim,\textrm{\'et}} \longrightarrow Aff_{\textrm{\'et}/A}^{\sim,\textrm{\'et}}$$  
preserves equivalences (one can apply for example Lemma \ref{l13}).  
Furthermore, if $Spec\, B \longrightarrow Spec\, A$ is  
a fibrant object in $Aff_{\textrm{\'et}/A}$, then the pseudo-representable  
hypercovers (see Definition \ref{dpseudorep})   
of $Spec\, B$ are the same in $Aff_{\textrm{\'et}/A}$ and in $Aff_{\mathbb{S}/A}$   
(because each structure map of a pseudo-representable hypercover is \'etale). This implies by  
Corollary \ref{c62}, that the functor $f^{*}$ preserves fibrant objects. In particular, if  
$K \longrightarrow K_{\textrm{\'et}}$ is a fibrant replacement in $Aff_{\mathbb{S}}^{\sim,\textrm{\'et}}$, so  
is its restriction to $Aff_{\textrm{\'et}/A}^{\sim,\textrm{\'et}}$. \hfill $\Box$  
  
As a consequence, we get the following comparison result to algebraic   
\'etale $K$-theory for fields; if $R$ is a (discrete) commutative   
ring, we denote by $K_{\textit{\'{e}t}}(R)$ its \'etale $K$-theory space (e.g. \cite{ja}).  
  
\begin{cor}\label{restrictedcomparison}  
For any discrete commutative ring $k$, we have an isomorphism $\mathbb{K}_{|\textrm{\'{e}t}}(Hk)   
\simeq K_{\textit{\'{e}t}}(k)$ in $\mathrm{Ho}(SSet)$.  
\end{cor}  
  
\textit{Proof:} This follows from corollaries \ref{sesae}, \ref{smallsitescompute} and from  
the comparison between algebraic $K$-theory of a commutative ring $R$ and algebraic $K$-theory  
of the $\mathbb{S}$-algebra $HR$ (see \cite[VI, Rmk. $6.1.5 (1)$]{ekmm}).  
\hfill \textbf{$\Box$}  

\end{subsection}  
  
\end{section}  
  
\begin{appendix}
  
\begin{section}{Model categories and universes}  
  
In this appendix we have collected the definitions of $\mathbb{U}$-cofibrantly generated,   
$\mathbb{U}$-cellular and $\mathbb{U}$-combinatorial model categories for a universe  
$\mathbb{U}$, that have been used all along this work.   
  
Throughout this appendix, we fix a universe $\mathbb{U}$.  
  
\begin{subsection}{$\mathbb{U}$-cofibrantly generated model categories}  
  
Recall that a category is a $\mathbb{U}$-category, or equivalently a   
locally $\mathbb{U}$-small category, if for any pair of objects  
$(x,y)$ in $C$ the set $Hom_{C}(x,y)$ is a $\mathbb{U}$-small set.  
  
\begin{df}\label{db1}  
A $\mathbb{U}$-model category is a category $M$ endowed with a model structure  
in the sense of \cite[Def. $1.1.3$]{ho} and satisfying the following two conditions  
\begin{enumerate}  
\item The underlying category of $M$ is  
a $\mathbb{U}$-category.  
\item  
The underlying category of $M$ has all kind of $\mathbb{U}$-small limits and colimits.  
\end{enumerate}  
\end{df}  
  
Let $\alpha$ be the cardinal of a $\mathbb{U}$-small set (we will simply say   
\textit{$\alpha$ is a $\mathbb{U}$-small cardinal}).  
Recall from \cite[Def. $2.1.3$]{ho} that an object $x$ in a $\mathbb{U}$-model category $M$,   
is $\alpha$\textit{-small}, if for any $\mathbb{U}$-small $\alpha$-filtered ordinal $\lambda$, and any  
$\lambda$-sequence  
$$y_{0} \rightarrow y_{1} \rightarrow \dots y_{\beta} \rightarrow y_{\beta+1} \rightarrow \dots$$  
the natural map  
$$\mathrm{colim}_{\beta < \lambda}Hom(x,y_{\beta}) \longrightarrow Hom(x,\mathrm{colim}_{\beta<\lambda}y_{\beta})$$  
is an isomorphism.   
  
We will use (as we did in the main text) the following   
variation of the notion of \textit{cofibrantly generated model category}  
of \cite[Def. $2.1.17$]{ho}.  
  
\begin{df}\label{db2}  
Let $M$ be a $\mathbb{U}$-model category. We say that $M$ is $\mathbb{U}$\emph{-cofibrantly generated}  
if there exist $\mathbb{U}$-small sets $I$ and $J$ of morphisms in $M$, and a $\mathbb{U}$-small cardinal $\alpha$,   
such that the following three conditions are satisfied  
\begin{enumerate}  
\item The domains and codomains of the maps of $I$ and $J$ are $\alpha$-small.  
\item The class of fibrations is $J$-inj.  
\item The class of trivial fibrations is $I$-inj.  
\end{enumerate}  
\end{df}  
  
The main example of a $\mathbb{U}$-cofibrantly generated model category is the model categories  
$SSet_{\mathbb{U}}$ of $\mathbb{U}$-small simplicial sets.  
  
The main ``preservation'' result is the following easy proposition (see \cite[\S $13.8$, $13.9$, $13.10$]{hi}).  
  
\begin{prop}\label{pb1}  
Let $M$ be a $\mathbb{U}$-cofibrantly generated model category.  
\begin{enumerate}  
\item If $C$ is a $\mathbb{U}$-small category, then the category $M^{C}$ of $C$-diagrams in $M$ is   
again a $\mathbb{U}$-cofibrantly generated model category in which equivalences and   
fibrations are defined objectwise.  
  
\item Let us suppose that $M$ is furthermore a $SSet_{\mathbb{U}}$-model category in the sense  
of \cite[Def. $4.2.18$]{ho} (in other words, $M$ is a simplicial $\mathbb{U}$-cofibrantly generated model category),   
and let $T$ be a $\mathbb{U}$-small $S$-category. Then, the category $M^{T}$ of simplicial functors  
from $T$ to $M$ is again a $\mathbb{U}$-cofibrantly generated model category   
in which the equivalences and fibrations are defined objectwise. The model category   
$M^{T}$ is furthermore a $SSet_{\mathbb{U}}$-model category in the sense of \cite[Def. $4.2.18$]{ho}.  
  
\end{enumerate}  
\end{prop}  
  
A standard construction we have been using very often in the main text is the following. We start by the model category   
$SSet_{\mathbb{U}}$ of $\mathbb{U}$-small simplicial sets. Now, if $\mathbb{V}$ is a universe  
with $\mathbb{U} \in \mathbb{V}$, then the category $SSet_{\mathbb{U}}$ is $\mathbb{V}$-small. Therefore,   
the category   
$$SPr(SSet_{\mathbb{U}}):=SSet_{\mathbb{V}}^{SSet_{\mathbb{U}}^{op}}$$  
of $\mathbb{V}$-small simplicial presheaves on $SSet_{\mathbb{U}}$, is a   
$\mathbb{V}$-cofibrantly generated model category.   
  
This is the way we have considered, in the main text, model categories of diagrams over a base model category  
avoiding any set-theoretical problem.  
  
\end{subsection}  
  
\begin{subsection}{$\mathbb{U}$-cellular and $\mathbb{U}$-combinatorial model categories}  

The following notion of   
\textit{combinatorial model category} is due to Jeff Smith (see, for example, \cite[\S 2]{du2}, \cite[I, \S 1]{bk}).  
  
\begin{df}\label{db3}  
\begin{enumerate}  

\item A category $C$ is called $\mathbb{U}$\emph{-locally presentable} (see \cite{du2}) if 
there exists a $\mathbb{U}$-small set of objects $C_0$ in $C$, which are all
$\alpha$-small for some cardinal $\alpha$ in $\mathbb{U}$ and such that  
any object in $C$ is an $\alpha$-filtered colimit of objects in $C_0$.
  
\item A $\mathbb{U}$\emph{-combinatorial model category} is a $\mathbb{U}$-cofibrantly generated model  
category whose underlying category is $\mathbb{U}$-locally presentable.  
  
\end{enumerate}  
\end{df}  
  
The following localization theorem is due to J. Smith (unpublished).   
Recall that a model category is \textit{left proper} if the   
equivalences are closed with respect to pushouts along cofibrations.  
  
\begin{thm}\label{tb1}  
Let $M$ be a left proper, $\mathbb{U}$-combinatorial model category, and $S \subset M$ be a   
$\mathbb{U}$-small subcategory. Then the left Bousfield localization $L_{S}M$  
of $M$ with respect to $S$ exists.  
\end{thm}  
  
Let us recall from \cite[\S $12.7$]{hi} the notion of \textit{compactness}. We will say that an object $x$ in a   
$\mathbb{U}$-cofibrantly generated model category $M$ is \textit{compact} is there exists  
a $\mathbb{U}$-small cardinal $\alpha$ such that $x$ is $\alpha$-compact  
in the sense of \cite[Def. $13.5.1$]{hi}.  
The following definition is our variation of the notion of \textit{cellular model category} of \cite{hi}.  
  
\begin{df}\label{db4}  
A $\mathbb{U}$-cellular model category $M$ is a $\mathbb{U}$-cofibrantly generated model category  
with generating $\mathbb{U}$-small sets of cofibrations $I$ and of trivial cofibrations $J$,  
such that the following two conditions are satisfied  
\begin{enumerate}  
\item The domains and codomains of maps in $I$ are compact.  
  
\item Monomorphisms in $M$ are effective.  
  
\end{enumerate}  
\end{df}  
  
The main localization theorem of \cite{hi} is the following.  
  
\begin{thm}{(P. Hirschhorn, \cite[Thm. 4.1.1]{hi})}\label{tb2}   
Let $M$ be a left proper, $\mathbb{U}$-cellular model category and $S\subset M$ be a   
$\mathbb{U}$-small subcategory. Then the left Bousfield localization $L_{S}M$  
of $M$ with respect to $S$ exists.  
\end{thm}  
  
Finally, let us mention the following ``preservation'' result.  
  
\begin{prop}\label{pb2}  
If in Proposition \ref{pb1}, $M$ is $\mathbb{U}$-combinatorial (resp. $\mathbb{U}$-cellular), then   
so are $M^{C}$ and $M^{T}$.  
\end{prop}  
  
\end{subsection}  
  
\end{section}  
  
\end{appendix}

\end{document}